\documentclass[final,3p,times,11pt]{elsarticle}

\usepackage{amsmath,amssymb,amsthm, mathrsfs}
\usepackage{mathtools}
\usepackage{graphicx}
\usepackage{stmaryrd}
\usepackage[dvipsnames]{xcolor}
\usepackage{cancel}
\usepackage{ulem}
\usepackage{showlabels}
\usepackage[ruled,vlined]{algorithm2e}
\usepackage{hyperref}
\usepackage{marginnote}
\usepackage{multirow}
\hypersetup{
    colorlinks=true,       
}

\newcommand{\RefereeA}[1]{{\textcolor{black}{#1}}}
\newcommand{\RefereeB}[1]{{\textcolor{black}{#1}}}

\journal{}
\makeatletter
\def\ps@pprintTitle{%
 \let\@oddhead\@empty
 \let\@evenhead\@empty
 \def\@oddfoot{}%
 \let\@evenfoot\@oddfoot}
\makeatother

\begin{document}

\begin{frontmatter}


\title{Surrogate Approximation of the Grad--Shafranov Free Boundary Problem via
Stochastic Collocation on Sparse Grids}



\author[umdcs]{Howard C.\ Elman}
\ead{helman@umd.edu}
\address[umdcs]{Department of Computer Science and Institute for Advanced Computer
Studies, University of Maryland, College Park.}
\author[umdm]{Jiaxing Liang}
\ead{jliang18@umd.edu}
\address[umdm]{Applied Mathematics \& Statistics and Scientific Computing Program, University of Maryland, College Park.}
\author[nyu,UA]{Tonatiuh S\'anchez-Vizuet}
\ead{tonatiuh@math.arizona.edu}
\address[nyu]{Courant Institute of Mathematical Sciences. New York University}
\address[UA]{Department of Mathematics, The University of Arizona}
\begin{abstract}
In magnetic confinement fusion devices, the equilibrium configuration of a plasma is determined by the balance between the hydrostatic pressure in the fluid and the magnetic forces generated by an array of external coils and the plasma itself. The location of the plasma is not known \textit{a priori} and must be obtained as the solution to a free boundary problem. The partial differential equation that determines the behavior of the combined magnetic field depends on a set of physical parameters (location of the coils, intensity of the electric currents going through them, magnetic permeability, etc.) that are subject to uncertainty and variability. The confinement region is in turn a function of these stochastic parameters as well. In this work, we consider variations on the current intensities running through the external coils as the dominant source of uncertainty. This leads to a parameter space of dimension equal to the number of coils in the reactor. With the aid of a surrogate function built on a sparse grid in parameter space, a Monte Carlo strategy is used to explore the effect that stochasticity in the parameters has on important features of the plasma boundary such as the location of the $x$-point, the strike points, and shaping attributes such as triangularity and elongation. The use of the surrogate function reduces the time required for the Monte Carlo simulations by factors that range between 7 and over 30.  
\end{abstract}

\end{frontmatter}


\section{Introduction}\label{sec:intro}
In a fusion reactor, a gas of hydrogen isotopes is injected into a vacuum vessel where it is heated to very high temperatures until it ionizes into a plasma. The very high temperatures induce collisions between the ionized nucleii which, if sustained for long enough at a high enough temperature, eventually lead to the merger of nucleii and the release of large amounts of energy through \textit{nuclear fusion}. One of the most successful procedures to force the plasma to remain confined for long enough to allow for a thermonuclear reaction to take place, and to prevent it from damaging the walls of the reactor, is to use magnetic fields. To achieve magnetic confinement, external coils are used to generate a magnetic field aimed at balancing the internal forces that drive the plasma to expand. In this study, we explore the impact that the intensity of the electric currents going through the coils, which we treat as stochastic parameters, has on properties of the confined plasma.

We will focus on equilibrium computations on \textit{tokamaks}, a family of axially symmetric reactors whose design dates back to the work of Soviet scientists in the late 1960's \cite{ArBoGo:1969, Wesson:2011}. 
In the axially symmetric setting, described by cylindrical coordinates  $(r,z,\phi)$, the magnetic field can be expressed in terms of two scalar functions which, following the nomenclature from \cite[page 94]{Jardin2010}, will be referred to as the \textit{poloidal flux function $\psi(r,z)$ and the \textit{toroidal field function}, $g(r,z)$} (which is proportional to the total current flowing in the poloidal direction). \RefereeB{Due to the cylindrical symmetry of the device, the problem is independent from the toroidal angle $\phi$ and it might be considered as if it were posed on a two-dimensional plane. For this reason, in this paper we will drop the use of $(r,z)$ to denote the position on the plane for any fixed angle $\phi$, and will instead use the notation $(x,y)$ to denote the distance from the axis of symmetry and the height, respectively.} 

The equilibrium between magnetic and hydrodynamic pressure in a plasma leads to 
a differential relation (see Section \ref{sec:mag-equal} for details) between the unknown poloidal flux and the hydrodynamic pressure $p$ and current densities, giving the following
free boundary problem:
\begin{subequations}\label{eq:FreeBoundary}
\begin{equation}\label{eq:FreeBoundarya}
-\nabla\,\cdot\,\left(\frac{1}{\mu x}\nabla\psi\right) = \left\{ \begin{array}{ll}
\RefereeA{x\frac{d}{d\psi} p(\psi) + \frac{1}{2\,\mu x} \frac{d}{d\psi} g^2(\psi)} & \text{ in } \Omega_p(\psi) \\
I_i/S_i & \text{ in } \Omega_{C_i} \\
0 & \text{ elsewhere. } 
\end{array}\right.
\end{equation}
The magnetic permeability $\mu$ may be considered to be a constant, although in the presence of an iron component, the permeability becomes a function of the magnitude of the magnetic field within the iron structure,  giving
\[
\mu = \left\{\begin{array}{cr}
\mu_0 & \text{in vacuum}, \\
\mu(|\nabla \psi|^2/x^2) & \text{in iron},
\end{array}\right.
\]
where $\mu_0$ is the magnetic permeability of the vacuum. The specific form of the free functions  $p$ and $g$, as well as the concrete values of the currents going through the coils and the locations of the coils are all user-provided problem data.

The first row of the second order, semilinear, elliptic partial differential equation \eqref{eq:FreeBoundarya} is known as the \textit{Grad--Shafranov equation}. It was derived independently by Grad and Rubin \cite{GrRu:1958}, Shafranov \cite{Shafranov:1958}, and L\"{u}st and Schl\"{u}ter \cite{LuSc:1957}, and describes the magnetic equilibrium in the interior of the plasma. When considering the equilibrium problem without \textit{a priori} knowledge of the location of the plasma, $\Omega_p$, and the shape of its boundary, $\partial \Omega_p$, the source term must be augmented as in \eqref{eq:FreeBoundarya} to include the contributions due to the internal forces on the plasma (in the confinement region $\Omega_p$) and to the currents $I_i$ going through the external coils -- located at the domains $\Omega_{C_i}$ with cross sectional area in the $(x,y)$ plane $S_i$. 

The boundary of the plasma is determined by the largest closed level set of $\psi$ that does not intersect any structure in the interior of the reactor. In equation \eqref{eq:FreeBoundarya} we make this dependence explicit by denoting the confinement region as $\Omega_p(\psi)$. This dependence implies that the boundary of the confinement region is also an unknown, making this a \textit{free boundary problem}. When $\partial\Omega_p$ happens to go through a saddle point of $\psi$, it is often referred to as a \textit{separatrix}, and the saddle point---which becomes of physical relevance---is known as an \textit{x-point}. In Section \ref{sec:LimitedAndDiverted} we will discuss other scenarios where the plasma boundary does not coincide a separatrix.

Well-posedness requires the following conditions on the solution at the symmetry axis and at infinity: 
\begin{equation}\label{eq:FreeBoundaryb}
\psi(0,y) = 0 ; \qquad \qquad \lim_{\|(x,y)\|\to\infty}\psi(x,y) = 0. 
\end{equation}
\end{subequations}
We are interested in the impact of uncertainty on this model.  Perhaps the most obvious sources of uncertainty are the values of the currents $\{I_i\}$ that go through the coils, which are subject to variations due, for instance, to small oscillations in the power supply or temperature variations and material impurities in the conducting wire that affect the conductivity. Even if the reactor has been built to very high engineering standards, the precise location of the coils within the reactor with respect to its projected location in the theoretical design will also be subject to a certain amount of uncertainty. Yet one more source of uncertainty is the value of the magnetic permeability, which can only be determined approximately through measurements, and under experimental conditions its value will most certainly differ from the one considered in 
theoretical computations.

Uncertainty of this type adds complexity to an already complicated mathematical and computational problem, but its potential impact on the numerical computation of an equilibrium configuration cannot be avoided. The variability in the parameters will have an impact on the properties of $\psi$ and, since the plasma is confined to within a region determined by its level sets, it will also have an effect on the confinement properties and the shape of the plasma boundary. We will focus on the effect of current variability on the coils producing the poloidal magnetic field only and will take the coil locations and the magnetic permeability to be certain.

For a given range of current variability in the coils, we seek to answer questions on three different levels.  First, we would like to determine the probability that the magnetic field will be conducent to confinement, i.e., whether the level sets of $\psi$ produce closed lines completely contained in the vacuum chamber. Whenever the answer to this question is affirmative, then the second level of question pertains to the nature of the plasma boundary: is it determined by an uninterrupted separatrix (in which case the locations of the x-point and strike points are of relevance) or is it determined by contact with any structure within the reactor (in which case the location of the contact point  is also of interest)?  Finally, for all cases where confinement is achievable (regardless of the presence of an x-point or a contact point), geometric properties of the plasma boundary such as aspect ratio, elongation and triangularity --- see Section \ref{sec:mag-equal} for definitions --- are of interest.

To explore these issues using computational methods, it is of course necessary to work with a discrete solution to (\ref{eq:FreeBoundary}), $\psi_h$, which we will construct using finite element methods, including use of adaptive mesh refinement to resolve complex phenomena near the free boundary. A straightforward way to proceed is then to perform a Monte Carlo simulation: solve the discrete system for a number of realizations of the parameter set and  compute sample statistics such as means, variances and probability estimates using the sample solutions.  Since these computations require solution of the discrete version of equation (\ref{eq:FreeBoundary}) for every realization, simulation done this way is expensive.  We consider as an alternative the use of the \textit{sparse-grid stochastic collocation method} \cite{Barthelmann-Novak-Ritter,Sm:1963,CWM:xiu} to construct a surrogate approximation $\hat \psi_h$ of $\psi_h$, which is less expensive to use in a simulation. Thus, the goal of this work is to examine the effectiveness in terms of accuracy and efficiency of the stochastic collocation method for understanding the effects of current variability on free boundary magnetic equilibrium computations in tokamaks. As we shall show, the use of a surrogate function can reduce the computational time by factors that range from 7 to larger than 30 without significant loss of accuracy in the computations.

An overview of the rest of the paper is as follows.  In Section \ref{sec:mag-equal}, we provide details on the model of magnetic equilibrium in fusion reactors that give rise to the Grad-Shafranov equation (\ref{eq:FreeBoundary}). In Section \ref{sec:Variational}, we present the variational formulation of this equation and outline the discretization strategy used. In Section \ref{sec:Collocation}, we describe the stochastic collocation method for constructing the surrogate approximation $\hat \psi_h$. In Section \ref{sec:NumericalExperiments}, we show the results of numerical experiments using this approach, including a demonstration of its efficiency.

\section{Magnetic equilibrium in fusion reactors} \label{sec:mag-equal}

The equilibrium between magnetic and hydrodynamic pressure in a plasma can be expressed through the equation
\begin{equation}\label{eq:equilibrium}
\text{grad } p = \mathbf J \times \mathbf B,
\end{equation}
where $p$ represents the hydrodynamic pressure in the plasma, $\mathbf J$ is the total current density in the system and $\mathbf B$ is the total magnetic field. Therefore, the region where the plasma can be confined is delimited by the surface of constant pressure where this relation holds. This region is, in general, not known a priori and must be determined as part of an equilibrium calculation.
In addition, the plasma must satisfy the following two fundamental equations:
\begin{align}
\label{eq:AmperesLaw}
\text{curl } \left(\mu^{-1}\,\mathbf B \right) =\;& \mathbf J   \quad \text{(Amp\`ere's Law),} \\
\label{eq:NoMagneticMonopole}
\text{div } \mathbf B =\;& 0.
\end{align}
Above, $\mu$ represents the magnetic permeability; note that $\mathbf J$ and $\mathbf B$ account for the joint contributions coming from both the array of coils used to generate the external magnetic field, and the plasma itself. 

\subsection{The Grad-Shafranov equation}\label{sec:Grad-Shaf}
In an axially symmetric setting, described by cylindrical coordinates $(r,z,\phi)$, the magnetic field $\mathbf B$ can be expressed in terms of two scalar functions $\psi(r,z)$ and $g(r,z)$ as
\begin{equation}\label{eq:representation}
\mathbf B = \widehat{\boldsymbol\phi} \times \frac{\text{grad } \psi}{r} + \widehat{\boldsymbol \phi}\,\frac{g}{r},
\end{equation}
where $\widehat{\boldsymbol\phi}$ is the {unit} vector in the
direction of increasing angle $\phi$. In this setting, $r$ and $z$ merely correspond to Cartesian coordinates in the poloidal plane, which justifies the subsequent simplifying use of the variables $(x,y)$ and of the symbols $\nabla$ and $\nabla\cdot$ to represent the Cartesian gradient and divergence operators respectively.

Both scalar functions appearing in \eqref{eq:representation} have very precise physical interpretations. The poloidal flux $\psi$ determines the component of the magnetic field on a plane of constant angle $\phi$ --- the \textit{poloidal component} --- and the function $g$ is related to the total current flowing through the system in the poloidal direction, $I_p$, through the relation $I_p = -2\pi g$ (see \cite{Freidberg2014,GoKePo2010}). It can be shown that both $g$ and the pressure $p$ are functions of the poloidal flux $\psi$ alone. This implies that the lines of constant pressure correspond to the level sets of 
$\psi$, which will also be referred to in the sequel as the \textit{stream function}.
The plasma can only remain confined in the region of the plane where the total field lines are closed, $\Omega_p$. However, since the total magnetic field has a contribution due to the plasma, whose precise location is unknown, the problem becomes a non-linear free boundary problem, where both the magnetic field and $\Omega_p$ must be determined simultaneously.

Using the representation of the magnetic field given by \eqref{eq:representation} in combination with conditions \eqref{eq:AmperesLaw} and \eqref{eq:NoMagneticMonopole}, the equilibrium condition \eqref{eq:equilibrium} can be expressed as the differential relation (\ref{eq:FreeBoundary}),
the Grad-Shafranov equation.
The form of the source term in the plasma that we will use in this paper is
\begin{equation}\label{eq:source}
\frac{d}{d\psi}p(\psi) = \RefereeA{\lambda}\frac{\beta}{x_0}\left(1-\psi_N^\alpha\right)^\gamma,  \qquad
\frac{1}{2}\frac{d}{d\psi}g^2(\psi) = \RefereeA{\lambda}\mu_0x_0(1-\beta)\left(1-\psi_N^\alpha\right)^\gamma.
\end{equation}
\RefereeA{The parameter $\lambda$ appearing above is a scaling factor introduced to enable a target value of $I_p$ in the equilibrium computation.} These forms were proposed by Luxon and Brown \cite{LuBr:1982} as the simplest model that was statistically adequate to describe experimental results and have been widely used for theoretical studies ever since. The term
\[
\psi_N : = \frac{\psi - \psi_{BD}}{\psi_{MA} - \psi_{BD}}
\]
is a normalization of the unknown with respect to the value of $\psi$ at the plasma boundary, $\psi_{BD}$, and the maximum value of $\psi$ within the plasma, $\psi_{MA}$, which determines the location of the point known as the \textit{magnetic axis}. The value $x_0$ corresponds to the outer radius of the vacuum chamber, the parameters $\alpha$ and $\gamma$ control how sharply the current peaks near the magnetic axis; $\beta$ is known as \textit{poloidal beta}, a parameter measuring the ratio between the hydrostatic pressure in the plasma due to temperature (usually called \textit{plasma pressure}) and the pressure due to the poloidal magnetic field (or \textit{magnetic pressure}). Since the aim of magnetic confinement is to balance the effects of the internal hydrostatic pressure with those of the external magnetic field, the value of $\beta$ effectively measures the strength of confinement of a reactor. All these parameters are determined by the user, depending of the design of the reactor or the particularities of an experiment, mainly the desired current and pressure profiles in the plasma. 

\subsection{Limited and diverted configurations}\label{sec:LimitedAndDiverted}
\begin{figure}\centering
\begin{tabular}{llll}
\includegraphics[height = 0.19\linewidth]{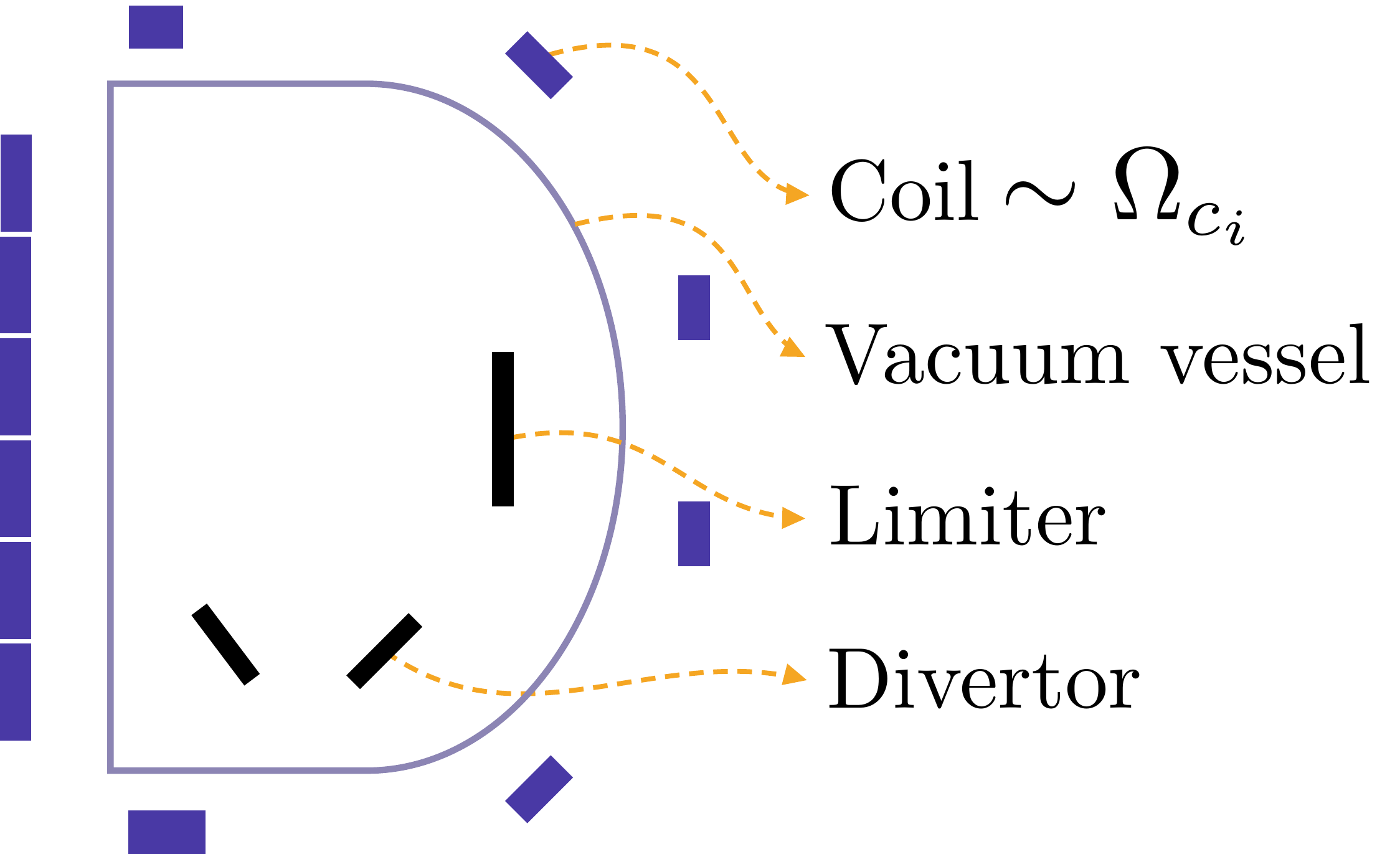} \qquad \qquad &
\includegraphics[height = 0.19\linewidth]{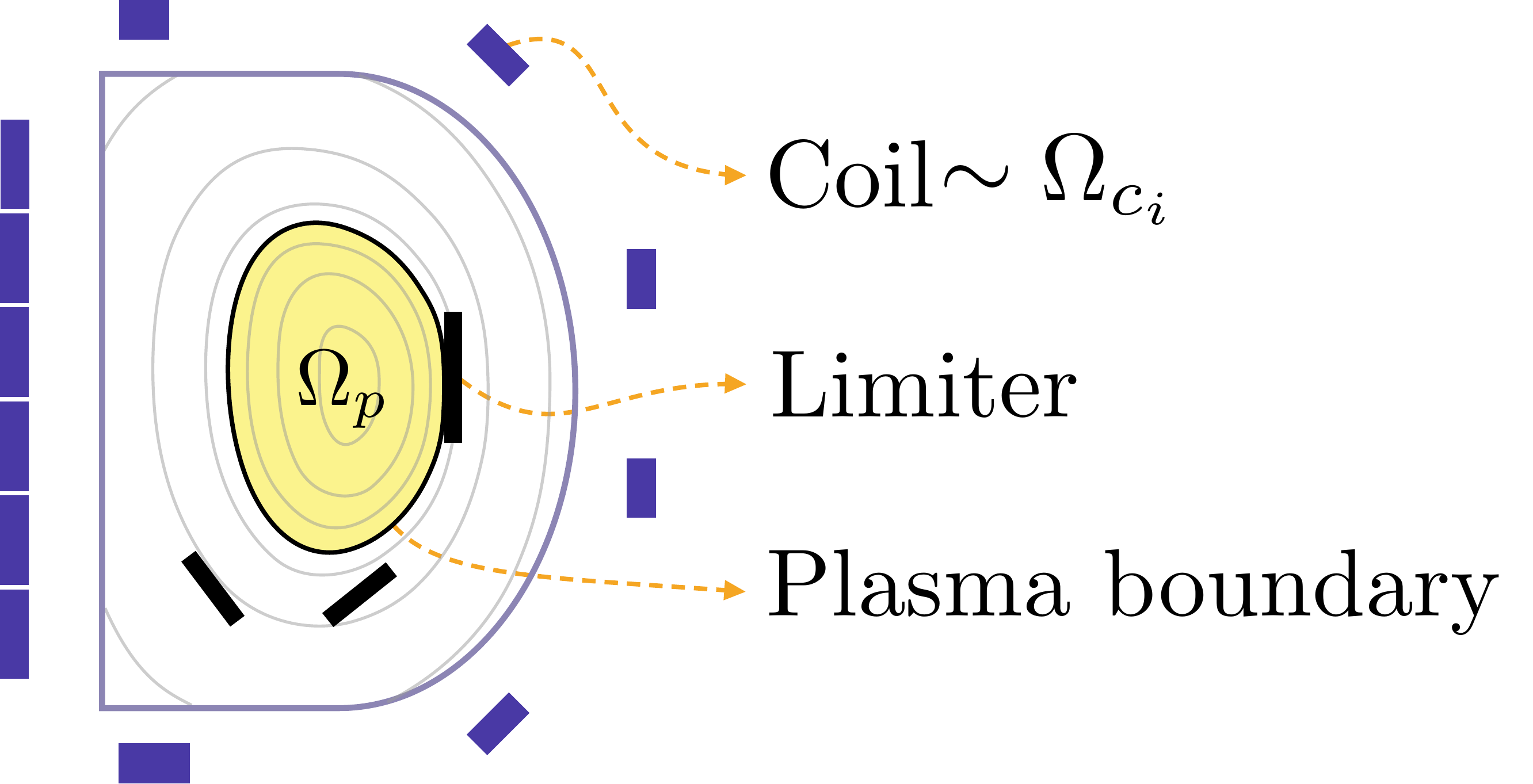}    \\ [2ex]
\includegraphics[height = 0.19\linewidth]{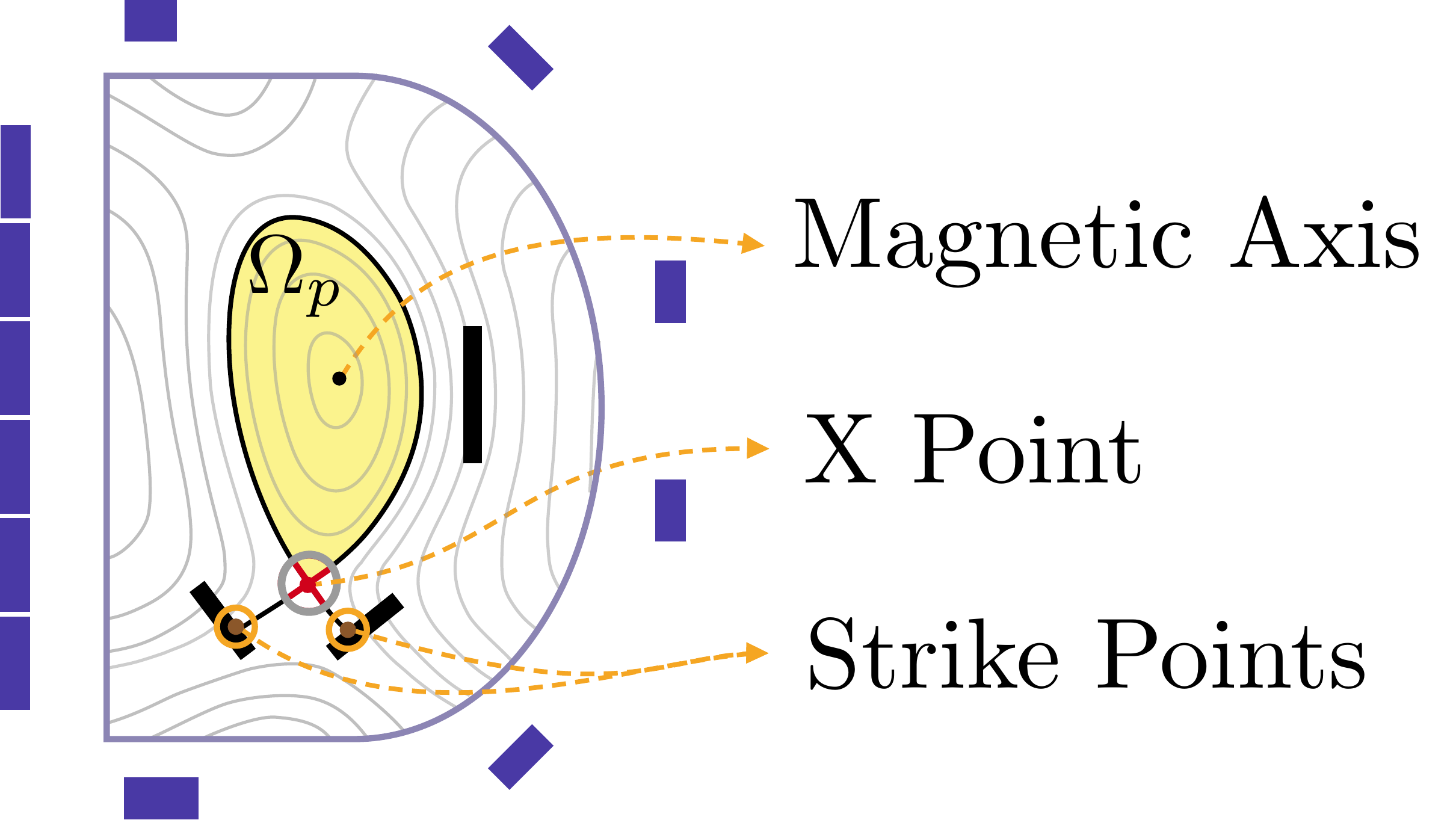} \qquad \qquad &
\includegraphics[height = 0.19\linewidth]{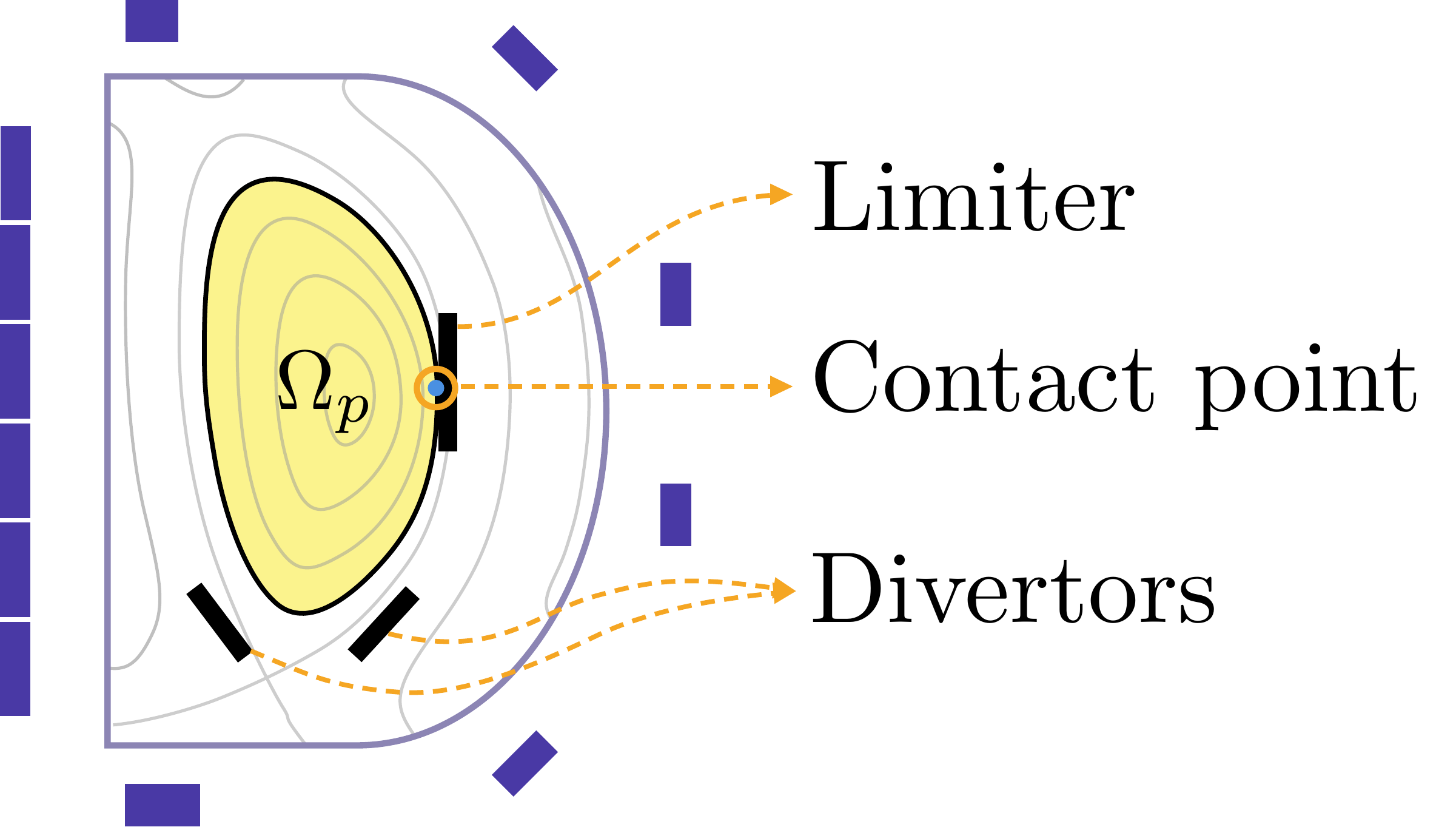}
\end{tabular}
\caption{Schematics of a cross section of a tokamak and possible plasma configurations. Top left:  simplified cross-section of a tokamak depicting some of the essential components for magnetic confinement. Top right: a \textit{limited configuration}. Bottom left: a \textit{diverted configuration}. Bottom right: the separatrix intersects the limiter plate.}
\label{fig:tokamak}
\end{figure}

Several scenarios can occur for plasma contained in a tokamak.
Cross-sections of typical tokamak configurations are shown in Figure \ref{fig:tokamak}. 
Depending on the configuration, in the interior of the vacuum vessel either divertor or limiter plates (solid black) are installed to protect the walls. 
In the generic top left image, the plasma is contained in a hollow toroidal vacuum vessel (violet line) and the confinement field is generated by an external array of coils (purple blocks). 

Depending on the location of the coils and the intensity of the current passing through them, there are three possible configurations of the magnetic field that are relevant for confinement. The top right image in Figure \ref{fig:tokamak} shows a \textit{limited} configuration, where the contour lines of $\psi$ contained in the vacuum vessel are all closed and nested, interrupted only by the walls of the chamber. In this situation, common in older reactors, the plasma could come in contact with the walls and therefore a physical barrier---typically a tungsten plate depicted as a solid black line in the schematic---has to be positioned to prevent the plasma from reaching the walls \cite{Wesson:2011}. The level set of $\psi$ tangential to the limiter  determines the boundary of the confinement region $\Omega_p$. While this solution succeeds at shielding the walls of the reactor, the fact that the plasma is in direct contact with the limiter has some undesirable effects including a reduced range of temperatures attainable and the introduction of pollution due to the ablation of the limiter plate.

Alternatively, as in the bottom left image of Figure \ref{fig:tokamak}, the magnetic field can give rise to a situation where the interior of the vessel is divided into a sector where the contour lines are closed and one in which they are open. 
In this \textit{diverted} configuration,  the stream function has a separatrix that divides the interior of the vessel into these two regions.
As confinement is possible only in regions where the field lines are closed, if this line is not obstructed by any structure, then the plasma remains confined within the separatrix without the need of a limiter. The point where the separatrix intersects itself is called an \textit{x-point}. Since heavy ions or impurities leaving the plasma do so by moving along magnetic stream lines, the vicinity of the x-point, where the separatrix intersects itself, is an ideal location for placing a filtering device known as \textit{divertor} that removes such ions and absorbs the excess heat, protecting the walls without coming directly in contact with the bulk of the plasma \cite{DiChAnFeJaMaPaTi:1995, JaBoFeIgKuPaPaSu:1995}. 
The area surrounding the points where the separatrix intersects the divertor plates---the \textit{strike points}---is also of interest, since it is the most likely to get in contact with fast particles and impurities leaving the plasma.

Owing to these and other advantages, reactors with magnetic x-points and diverted configurations are favored nowadays; designs with two or more divertors have been proposed and studied in \cite{BoOkFo:1985,Kesner:1990,RySo:2015}. 
It is clear that  accurate determination of the magnetic equilibrium configuration and in particular the location of a magnetic x-point on a given reactor is of great relevance for theoretical, experimental and design purposes.

A third possibility, shown in the bottom right schematic of Figure \ref{fig:tokamak},
is that the stream function gives rise to a separatrix that intersects the walls of the vacuum chamber or some other structure contained in it, such as the divertor or a limiter.  In this case,  although confinement is still possible, the separatrix no longer determines the plasma boundary, as the obstructing structure acts effectively as a limiter. Instead, the plasma boundary is determined by the first closed contour inside the separatrix that is tangential to said structure
at a contact point. This scenario is  undesirable, for the plasma walls or structures inside of the vacuum chamber, other than the limiter and divertor plates, are not designed to be in direct contact with the plasma and can be severely damaged by such interactions. As will be seen later, even if reactors are designed to avoid this scenario, unexpected variations in current intensity can give rise to such undesired contacts. Being able to quantify the probability of this happening (and the most likely location of contact with any structure) is therefore of paramount importance.

\subsection{Plasma shaping} \label{sec:PlasmaShaping}
Some geometric properties of the plasma boundary have been found to have significant influence in important physical processes such as stability and transport. Substantial theoretical and experimental efforts have been devoted to the design and realization of particular confinement geometries in what has come to be known as \textit{plasma shaping} \cite{shaping3,shaping2,shaping5,shaping1,shaping4}.

\begin{figure}\centering
\includegraphics[height = 0.25\linewidth]{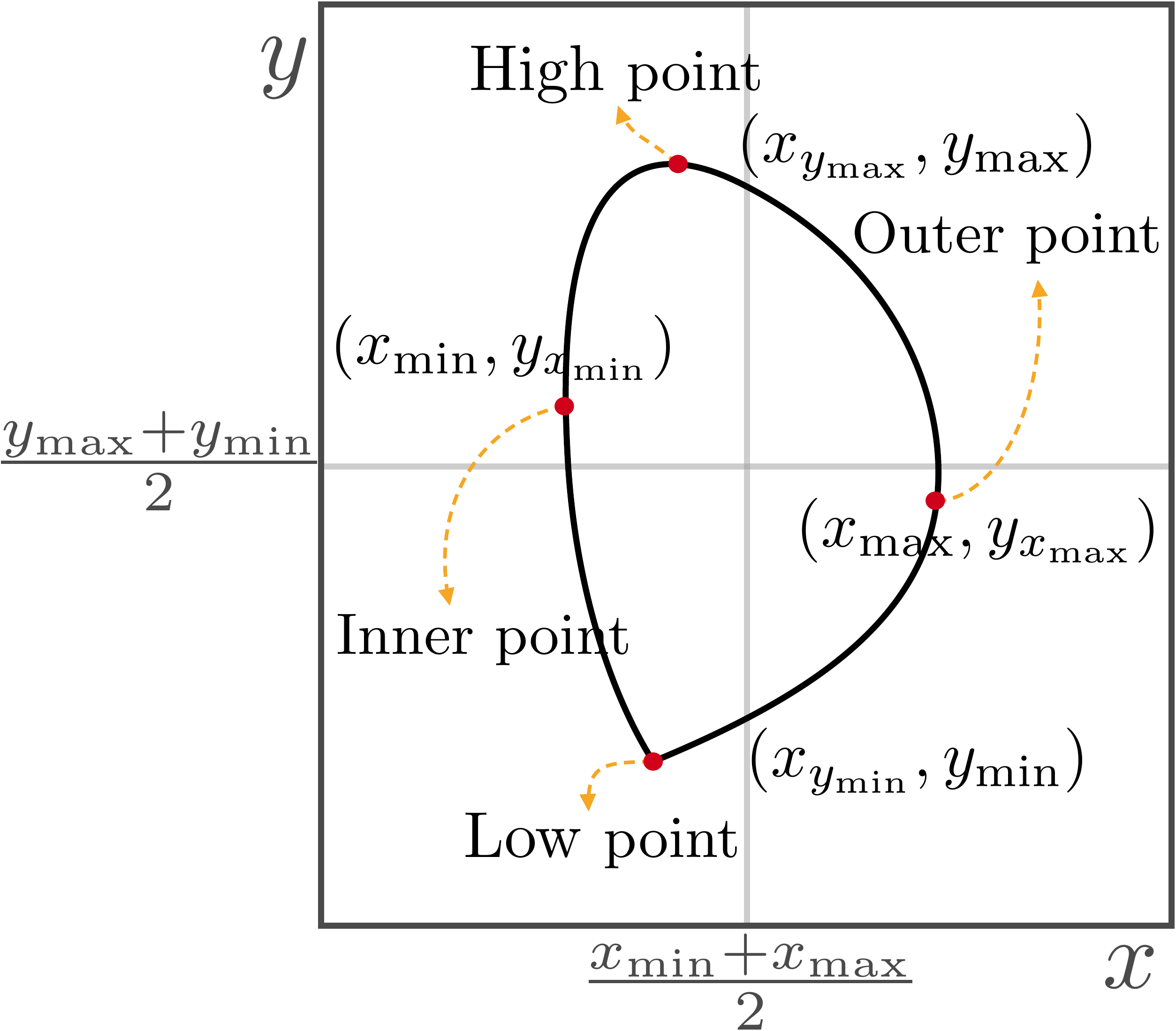} \qquad \qquad \qquad \qquad
\includegraphics[height = 0.25\linewidth]{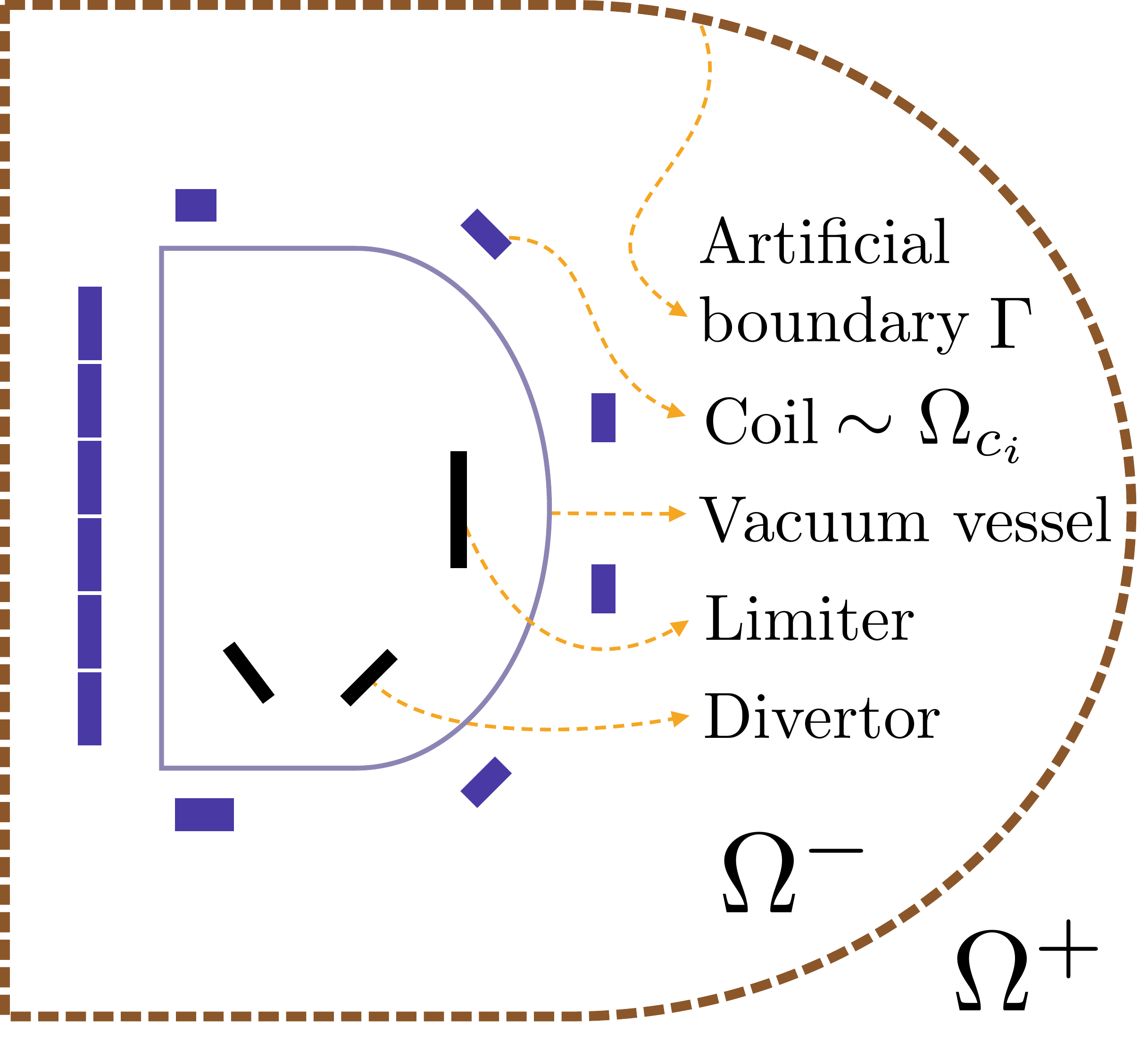}
\caption{Left: Locations of the high, low, inner and outer points in a typical plasma boundary. These points are used to define the shaping parameters in \eqref{eq:ShapeParameters}: major and minor radii, aspect ratio, elongation and triangularity. Right: An artificial boundary $\Gamma$ enclosing the reactor and all the relevant geometric structures is introduced. The region contained within $\Gamma$ is meshed with a geometry-fitting triangulation, where the equation \eqref{eq:FreeBoundarya} is discretized. For the radiation condition \eqref{eq:FreeBoundaryb} to be satisfied, an integral representation of the solution is used for the region exterior to $\Gamma$. If the artificial boundary is chosen to be a circle, the integral equation can be solved analytically and substituted back into the finite element system as an additional term. }\label{fig:ShapingParameters}
\end{figure}

The geometric features that have been found to be most relevant are referred to as \textit{elongation} and \textit{triang\-ularity}. These terms are {largely} self-explanatory, but to define them properly we need to introduce a few auxiliary quantities---the notation and terminology are standard in the fusion literature, see for instance  \cite{CeFr:2010,Luce:2013}:  
\begin{subequations}\label{eq:ShapeParameters}
\begin{alignat}{6}
R_{geo} :=\,& (x_{max} + x_{min}) / 2 \qquad& \text{ Major radius}, \\
a :=\,& (x_{max} - x_{min}) / 2 \qquad& \text{ Minor radius}, \\
\epsilon :=\,& a/R_{geo} \qquad& \text{Inverse aspect ratio},\\
\kappa :=\,& (y_{max} - y_{min}) / 2a \qquad& \text{ Elongation},\\
\delta_u : =\,& (R_{geo} - x_{y_{max}})/ a \qquad& \text{ Upper triangularity}, \\
\delta_l : =\,& (R_{geo} - x_{y_{min}})/ a \qquad& \text{ Lower triangularity}.
\end{alignat}
\end{subequations}
The left panel of Figure \ref{fig:ShapingParameters} depicts the points that determine the height and major and minor radii on a typical plasma boundary with an x-point. Note that in general the locations of the Inner/Outer and High/Low points need not be aligned, nor does the shape of the plasma need to be symmetric.

\section{Variational formulation and finite element discretization}\label{sec:Variational}
Due to its relevance in fusion science, the mathematical and numerical aspects of the free boundary problem \eqref{eq:FreeBoundary} have been a continuous subject of interest. Among the most relevant recent efforts have been those by Heumann and collaborators \cite{FaHe:2017,Heumann:feeqsm,CEDRES,HeRa:2017}; here we rely on their contributions for the numerical computation of the equilibrium configuration, {in particular}, the {\tt Matlab} implementation of the free boundary solver described in \cite{CEDRES} and kindly provided by the developers to the authors of this paper. We refer the reader interested in the details of the numerical solution of \eqref{eq:FreeBoundary} to those references; here we describe the process in very general terms.
 
The main challenges posed by \eqref{eq:FreeBoundary} are: (1) the unboundedness of the domain where the problem is posed, (2) the nonlinearity of the partial differential equation stemming from the source term and---in the presence of iron components---the permeability coefficient, and (3) the nonlinearity of the problem owing to the unknown nature of the plasma region $\Omega_p$. In the code {\tt FEEQS.m} used for this work, point (1) is dealt with by coupling with a boundary integral equation, while points (2) and (3) are addressed through an iterative Newton approach by linearizing around an appropriately chosen initial guess.

We note that there are other free boundary solvers available in the literature that take different computational approaches \cite{Albanese2015, Fitzgerald2013,Hofmann1988,Mccarthy1999}, and the stochastic collocation approach that we describe in Section \ref{sec:Collocation} can be easily interfaced with any of them.  The most notable differences that set {\tt FEEQS.m} (and the closely related solvers NICE \cite{Faugueras2020} and CEDRES++  \cite{CEDRES}) apart from the rest are:  (a) the use of Newton’s method to resolve the nonlinearities, as the majority of free-boundary solvers used in the magnetic confinement community resort to a fixed-point approach; and (b) the coupling with a boundary integral formulation to address the unboundedness — the most common approach is the von Hagenow–Lackner coupling \cite{Lackner1976,vHaLa1975}. From a mathematical point of view, use of the boundary integral method is perhaps the difference that requires the most elaborate treatment and we will describe it in some more detail in what follows.  In doing so, we will make use of some results and terminology from the theory of boundary integral equations. We refer to \cite{HsWe:2004} for a concise overview of the subject.

\subsection{The variational formulation}
An artificial boundary $\Gamma$ enclosing the entire reactor and remaining structures is introduced, as shown in the right panel of Figure \ref{fig:ShapingParameters}. $\Gamma$ must be such that the support of the right hand side of \eqref{eq:FreeBoundarya} is compactly contained inside. This curve divides the half plane into the region enclosed by it, which we will denote by $\Omega^-$,  and its unbounded complement $\Omega^+:=\mathbb{R}^2\setminus\overline{\Omega^-}$. For simplicity we will leave $\Gamma$ as general as possible, but its precise form will be made explicit in due time.

The solution to the free boundary problem can then be considered to be of the form $\psi=\psi^- + \, \psi^+$, where the two parts $\psi^-$ and  $\psi^+$ are supported inside $\Omega^-$ and $\Omega^+$ respectively. The continuity of $\psi$ is preserved by requiring both interior and exterior solutions $\psi^-$ and  $\psi^+$ (and their normal derivatives) to match on the interface $\Gamma$. More precisely, problem \eqref{eq:FreeBoundary} can be reformulated as a system consisting of the interior problem
\begin{subequations}\label{eq:InteriorExteriorFormulation}
\begin{alignat}{6}
\label{eq:InteriorExteriorFormulationA}
-\nabla\,\cdot\,\left(\frac{1}{\mu x}\nabla\psi^-\right) =\,& F(\psi^-,x,\Omega_p) := \left\{ \begin{array}{ll}
\RefereeA{x\frac{d}{d\psi^-} p(\psi^-) + \frac{1}{2\,\mu x} \frac{d}{d\psi^-} g^2(\psi^-)} & \text{ in } \Omega_p(\psi^-) \\
I_i/S_i & \text{ in } \Omega_{C_i} \\
0 & \text{ Elsewhere. } 
\end{array}\right. &\qquad& \text{ in } \Omega^- \\
\label{eq:InteriorExteriorFormulationB}
\psi^-(0,y) =\,& 0, &\qquad& 
\end{alignat}
and the exterior problem
\begin{alignat}{6}
\label{eq:InteriorExteriorFormulationC}
-\nabla\,\cdot\,\left(\frac{1}{\mu x}\nabla\psi^+\right) =\,& 0 &\qquad& \text{ in } \Omega^+ \\
\nonumber
\psi^+(0,y) =\,& 0, &\quad &\lim_{\|(x,y)\|\to\infty}\psi^+(x,y) = 0,  
\end{alignat}
connected through the continuity conditions
\begin{alignat}{6}
\label{eq:InteriorExteriorFormulationF}
\psi^- =\,& \psi^+ &\qquad & \text{ on } \Gamma \\
\label{eq:InteriorExteriorFormulationG}
\partial_{\nu}\psi^- =\,& \partial_{\nu}\psi^+ &\qquad& \text{ on } \Gamma.
\end{alignat}
\end{subequations}
Testing with a function $\varphi$ satisfying the boundary condition \eqref{eq:InteriorExteriorFormulationB} but otherwise arbitrary, one arrives at
\begin{subequations}\label{eq:WeakFormulation}
\begin{equation}\label{eq:WeakFormulationA}
\int_{\Omega^-}\frac{1}{\mu x}\nabla\psi^-\cdot\nabla\varphi\,dA -\int_\Gamma\frac{1}{\mu x}\partial_{\nu}\psi^-\,\varphi\,d\Gamma  = \int_{\Omega^-}F(\psi^-,x,\Omega_p)\,\varphi\,dA .
\end{equation}
 We now turn our attention to the exterior problem. In order to deal with the unbounded domain we will use an integral representation of the exterior solution of the form 
\[
\psi^+ = \int_\Gamma G \, \partial_{\nu}\psi^+ \, d\Gamma - 
\int_\Gamma  \partial_{\nu}G\, \psi^+ \,d\Gamma, 
\]
where $G = G(\boldsymbol x_1,\boldsymbol x_2)$ is the Green's function for \eqref{eq:FreeBoundary} vanishing on $\Gamma$. To avoid cumbersome notation, the explicit expression of the Green's function will not be given in this discussion. The terms in the right hand side are known respectively as \textit{single layer potential} and \textit{double layer potential}. Note that the values of the restrictions to the boundary and the normal derivatives appearing in the right hand side are not known a priori and have to be determined as part of the solution process.

Using the expression above to take the restriction of $\psi^+$ to the boundary $\Gamma$ enables the use of the \textit{jump properties} of the single and double layer potentials \cite{HsWe:2004}, which yield the following integral equation on $\Gamma$
\[
\psi^- = \frac{1}{2}\psi^- + \int_\Gamma \partial_{\nu}G\,\psi^- \,d\Gamma + \int_\Gamma G\, \partial_{\nu}\psi^-\,d\Gamma
\]
where \eqref{eq:InteriorExteriorFormulationF} and \eqref{eq:InteriorExteriorFormulationG} were used to connect the exterior and interior problems  by writing the integral equation in terms of the restrictions of the interior solution only. This equation can be posed weakly by testing with an arbitrary function $\lambda$ defined on $\Gamma$, yielding
\begin{equation}\label{eq:WeakFormulationB}
{\renewcommand{\arraystretch}{2.2}
\begin{array}{l}
{\displaystyle 
\int_\Gamma\frac{1}{2}\psi^-\,\lambda\,d\Gamma - \int_\Gamma \int_\Gamma \partial_{\nu(\boldsymbol x_1)}G(\boldsymbol x_1,\boldsymbol x_2)\,\psi^- (\boldsymbol x_1)\,\lambda(\boldsymbol x_2) \, d\Gamma(\boldsymbol x_1)\,d\Gamma(\boldsymbol x_2)} \\
\hspace{1.2in}
\displaystyle{- \int_\Gamma \int_\Gamma G(\boldsymbol x_1, \boldsymbol x_2)\, \partial_{\nu(\boldsymbol x_1)}\psi^-(\boldsymbol x_1)\,\lambda(\boldsymbol x_2)\,d\Gamma(\boldsymbol x_1)\,d\Gamma(\boldsymbol x_2) = 0.}
\end{array}
}
\end{equation}
\end{subequations}
Equations \eqref{eq:WeakFormulationA}-\eqref{eq:WeakFormulationB} involve only the unknowns $\psi^-$ and  $\partial_{\nu}\psi^-$ so, if they can be solved simultaneously, the solution $\psi$ to the original problem can be recovered. 

So far, we have followed the general strategy introduced by Johnson and N\'ed\'elec in their celebrated article \cite{JoNe:1980}. However, as shown by Hsiao and collaborators \cite{GaHs:1995,HsZh:1994}, if $\Gamma$ is taken to be a semi-circle, equation \eqref{eq:WeakFormulationB} can be solved exactly for $\partial_{\nu}\psi^+$ in terms of the tangential derivative of $\psi^-$ along $\Gamma$. Substitution of this expression for $\partial_{\nu}\psi^+$ back into \eqref{eq:WeakFormulationA} gives the following uncoupled equation for $\psi^-$, which must hold for all sufficiently regular test functions $\varphi$ vanishing on the $y$ axis:
\[
\int_{\Omega^-}\frac{1}{\mu x}\nabla\psi^-\cdot\nabla\varphi\,dA -2\int_\Gamma \!\int_\Gamma\frac{1}{\mu x_1}G(\boldsymbol x_1, \boldsymbol x_2)\,\partial_{\tau(\boldsymbol x_1)}\psi^-(\boldsymbol x_1)\,\partial_{\tau(\boldsymbol x_2)}\varphi(\boldsymbol x_2)\,d\Gamma(\boldsymbol x_1)\,d\Gamma(\boldsymbol x_2)  = \!\! \int_{\Omega^-}\!\!\!F(\psi^-\!,x,\Omega_p)\,\varphi\,dA.
\]
The symbol $\partial_\tau$ denotes (weak) tangential differentiation along the curve $\Gamma$. In this equation, the unknowns are the stream function inside of the fictitious boundary, $\psi^-$, and the support of the plasma, $\Omega_p$, which appears through the source term of the equation $F(\psi^-,x,\Omega_p)$---see \eqref{eq:InteriorExteriorFormulation}. The boundary of the plasma region is the \textit{closed} level set of $\psi^-$ that either: a) passes through a saddle point of $\psi^-$ located inside of the reactor, or b) is tangential to the inner wall of the reactor---if the curve defined by a) intersects any structure inside the reactor or if no curve satisfies the condition a). See Figure \ref{fig:tokamak} for graphical depictions of these conditions. 

Once $\psi^-$ is found, the continuity conditions \eqref{eq:InteriorExteriorFormulationF} and \eqref{eq:InteriorExteriorFormulationG} along with the integral representation for $\psi^+$ can be used to evaluate $\psi$ in the unbounded component $\Omega^+$ if needed. Note that even if it is supported on the entire component $\Omega^+$, the function $\psi^+$ appears in these equations only through the integrals over $\Gamma$, which avoids the discretization of an unbounded domain.

The explicit expression for the equation above, derived tersely in \cite{AlBlDe:1986}, is obtained by substituting the appropriate Green's function and integrating by parts the term involving the tangential derivatives. We can now state the variational problem that must be solved. 

\vspace{.1in}
\noindent
{\bf Variational Formulation.}
Consider the space 
\[
V := \left\{ \psi:\Omega^- \to \mathbb R \,\Bigg| \, \int_{\Omega^-} (\psi)^2 x\,dx\,dy < \infty\,;\,\int_{\Omega^-} (\nabla \psi)^2 x^{-1}\,dx\,dy < \infty\,;\, \psi(0,y) = 0 \right\} \cap \mathcal C^0\left(\overline{\Omega^-}\right).
\]
Let $\Gamma$ be a circle of radius $\rho_\Gamma$, centered at the origin and completely containing the reactor geometry. 
The problem is that of finding $\psi$ belonging to $V$ and such that for all test functions $\varphi \in V$ it  holds that
\begin{align}
\nonumber
& \int_{\Omega^-}\frac{1}{\mu x}\nabla\psi\cdot\nabla\varphi\,dA +\int_\Gamma \psi\,N\,\varphi\,d\Gamma \\ 
\label{eq:VariationalFormulation}
&\hspace{.2in}
 + \int_\Gamma \int_\Gamma  \left(\psi(\boldsymbol x_1) - \psi(\boldsymbol x_2)\right) M(\boldsymbol x_1,\boldsymbol x_2)  \left(\varphi(\boldsymbol x_1) - \varphi(\boldsymbol x_2)\right) \,d\Gamma(\boldsymbol x_1)\,d\Gamma(\boldsymbol x_2) =  \int_{\Omega^-}F(\psi,x,\Omega_p)\,\varphi\,dA,
\end{align}
where $\boldsymbol x_i = (x_i,y_i)$ and
\begin{alignat*}{6}
N :=\,& \frac{1}{x} \left(\frac{1}{\delta_+} + \frac{1}{\delta_-} - \frac{1}{\rho_\Gamma}\right)\, , \\
\delta_\pm :=\,& \sqrt{x^2 + (\rho_\Gamma \pm y)^2}\, ,
\\
M(\boldsymbol x_1,\boldsymbol x_2) :=\,& \frac{\kappa(\boldsymbol x_1,\boldsymbol x_2)}{2\pi(x_1x_2)^{3/2}} \left(\frac{2-\kappa^2(\boldsymbol x_1,\boldsymbol x_2)}{2-2\kappa^2(\boldsymbol x_1,\boldsymbol x_2)}E(\kappa(\boldsymbol x_1,\boldsymbol x_2))- K(\kappa(\boldsymbol x_1,\boldsymbol x_2))\right)\, ,\\
\kappa(\boldsymbol x_1,\boldsymbol x_2):=\,& \sqrt{\frac{4x_1x_2}{(x_1 + x_2)^2 + (y_1-y_2)^2}}\, ,
\end{alignat*}
and $E(\kappa(\boldsymbol x_1,\boldsymbol x_2))$ and $ K(\kappa(\boldsymbol x_1,\boldsymbol x_2))$ are, respectively, complete elliptic integrals of first and second kind. $\bullet$

\subsection{Finite element discretization
\label{sec:FEM}}

Computational studies of the free boundary problem
\eqref{eq:FreeBoundary} require an accurate numerical solution $\psi_h$.
For our computations, we used the package {\tt FEEQS.M} \cite{Heumann:feeqsm}, which approximates the solution to the variational formulation \eqref{eq:VariationalFormulation} specified above
using a piecewise linear finite element discretization built on a geometry-fitting mesh. It solves the discrete nonlinear system using Newton's method, starting from an initial value for $\psi$ given by a function of the form $\psi_0 = \frac{(x-x_0)^2}{a^2}+\frac{(y-y_0)^2}{b^2} + K$ for some point $(x_0,y_0)$ inside the vacuum chamber, together with an initial value for $\Omega_p$ taken to be one of the elliptical level sets fully contained in the reactor geometry. We provide additional details on the Newton iteration below. See also \cite{CEDRES} for a detailed description of the numerical solver.

The resolution and accuracy of the solution provided by {\tt FEEQS.M} can be enhanced through mesh refinement. The most natural strategy---uniform refinement---would result in a significant increase in the number of system unknowns. Keeping in mind that several steps of a Newton iteration are needed to compute a single realization of the solution, and that a large number of realizations is necessary 
to sample the parameter space, the additional cost associated with uniform refinement becomes unmanageable. 

\begin{figure}[htb]\centering
\begin{tabular}{ccc}
\includegraphics[width=0.16\linewidth]{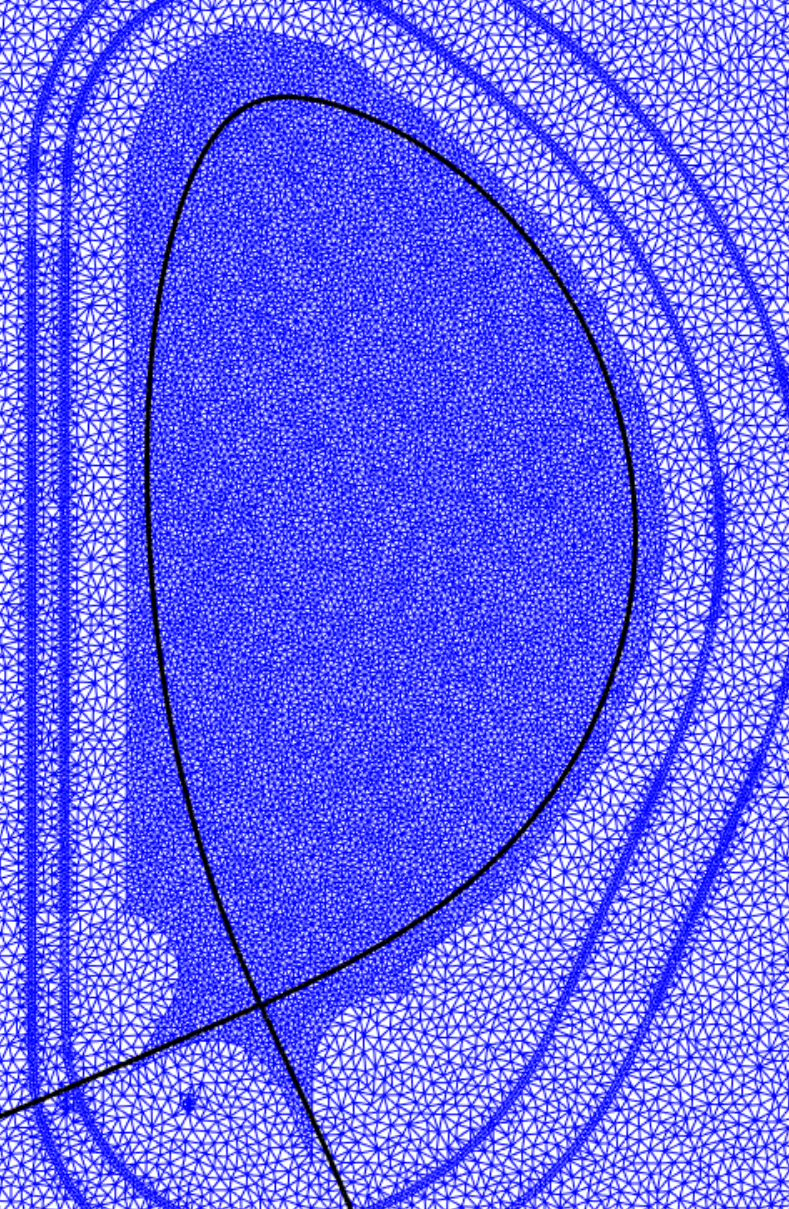} \qquad&
\includegraphics[width=0.16\linewidth]{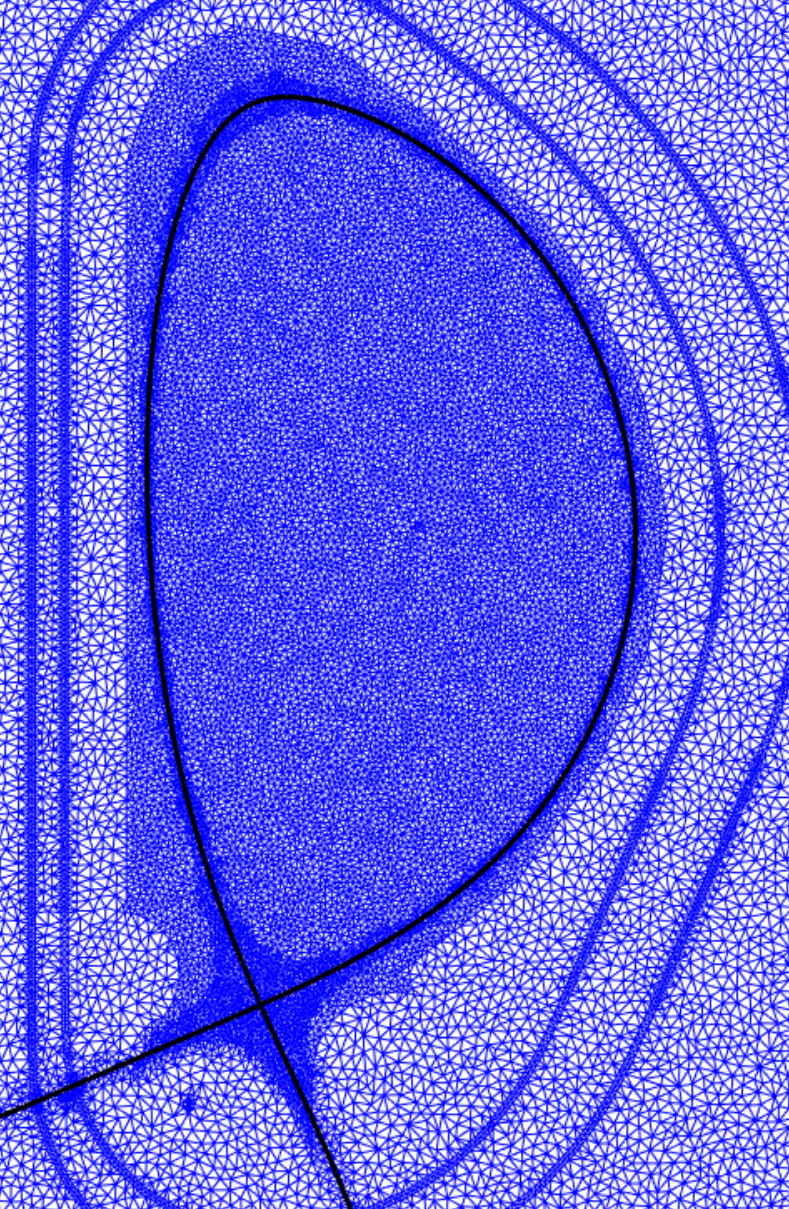} \qquad&
\includegraphics[width=0.16\linewidth]{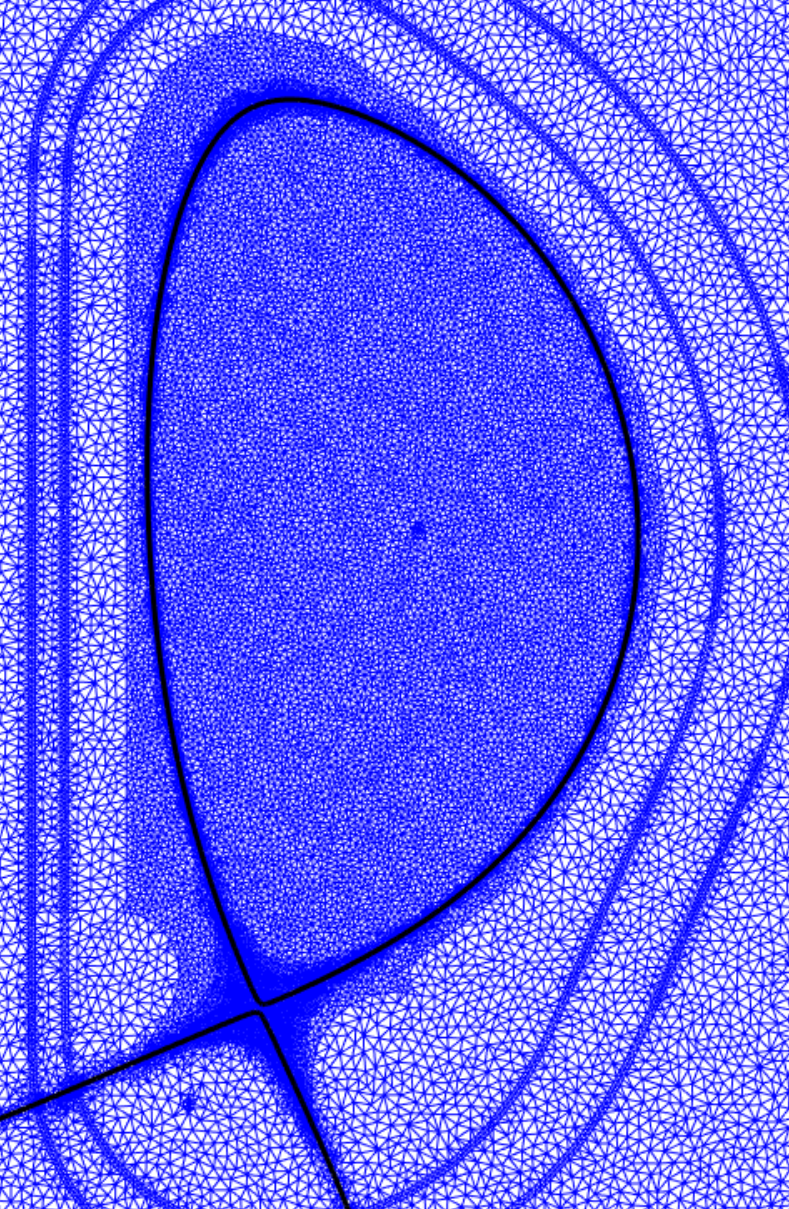} \\
\includegraphics[width=0.16\linewidth]{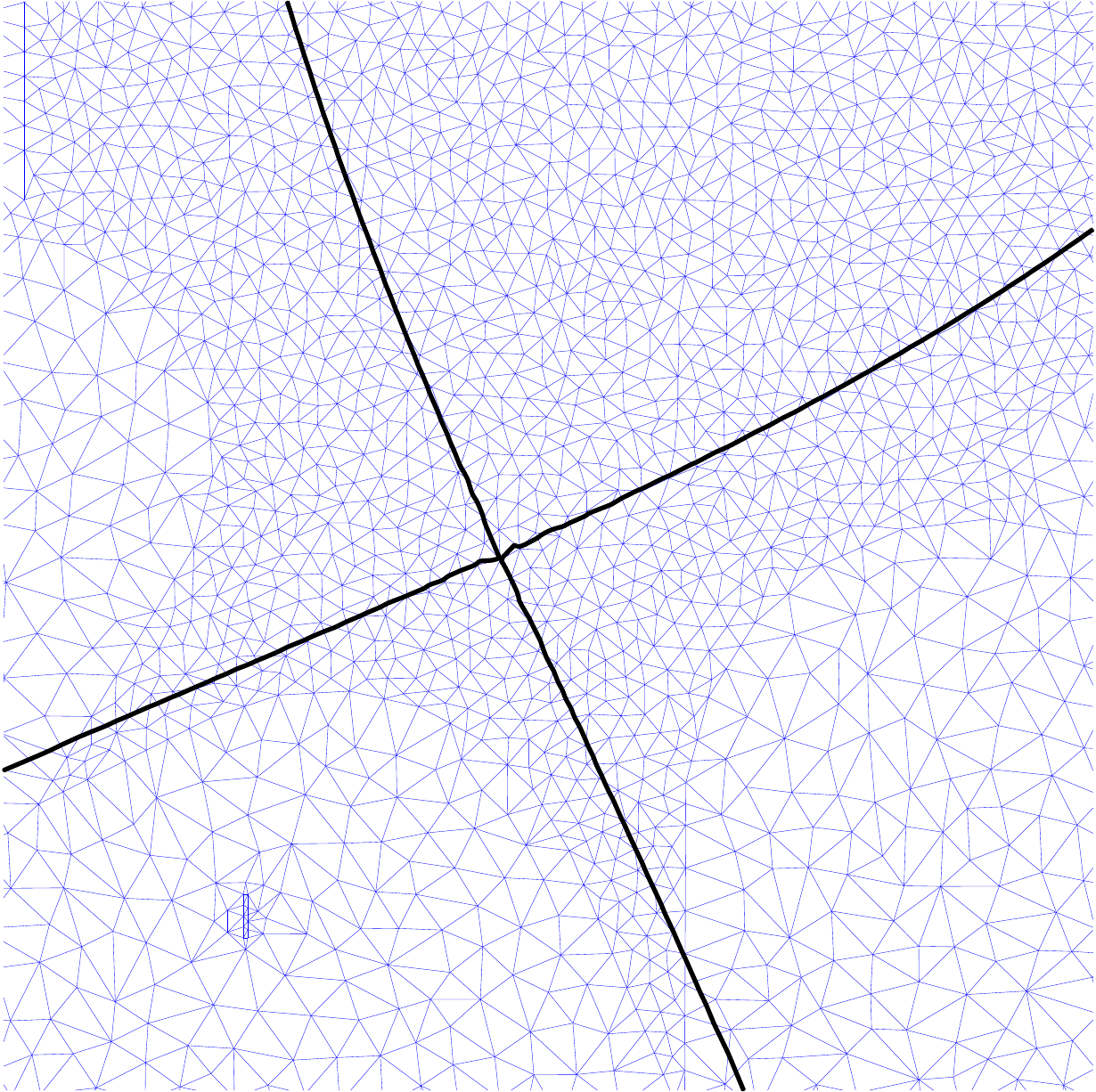} \qquad&
\includegraphics[width=0.16\linewidth]{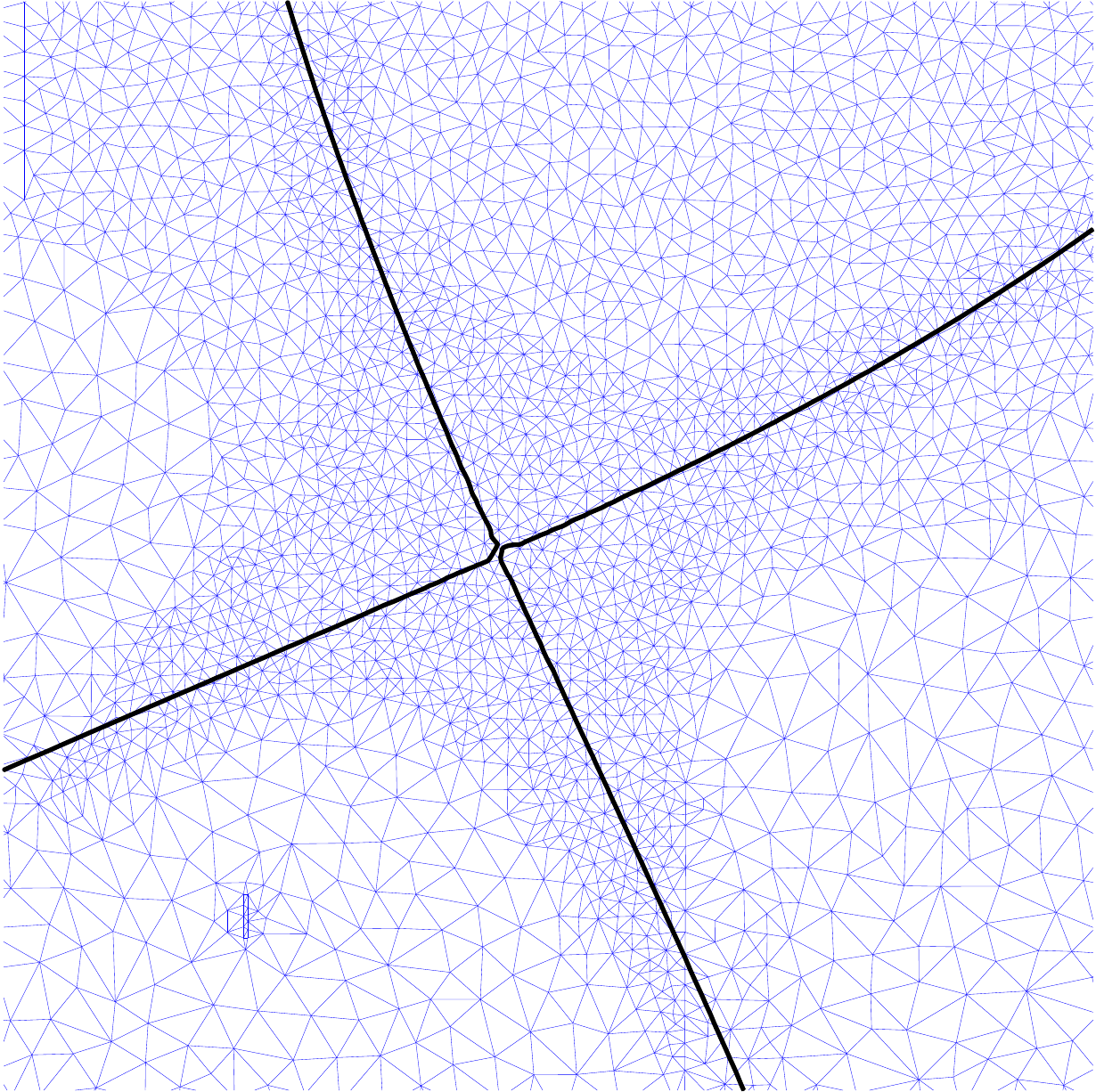} \qquad&
\includegraphics[width=0.16\linewidth]{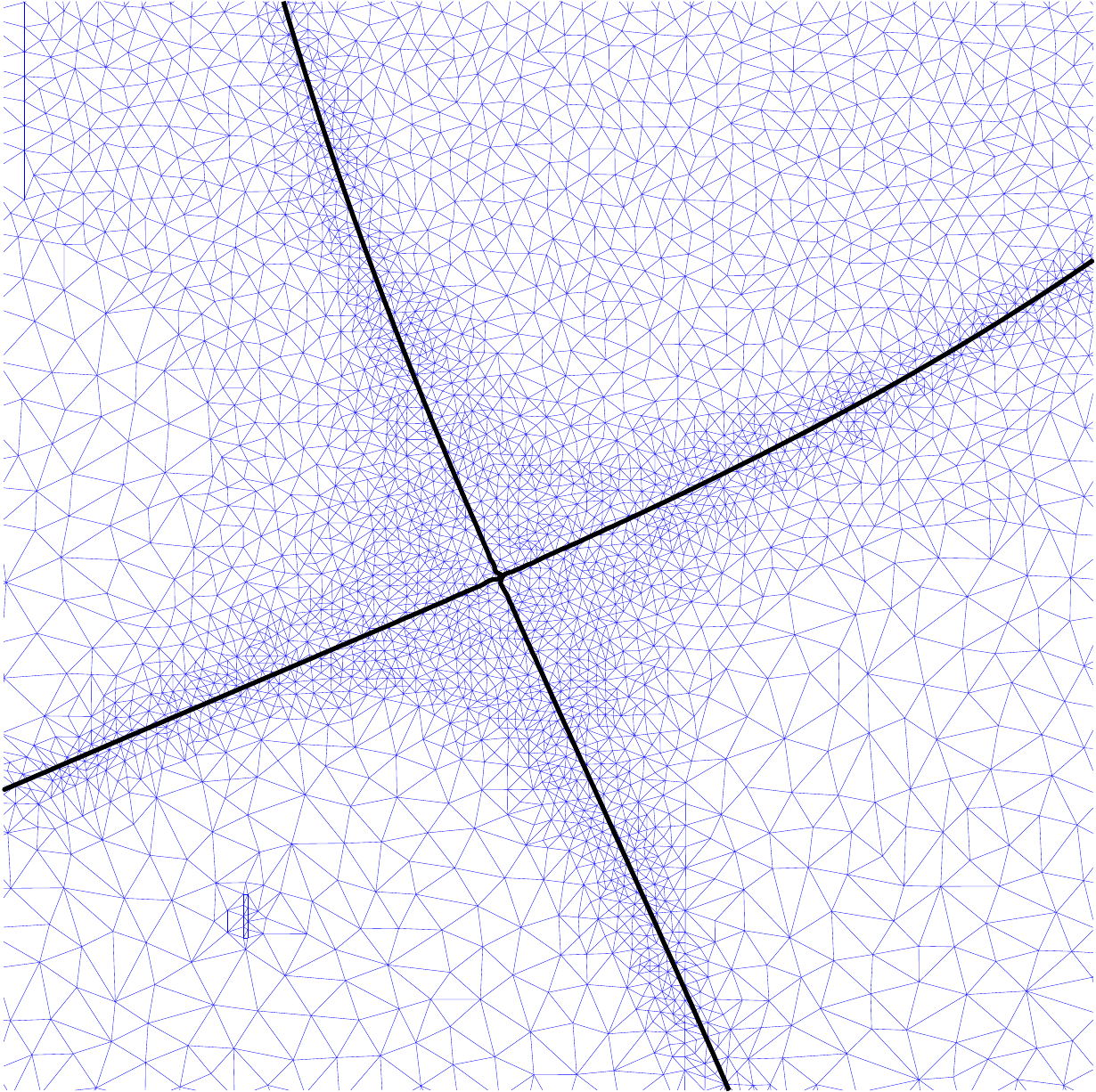}
\end{tabular}
\caption{Top row: The free boundary problem is solved in a sequence of refined meshes, starting from an initial coarse mesh (left). The location of the separatrix drives the refinement strategy and at every level a banded region around it is marked for refinement. Bottom row: Close-up of the neighborhood of the magnetic x-point for the three refinement levels shown in the top.}\label{fig:Adaptive}
\end{figure}

Given that in our application the most relevant quantities are related to the location of the plasma boundary, we resort instead to the following local refinement strategy. The free boundary problem is solved on a coarse grid and the approximate location of the plasma boundary is determined by first finding the location of either the x-point, in the case where the plasma is not in contact with the reactor, or of the point where the streamlines are tangential to an obstructing structure. Let us denote either of these points by $(x^*,y^*)$. Since the plasma boundary corresponds to the level set of $\psi_h(x^*,y^*)$, the elements in the neighborhood of the separatrix are estimated by fixing $\alpha\in(0,1)$ and identifying all the nodes $(x,y)$ where
\[
|\,\psi_h(x,y) - \psi_h(x^*,y^*)\,| \leq \alpha |\,\psi_h(x^*,y^*)\,|.
\]
All the elements having a vertex that satisfies the above condition are then marked and refined---in our experiments we set $\alpha=0.05$. The new mesh has a band of finer elements in the vicinity of the separatrix. The problem is then solved again on the finer mesh and the process can be restarted and repeated. In our experiments, the adaptive cycle is stopped after two iterations of local refinement. Given that, whenever they are present, the x-point, the strike points or the contact points are located in the vicinity of the separatrix, the local refinement leads to a better estimation of their positions. This is especially true for the location of the x-points, which are bound to element vertices; by enriching the mesh with additional vertices more possible locations become accessible to the x-points. Figure \ref{fig:Adaptive} shows a refinement cycle consisting of an initial coarse grid and two levels of adaptive refinement. With refinement limited to a neighborhood of the separatrix, the growth in system unknowns is curtailed while resolution is improved. If needed, the growth of the number of degrees of freedom can be limited even more by decreasing the value of the marking parameter $\alpha$ as the mesh gets finer---in our experiments, the value of $\alpha =0.05 $ was kept fixed across refinement levels.

Since the free boundary problem must be solved at every level of the adaptive algorithm, the solution obtained on a coarser level is interpolated to the new grid and used as initial guess for the Newton iteration. Moreover, the stopping tolerance for the Newton iterations at coarse refinement levels can be set to be relatively large and reduced dynamically as the mesh is refined.  As stopping criterion for the Newton process, the code {\tt FEEQS.M} measures the magnitude of the nonlinear updates to the solution operator relative to the norm of the matrix arising from the discretization of \eqref{eq:VariationalFormulation}. In our experiments, the tolerance for the adaptive solver with $R$ refinements was determined by the formula $\mathrm{TOL}=10^{-11((i+1)/(R+1))}$, where $i=1,2,\ldots,R$. In the sequel, we will refer to {\tt FEEQS.M} as the \textit{direct solver} and to the numerical solutions provided by this code as \textit{direct solutions}. 

Our ultimate goal is to build a surrogate model based on interpolation between numerical solutions associated with a small number of carefully selected current values (specific details are given in Section \ref{sec:Collocation}). However, since the values of the currents passing through the coils determine the location of the plasma boundary, different realizations will result in different adaptively refined meshes. Therefore, the corresponding numerical approximations will belong to different discrete spaces and the different solutions must be projected into a common space (a common mesh) prior to the construction of the surrogate. In our experiments, all the numerical solutions start from a common coarse grid that is adaptively refined around the boundary of the plasma in each realization. The numerical solutions are then $L^2$-projected into a finer common mesh built by uniformly refining the interior of the reactor from the initial, common, coarse mesh. The surrogate obtained in this way is then, by construction, defined on the finer common mesh.

This step introduces an additional projection error that is initially negligible compared to the surrogate error but can eventually become dominant as the surrogate gets refined. To curtail the projection error, in our experiments the common mesh is refined as many times as the prescribed depth of the adaptive algorithm. It is important to note that no numerical solutions are ever computed on the finer common grid; instead only surrogate evaluations are needed. As we shall show later, the time required for such evaluations is often considerably shorter than the time required for the solution of the free boundary problem, even on the coarsest grid.

\section{Stochastic collocation}\label{sec:Collocation}
The solution $\psi$ depends on the values of the current intensities on the coils, $\{I_i\}$, and we wish to understand the impact of these values on the properties of the solution, such as whether its level sets present closed lines completely contained in the vacuum chamber and whether the plasma boundary is defined by an uninterrupted separatrix.  A straightforward way to address this question is to perform a Monte Carlo simulation: solve the system for a large number of realizations of the parameter set and then compute sample statistics (means and variances) and probability estimates using the sample solutions.  Since the computation of each such solution entails solution of a discrete version of the system \eqref{eq:InteriorExteriorFormulation}, simulation done this way is expensive. 
As an alternative, we consider  use of the \textit{stochastic collocation method} \cite{Barthelmann-Novak-Ritter,Sm:1963} to construct a surrogate approximation $\hat \psi_h$ to the discrete solution $\psi_h$ that is less expensive to use in a simulation.

We give a generic overview of this methodology, generally following \cite{KlBa:2005};
this reference also contains a description of the software package  {\tt spinterp}, which we used for implementation.  
In our tests, we used version 5.1.1 of the code, available from \cite{spinterp-v5.1}. 

Suppose we wish to approximate a smooth function $f$ mapping a parameter space $[0,1]$ to a physical space contained in $\mathbb{R}$. One approach is to interpolate $f$ with a finite number of nodes in the parameter space. 
Let $X=\{x_j \in [0,1], 1\le j \le m\}$ be a set of nodes in $[0,1]$
and let $\{\phi_i\}_{i=1}^{m}$ be interpolating functions satisfying
\[
\phi_i(x_j)=
\left\{\begin{array}{lll}
1 & \text{ if } & i=j\\
0& \text{ if } & i\ne j
\end{array}.
\right.
\]
Letting $U^{X}$ be the one-dimensional interpolation operator associated with the set of nodes $X$, the interpolation formula for $f$ is
\begin{equation}
\label{eq: uni_interp}
U^{X}(f) = \sum_{j=1}^m f(x_j) \, \phi_j\,. 
\end{equation}
A $d$-dimensional multivariate interpolant can be defined in an analogous manner on $[0,1]^d$ as
\begin{equation}
\label{eq: Multi-interp_formula}
\mathcal U (f) = \left(U^{X_1}\otimes\cdots\otimes U^{X_d}\right) (f)
: = 
\sum_{j_1 =1}^{m_1}   \cdots  \sum_{j_d=1}^{m_d}
f\left(x_{j_1}^{(1)},\ldots,x_{j_d}^{(d)}\right)\, \left(\phi_{j_1}^{(1)}\otimes\cdots \otimes  \phi_{j_d}^{(d)}\right)
= \sum_{j=1}^M f({\bf x}_j)\, \Phi_j\,,
\end{equation}
where for every index $p=1,\ldots,d$, $X_p$ denotes the set of $m_p$ nodes along the $p$-th component, the superscripts tag the corresponding component in the $d$-dimensional space, and the interpolatory functions 
$\{\Phi_j\}_{j=1}^{M}$ have tensor product structure. Due to the tensor product construction, this formulation requires a total of $M=\Pi_{i=1}^d m_i$ interpolation nodes. We will refer to this set of support nodes as the \textit{full grid}. In the shorthand expression on the rightmost of \eqref{eq: Multi-interp_formula}, the sum ranges over each and every one of nodes on the full grid. It is easy to see that the size of the full grid grows rapidly as the dimension increases; this phenomenon is known as the \textit{curse of dimensionality}.

\begin{figure}
\begin{center}
\vspace{-.3in}
\includegraphics[width=0.5\linewidth]{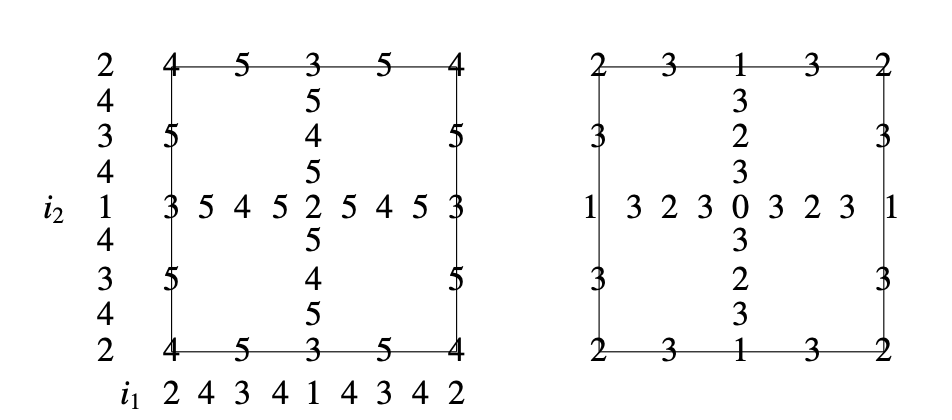}
\vspace{-.2in}
\end{center}
\caption{Left: multi-index values $(i_1,i_2)$ and contributions to $H_{q,d}$ from
$\bigcup_{|\boldsymbol{i}|=q} {\bf X}^{\boldsymbol{i}}$, for $d=2$ and $q=i_1+i_2=2$ through $5$.  
Right: Highest grid level associated with each node.}\label{fig:sparse-grid-picture}
\end{figure}

Sparse grid methods are designed to alleviate the curse of dimensionality while maintaining  accuracy comparable to that of the full grid. For $1\leq p\leq d$, let $X^{i_p}$ be a collection of nodes in $[0,1]$ of size $\xi(i_p)$ where 
$\xi(i_p)$ is an increasing function that relates the index $\boldsymbol{i}_p$ to the 
number of support nodes in component $p$.
Let $\boldsymbol{i} = (i_1,\ldots,i_d)\in \mathbb{N}_+^d$ be the $d$-dimensional multivariate index, 
let ${\bf X}^{\boldsymbol{i}}:= \left(X^{i_1}\times \cdots \times X^{i_d}\right) \subset [0,1]^d$ denote the 
associated collection of points in $[0,1]^d$, 
and for integer $q\ge d$, let 
\begin{equation}
\label{eq:NestedColPts}
H_{q,d} :=  \bigcup_{|\boldsymbol{i}|\le q} {\bf X}^{\boldsymbol{i}} = 
\bigcup_{|\boldsymbol{i}|\le q} \left(X^{i_1}\times \cdots\times X^{i_d}\right). 
\end{equation}
Insight into the structure of $H_{q,d}$ can be obtained from Figure \ref{fig:sparse-grid-picture}, which shows 
on the left examples of the individual contributions 
$\bigcup_{|\boldsymbol{i}|=q}  {\bf X}^{\boldsymbol{i}}$ to $H_{q,d}$ for $d=2$, and
$q=2$ through $5$.
In this image, the individual indices come from the values $i_p=1$ through $4$ and $\xi(1)=1$, 
$\xi(i_p)=2^{i_p-1}+1$ for $i_p>1$.  

We also refer to the \textit{level} of a sparse grid, where nodes at level $\ell$ are those with indices
$\boldsymbol{i}$ such that $|\boldsymbol{i}|-d\le \ell$.
The right side of Figure \ref{fig:sparse-grid-picture} shows the highest grid level for each node.
For example, the level-2 grid consists of nodes with labels $0$, $1$ and $2$ on the right side of the figure. The sparse grid interpolant reproduces $f$ at all nodes in the sparse grid. When it is generated using polynomials, it includes univariate polynomial factors of degree $\xi({i}_p)-1$; for example, there is a univariate polynomial of degree $\xi(4)-1=8$ determined from points with multi-index $(4,1)$ ($q=5$).

As is well known  \cite{Barthelmann-Novak-Ritter,KlBa:2005,Sm:1963}, if $q$ is small or 
$d$ is large, then $H_{q,d}$ will be much smaller than the size of the full grid.
For example, we will be interested in $d=12$, the number of coils in the reactor model we will study, which will produce sparse grids of size $2,649$ for level $3$ and $17,265$ for level $4$; these contrast with corresponding full grid sizes of sizes $2.8\times 10^{11}$ and $5.8\times 10^{14}$ for the corresponding full grids.

Sparse grid collocation uses a sparse version of the interpolation formula.
The grids comprising $H_{q,d}$ can be either nested or not.
We will use the Smolyak construction \cite{Sm:1963} as implemented in {\tt spinterp} \cite{KlBa:2005}, which uses nested grids, i.e, they satisfy in each component $X^{i_p} \subset X^{i_p+1}$. Let $\Delta H_{q,d}:=H_{q,d}\setminus H_{q-1,d}$, which contains all nodes ${\bf X}^{\boldsymbol{i}}$ with $|\boldsymbol{i}|=q$. It follows from the nestedness of the grids that \eqref{eq:NestedColPts} can be rewritten as 
\begin{equation} \label{eq:nested-grid}
H_{q,d} = H_{q-1,d}\cup \Delta H_{q,d}.
\end{equation}
Let $\Delta^{i_p} = \mathcal{U}^{i_p} - \mathcal{U}^{i_p - 1}$  represent an incremental 
operator associated with a one-dimensional nested grid, where $\mathcal{U}^0=0$. The Smolyak interpolation operator $\mathcal{S}_{q, d}$ can be defined by the formula
\begin{equation}
\label{eq: Smolyak_Interp_formula}
 \mathcal{S}_{q, d} (f): =\sum_{|\boldsymbol{i}|\le q}
 \left(\Delta^{i_1}\otimes \cdots\otimes \Delta^{i_d}\right)(f).
\end{equation}
This can be rewritten as
\begin{equation} \label{eq: SmolyakAlgorithm}
\mathcal{S}_{q, d} (f) 
= \sum_{|\boldsymbol{i}|\le q-1}\left(\Delta^{i_1}\otimes \cdots \otimes \Delta^{i_d}\right)(f)+ \sum_{|\boldsymbol{i}| = q} \left(\Delta^{i_1}\otimes \cdots\otimes \Delta^{i_d}\right)(f)
\ = \ \mathcal{S}_{q-1, d}(f) + \Delta  \mathcal{S}_{q, d}(f).
\end{equation}
This Smolyak construction uses the relation \eqref{eq:nested-grid} by building the interpolant on $H_{q,d}$ as the sum of the interpolant on $H_{q-1,d}$ (i.e., $\mathcal{S}_{q-1,d}(f)$) and a correction $\Delta \mathcal{S}_{q, d}(f)$ that interpolates $f({\bf x})-\mathcal{S}_{q-1,d}(f)({\bf x})$ at the nodes in $\Delta H_{q,d}$ and is $0$ on $H_{q-1,d}$. We elaborate on this using the left image of Figure \ref{fig:sparse-grid-picture}. Here, $q=5$ (and $d=2$). $\mathcal{S}_{4,2}(f)$ consists of polynomials that interpolate $f$ at all the nodes in the image with labels $4$ or smaller ($H_{4,2}$), and $\Delta  \mathcal{S}_{5, 2}(f)$ is a polynomial that has the value $0$ at all these nodes and value $f({\bf x})-\mathcal{S}_{4,2}(f)({\bf x})$ at all nodes with label $5$ ($\Delta H_{5,2}$). Evaluation of $\mathcal{S}_{5,2}(f)$ at any other point $\bf x$ in the domain is done recursively as $\mathcal{S}_{5,2}(f)({\bf x}) = \mathcal{S}_{4,2}(f)({\bf x}) + \Delta \mathcal{S}_{5,2}(f)({\bf x})$. There are $16$ terms in the sum contributing to $\Delta  \mathcal{S}_{5,2}(f)({\bf x})$, one for each of the nodes labeled $5$. Consider for example the bottom left of these nodes situated on the left boundary of the domain, which we refer to here as ${\bf x}_\ast=(x_\ast,y_\ast)$. The summand associated with ${\bf x}_\ast$ contributing to $\Delta  \mathcal{S}_{5,2}(f)$ has the form of a product $(f({\bf x}_\ast)-\mathcal{S}_{4,2}({\bf x}_\ast))\phi_1(x)\phi_2(y)$, where $\phi_1$ is a Lagrange polynomial of degree $2$ with value $1$ at $x_\ast$ and $0$ at the two other nodes for which $y=y_\ast$, and similarly $\phi_2$ is a degree $4$ polynomial with value $1$ at $y_\ast$ and $0$ at the other nodes along the vertical line. Use of this hierarchical structure allows the interpolant to be implemented using relatively simple data structures and bookkeeping while also enabling the re-use of all the information present in $H_{q,d}$ if a higher order interpolant $H_{q+i,d+j}$ is needed at a later time, thus speeding up the off-line component of the process; see \cite{KlBa:2005} for details.

In computational tests with stochastic collocation, presented in the next section, when $m_p=\xi(i_p)$  interpolation points are used in a component, the basis functions are built from interpolating polynomials of degree $m_p-1$, and interpolation points consist of the extrema of the Chebyshev polynomials of degree $m_p$, $x_j =\cos \left( \frac{(j-1)\pi}{m_p-1} \right) ,\, 1 \le j \le m_p$. With this strategy, the interpolation error is bounded as \cite[Theorem~8]{Barthelmann-Novak-Ritter}
\begin{equation} \label{eq:coll-error-bound}
  \big\|f-\mathcal{S}_{q, d} (f)\big\|_\infty = O\left(N^{-k} |\log N|^{(k+2)(d-1)+1}\right) 
\end{equation}
where $N$ is the number of nodes in the sparse grid and $k$, a measure of smoothness
of $f$, is the largest integer such that $D^{\beta}f$ is continuous if $\beta_i\le k$ for all $i=1,\ldots,d$. This bound shows that accuracy is good for large $N$ but there is a logarithmic factor that may limit its utility, as the number of nodes must be at least $N \sim e^{\left((k+2)(d-1)+1\right)/k}$ before the logarithmic growth is controlled. This makes for a long pre-asymptotic regime that becomes more acute as the number of parametric dimensions $d$ grows---an example of this can be seen in the right panel of Figure \ref{fig:GeometricConfiguration}.

We comment on one aspect of our particular application that is affected by the use of the surrogate based on collocation. A key step in the computations is the determination of the location of the x-point. In the absence of obstructing structures, this is the point where the streamlines of the magnetic field cross themselves, and therefore determine the separatrix between the region of the plane where the lines are closed (allowing for confinement) and the region where they are open. As such, the segment of the separatrix describing a closed curve determines the plasma boundary $\Gamma$. Geometrically, the x-point corresponds to a saddle point on the graph of the scalar field $\psi$; this characteristic is used to determine its location on the discrete approximations $\psi_h$ coming from the direct solver {\tt FEEQS.m} and $\widehat{\psi}_h$ coming from the evaluation of the surrogate. By definition, a saddle point must be a local maximum in one direction, and a local minimum in another one.

\begin{figure}[htb]\centering
\includegraphics[width=0.25\linewidth]{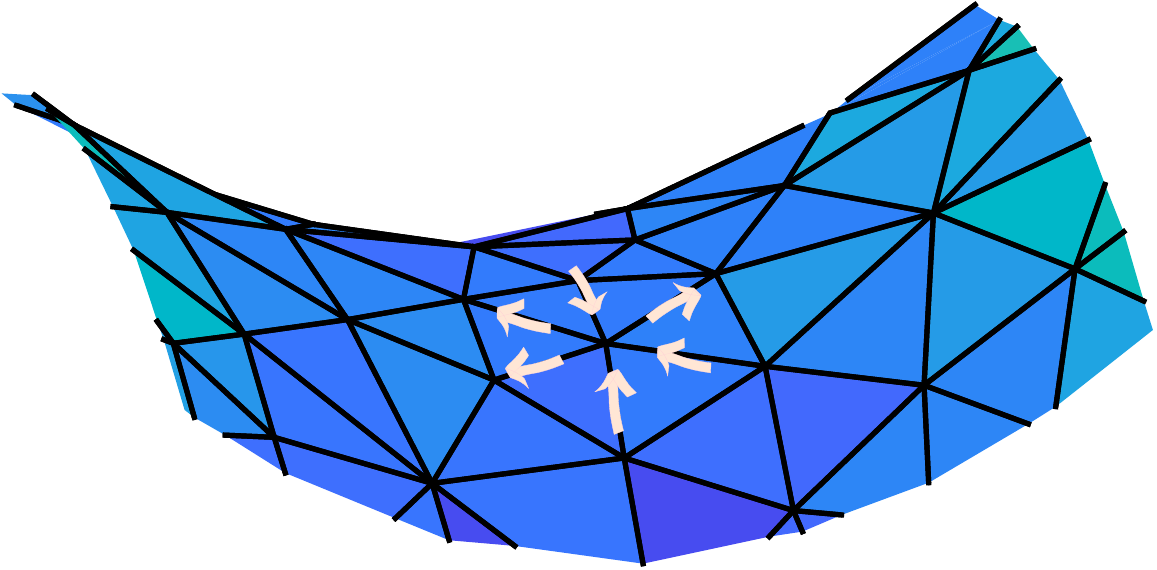}
\qquad \qquad \qquad
\includegraphics[width=0.2\linewidth]{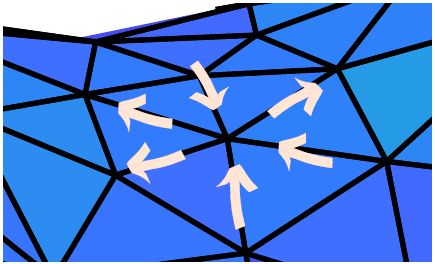}
\caption{In a piecewise linear approximation of $\psi$, the saddle points are located at those element vertices $\mathbf v$ for which the difference $\psi_h(\mathbf v)-\psi_h(\mathbf v_i)$ for all its neighboring vertices $\mathbf v_i$ changes sign at least four times. The sign of this difference indicates whether the discrete function $\psi_h$ along the element vertices is increasing or decreasing. The point surrounded by the arrows in the figure corresponds to a saddle, as it is a local maximum in one direction and a local minimum in another one. This is reflected by the number of times the arrows change direction.}\label{fig:DiscreteSaddlePoint}
\end{figure}

A discrete version of this fact is implemented in {\tt FEEQS.m} and used as the criterion to identify candidates for the location of the x-point. An approximation to the directional derivative of the piecewise linear approximation $\psi_h$ in the direction of each of the element edges originating at an element vertex $\mathbf v$ is given by the difference  $\psi_h(\mathbf v)-\psi_h(\mathbf v_i)$, where the subindex $i$ on $\mathbf v_i$ runs over all the mesh nodes connected to $\mathbf v$ by an edge. At a saddle point of a piecewise linear approximation, the directional derivative must change sign at least four times---this situation is depicted in Figure \ref{fig:DiscreteSaddlePoint}. The code computes the number of sign changes of the directional derivative for every node within the reactor and identifies those vertices where the sign changes at least four times as candidates for the saddle point. For numerical solutions coming out of the direct solver this criterion is typically enough to identify the discrete x-point correctly. However, due to the oscillatory nature of polynomial interpolation, evaluations of the surrogate require some additional work. Using a classical gradient recovery technique introduced by Zienkiewicz and Zhu \cite{ZiZh:1992a, ZiZh:1992b}, it is possible to build a second order approximation to $\nabla\psi$ from the nodal values of $\widehat{\psi}_h$ at each vertex and its neighbors. If more than one node satisfies the ``alternating sign" criterion, then the magnitude of the recovered gradient is evaluated at each candidate node and all of its neighbors are averaged. The node with the smallest local average magnitude is then identified as the saddle point. The logic behind this strengthened criterion is that the spurious oscillations originating from the surrogate evaluations are highly localized in nature, and therefore must present larger gradients than the ones associated with the underlying function $\psi$.   

\section{Numerical experiments}\label{sec:NumericalExperiments}

\begin{figure}\centering
\includegraphics[height=0.28\linewidth]{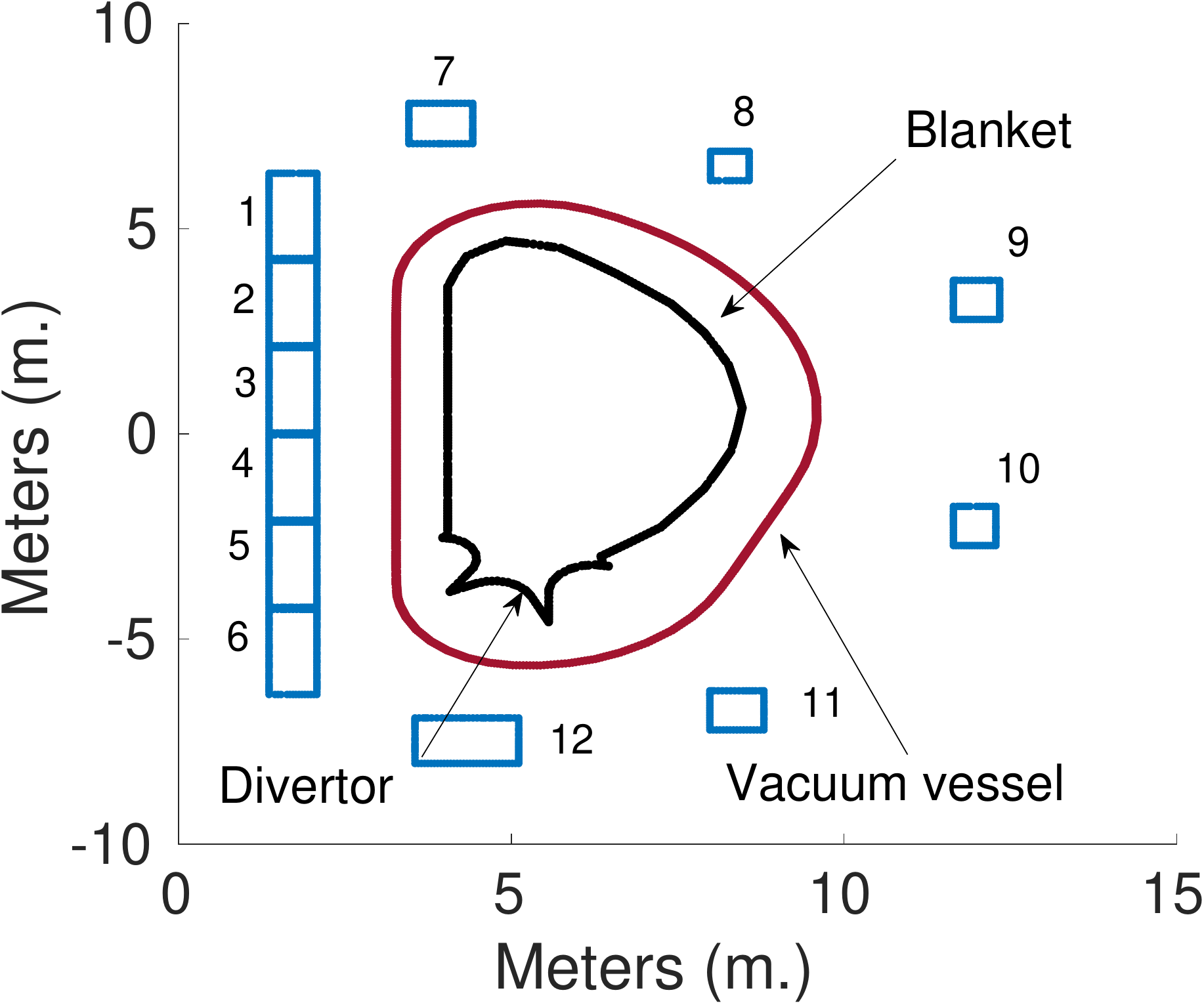} \qquad \qquad \qquad
\includegraphics[height=0.28\linewidth]{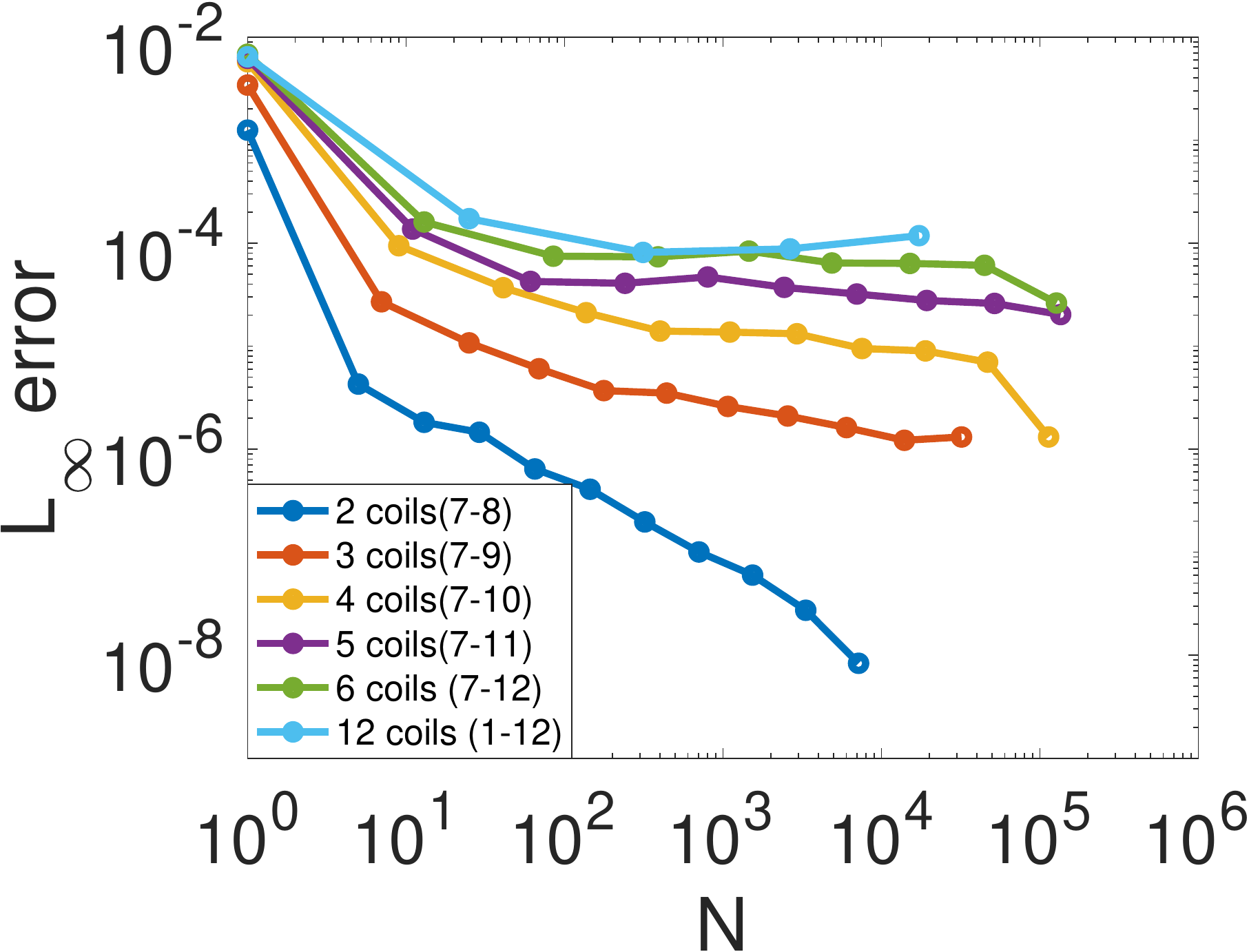} 
\caption{Left: Reactor geometry used for the numerical simulations. The divertor, inner wall (blanket), and outer wall of the vacuum vessel are shown along with the numbered coils. Right: Mean relative difference between surrogate evaluations and numerical solutions to the free boundary problem of 100 randomly chosen current values, with 1\% variability in currents. The number of collocation nodes ($N$) on the sparse grid is displayed on the horizontal axis. Results of the experiment for a variety of numbers of currents perturbed are displayed in different colors.}\label{fig:GeometricConfiguration}
\end{figure}

In this section, we present the results of numerical experiments that illustrate the
accuracy and cost savings obtained using the surrogate in place of direct solution of the 
discrete version of \eqref{eq:FreeBoundary}.
The simulations were done using the ITER-like geometry depicted in Figure 
 \ref{fig:Adaptive} and the left panel of Figure \ref{fig:GeometricConfiguration}, where the latter 
 figure shows the numbering of the coils; each component of the vector 
 $\mathbf I = [I_1,\ldots,I_{12}]$ corresponds to the current value going through the coil with the same number. We will consider the impact of a parametrized version of (\ref{eq:FreeBoundary}) starting from the deterministic equlibrium generated by the current values
\begin{equation}\label{eq:CentralCurrentValue}
{\renewcommand{\arraycolsep}{2pt}
\begin{array}{llll}
I_1 = -1.4 \times 10^{6} A, \quad & I_2 = -9.5 \times 10^{6} A, \quad & I_3 = -2.0388 \times 10^{7} A, \quad & I_4 = -2.0388 \times 10^{7} A, \\
I_5 = -9 \times 10^{6} A, \quad & I_6 = 3.564 \times 10^{6} A, \quad & I_7 = 5.469 \times 10^{6} A, \quad & I_8 = -2.266 \times 10^{6} A, \\
I_9 = -6.426 \times 10^{6} A, \quad & I_{10} = -4.82 \times 10^{6} A, \quad & I_{11} = -7.504 \times 10^{6} A, \quad & I_{12} = 1.724 \times 10^{7} A, 
\end{array}
}
\end{equation}
studied by Faugeras and Heumann \cite{FaHe:2017}. The parametrized  problem is constructed assuming that a uniformly distributed variability is present on the current flowing through each coil. For example, 1\% variability yields a 12-dimensional parameter space 
centered at the values $(I_1,\ldots,I_{12})$ of (\ref{eq:CentralCurrentValue}) that can be described as the Cartesian product of the closed intervals $(I_1\pm 0.01|I_1|)\ldots (I_{12}\pm 0.01|I_{12}|)$. Surrogate functions for the parameter-dependent discrete stream function $\psi_h$ were built using increasing levels of sparse grids as described in Section \ref{sec:Collocation}. All experiments were performed using {\tt MATLAB} R2020b on a System 76 Thelio Major with 256GB RAM and a 64-Core @4.6 GHz AMD Threadripper 3 Processor. 

We will first assess the approximation properties of the surrogate function by comparing the point-wise difference between direct solutions and evaluations of surrogates built with increasing levels of refinement in the sparse grid. We will then test the surrogate's accuracy of reproducing physically relevant quantities such as locations of x-points or plasma shaping parameters, as well as its ability to predict unwanted interactions between the plasma and components of the reactor. These comparisons are done for identical sets of randomly chosen current values. After establishing the approximation capabilities of the surrogate function, we will explore its efficiency by comparing the timings associated with surrogate evaluations and direct solutions. This will require us to introduce a dynamic sampling strategy designed to avoid unnecessarily large sample sizes, thus testing each of the two methods under ideal circumstances. Finally, we conclude with experiments that showcase the information that can be extracted from fast simulations made using the surrogate function.

We start with several tests to assess the accuracy of the surrogate approximation to $\psi_h$. To begin, we evaluated the surrogates at $n_s=100$ random sample values of the currents $\{{\mathbf I}_i\}_{i=1}^{n_s}$, 
solved the free boundary problem at the same sample currents, and then we
assess the error using the average of the relative error in the $L_{\infty}$ norm
\begin{equation}\label{eq:MeanError}
E_{n_s}: = \frac{1}{n_s}\sum_{i=1}^{n_s}
\frac{\|\widehat{\psi}_h(\mathbf I_i) - \psi_h(\mathbf I_i)\|_{\infty}} {\|\psi_h(\mathbf I_i)\|_{\infty}}
\end{equation}
For this test, we considered perturbed versions of all twelve currents, as well as perturbations of subsets of them of sizes $2$ through $6$. The smaller subsets comprise parameter sets of smaller dimension and enable the use of higher sparse-grid  levels. The results of these tests for 1\% variability are shown in the right panel of Figure \ref{fig:GeometricConfiguration} and Table \ref{Tab:Conv_test1}. They indicate some benefit of increasing the sparse grid level although these improvements are not monotone and they are not dramatic for levels greater than four. They also show that accuracy tends to be higher for smaller numbers of parameters; we attribute this to the logarithmic factor in the error bound (\ref{eq:coll-error-bound}). To single out the approximation properties of the surrogate and separate them from the projection error mentioned in Section \ref{sec:FEM}, these convergence studies were performed using only a common grid without any adaptive refinement around the plasma boundary.

\begin{table}[ht]
	\centering
	\scalebox{0.7}{
		\begin{tabular}{|c|c|c|c|c|c|c|c|c|c|c|c|c|c|}
		    \hline
			&\multicolumn{2}{c|}{Coils 7-8}&\multicolumn{2}{c|}{Coils 7-9}&\multicolumn{2}{c|}{Coils 7-10}&\multicolumn{2}{c|}{Coils 7-11}&\multicolumn{2}{c|}{Coils 7-12}&\multicolumn{2}{c|}{Coils 1-12}\\
			\hline	
			Level &\# Nodes&$L_{\infty}$ error&\# Nodes&$L_{\infty}$ error&\# Nodes&$L_{\infty}$ error&\# Nodes&$L_{\infty}$ error&\# Nodes&$L_{\infty}$ error&\# Nodes&$L_{\infty}$ error\\
			\hline
			 0&1&1.2f8e-03&1&3.42e-03&1&5.76e-03&1&6.22e-03&1&6.82e-03&1&6.44e-03\\
			 1&5&4.28e-06&7&2.69e-05&9&9.43e-05&11&1.37e-04&13&1.60e-04&25&1.73e-04\\
			 2&13&1.82e-06&25&1.08e-05&41&3.69e-05&61&4.25e-05&85&7.46e-05&313&8.13e-05\\
			 3&29&1.46e-06&69&6.04e-06&137&2.11e-05&241&4.07e-05&389&7.36e-05&2649&8.79e-05\\
			 4&65&6.42e-07&177&3.70e-06&401&1.40e-05&801&4.69e-05&1457&8.36e-05&17265&1.18e-04\\
			 5&145&4.08e-07&441&3.50e-06&1105&1.37e-05&2433&3.72e-05&4865&6.41e-05 &93489  &-\\
			 6&321&1.97e-07&1073&2.60e-06&2929&1.32e-05&6993&3.20e-05&15121&6.35e-05 &442001  &-\\
			 7&705&1.01e-07&2561&2.10e-06&7537&9.48e-06&19313&2.77e-05&44689&6.09e-05 &1887377  &-\\
			 8&1537&5.99e-08&6017&1.62e-06&18945&8.99e-06&51713&2.60e-05&127105&2.62e-05 &7451393  &-\\
			 9&3329&2.73e-08&13953&1.22e-06&46721&6.99e-06&135073&2.03e-05&350657 & - &27649409  &-\\
			 10&7169&8.35e-09&32001&1.31e-06&113409&1.31e-06&345665&-&943553 & -  &97566977  &-\\
			\hline	
	\end{tabular}}
			\caption{Mean relative $L_{\infty}$ errors (\ref{eq:MeanError}) of surrogate solutions for 100 random sample currents, various sparse-grid levels and numbers of perturbed currents, with 1\% noise level.}
		\label{Tab:Conv_test1}
\end{table}

\begin{figure}[htb]\centering
\begin{tabular}{cccc}
\includegraphics[width=0.2\linewidth]{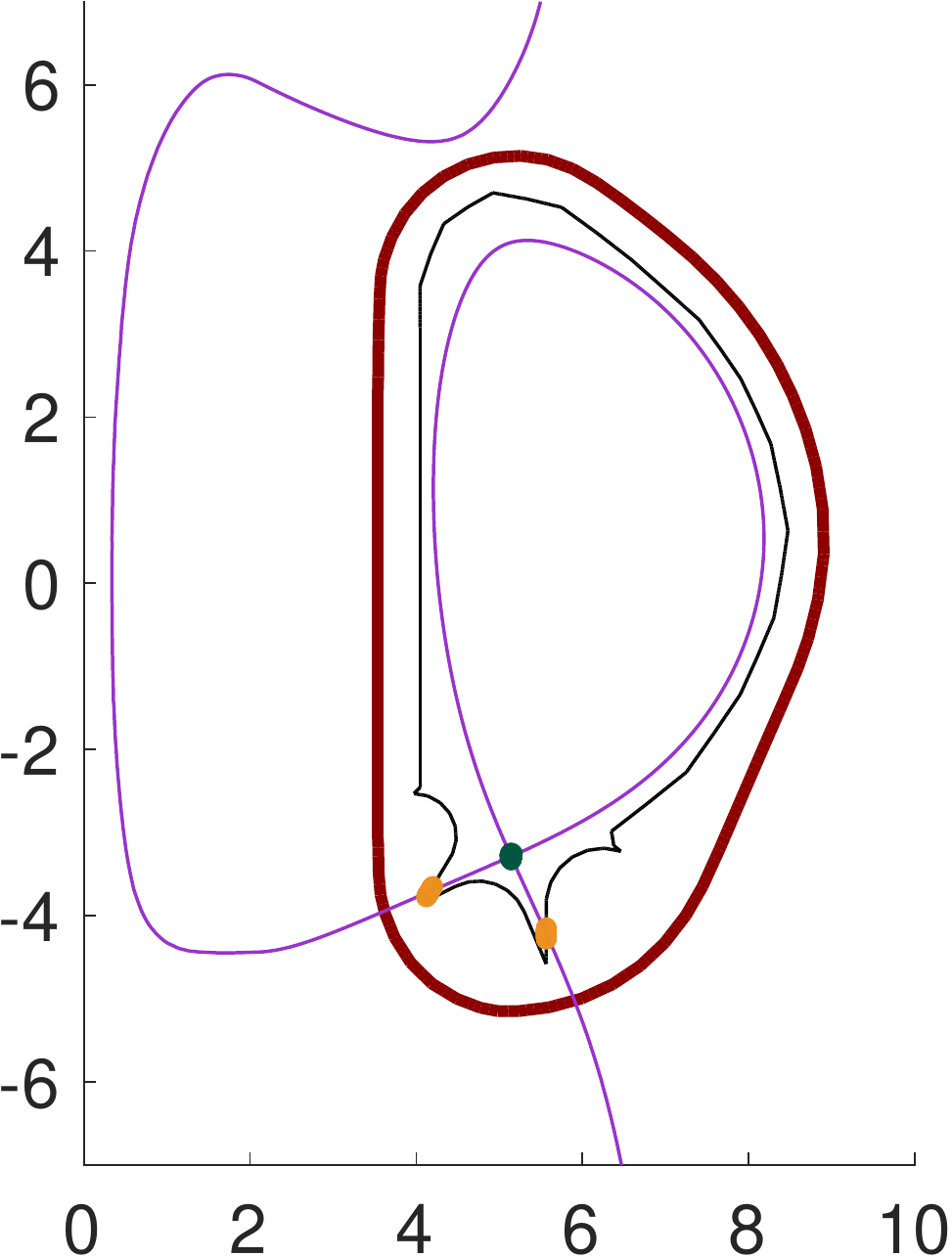} \quad&
\includegraphics[width=0.2\linewidth]{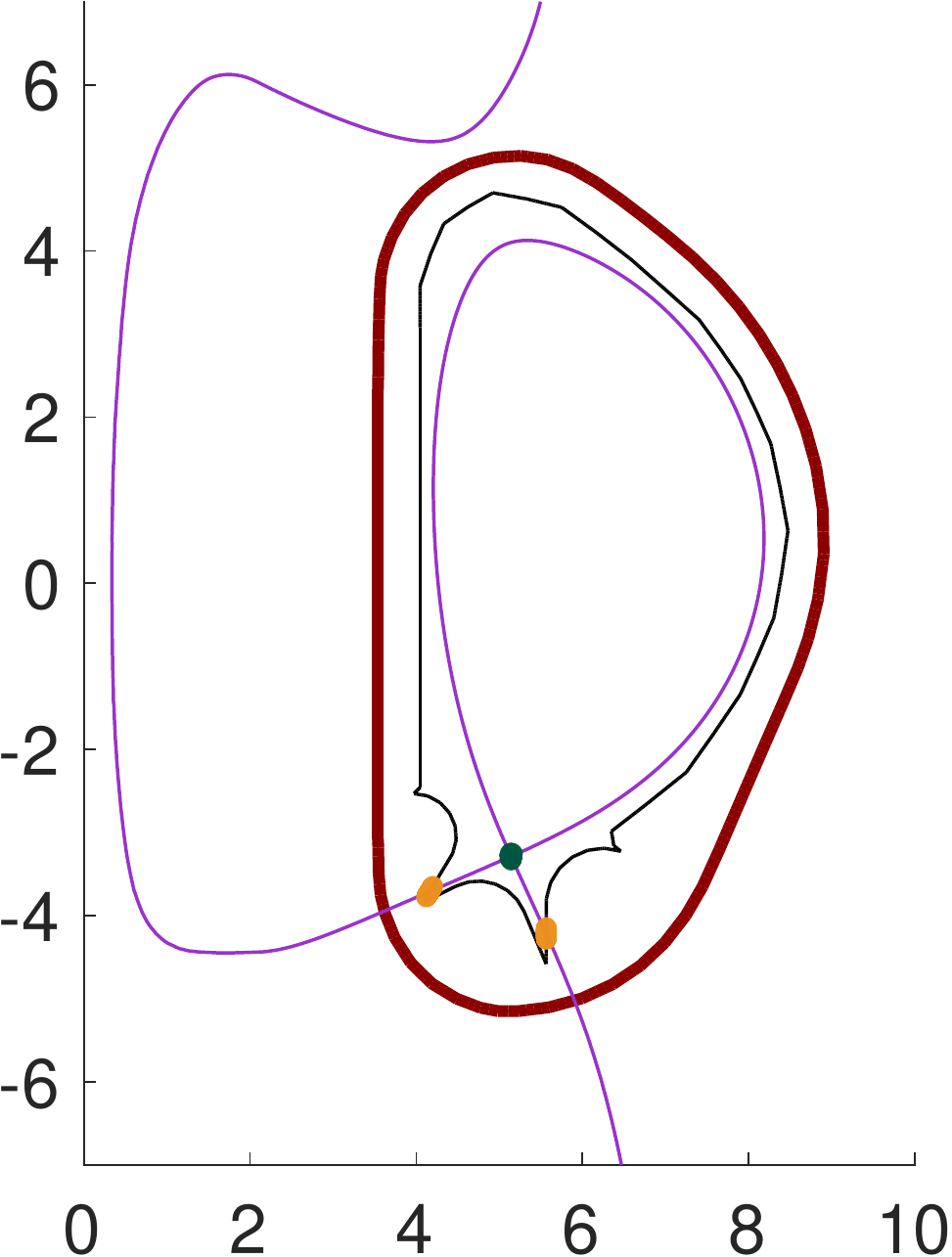} \quad&
\includegraphics[width=0.2\linewidth]{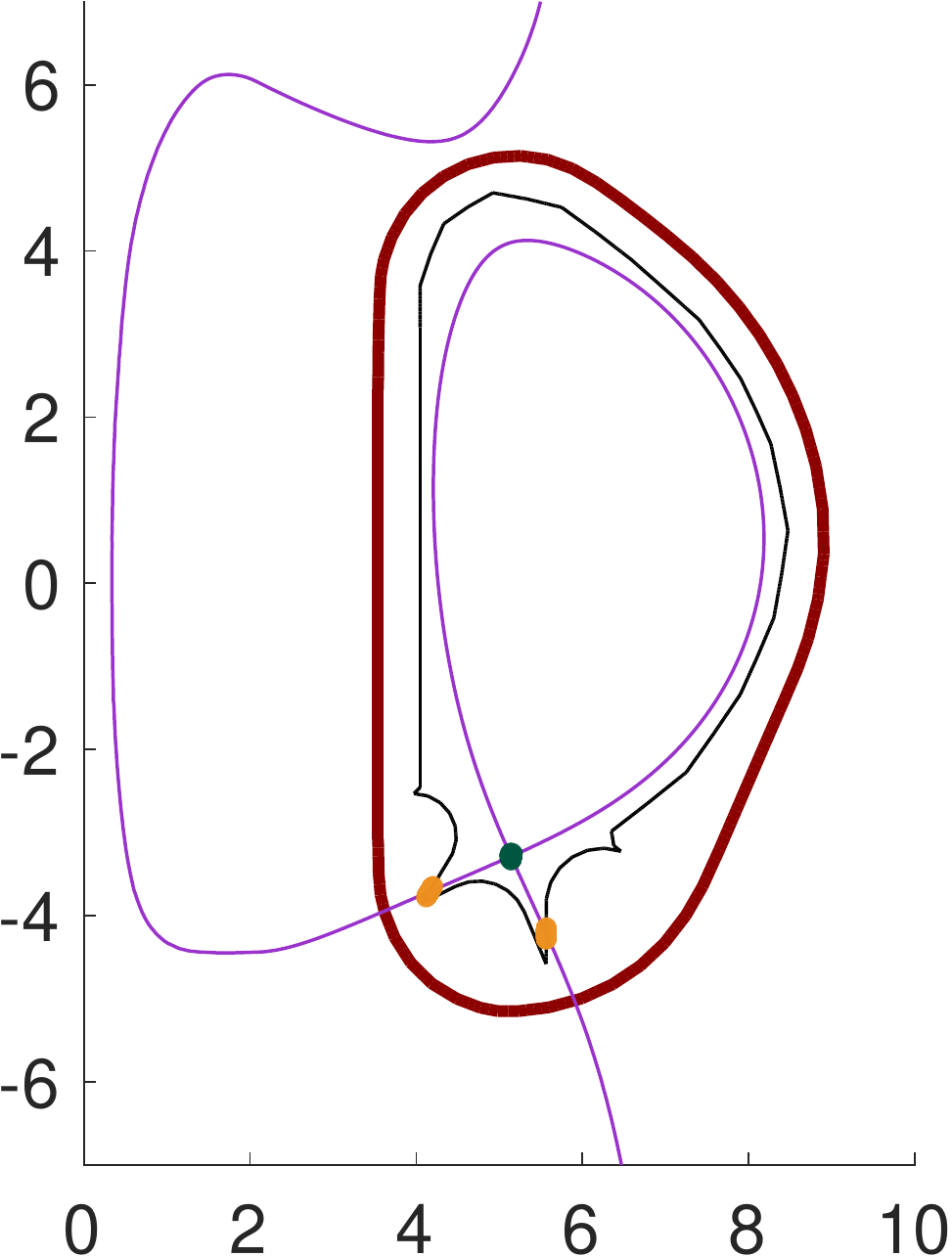} \quad&
\includegraphics[width=0.2\linewidth]{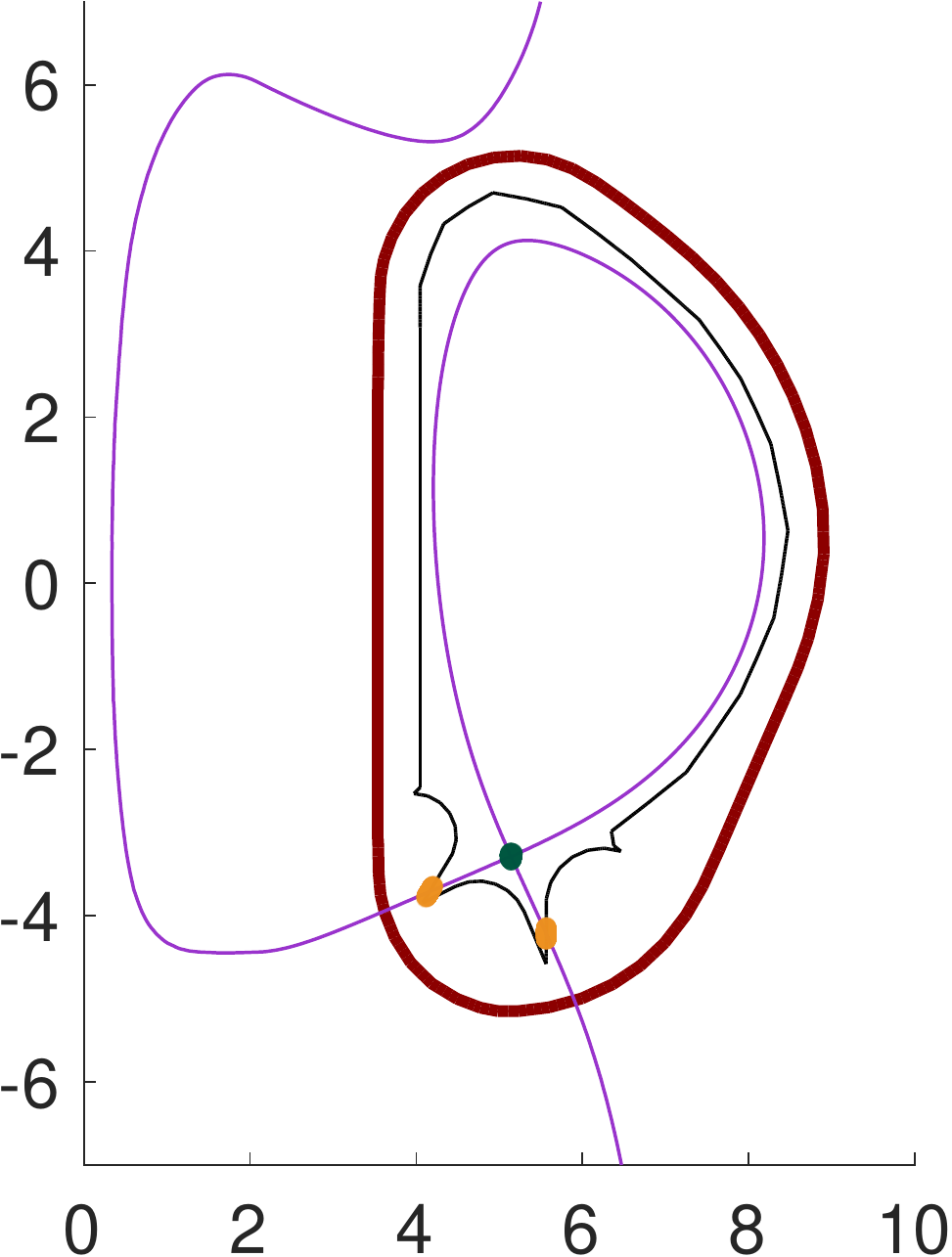}\\ 
\includegraphics[width=0.2\linewidth]{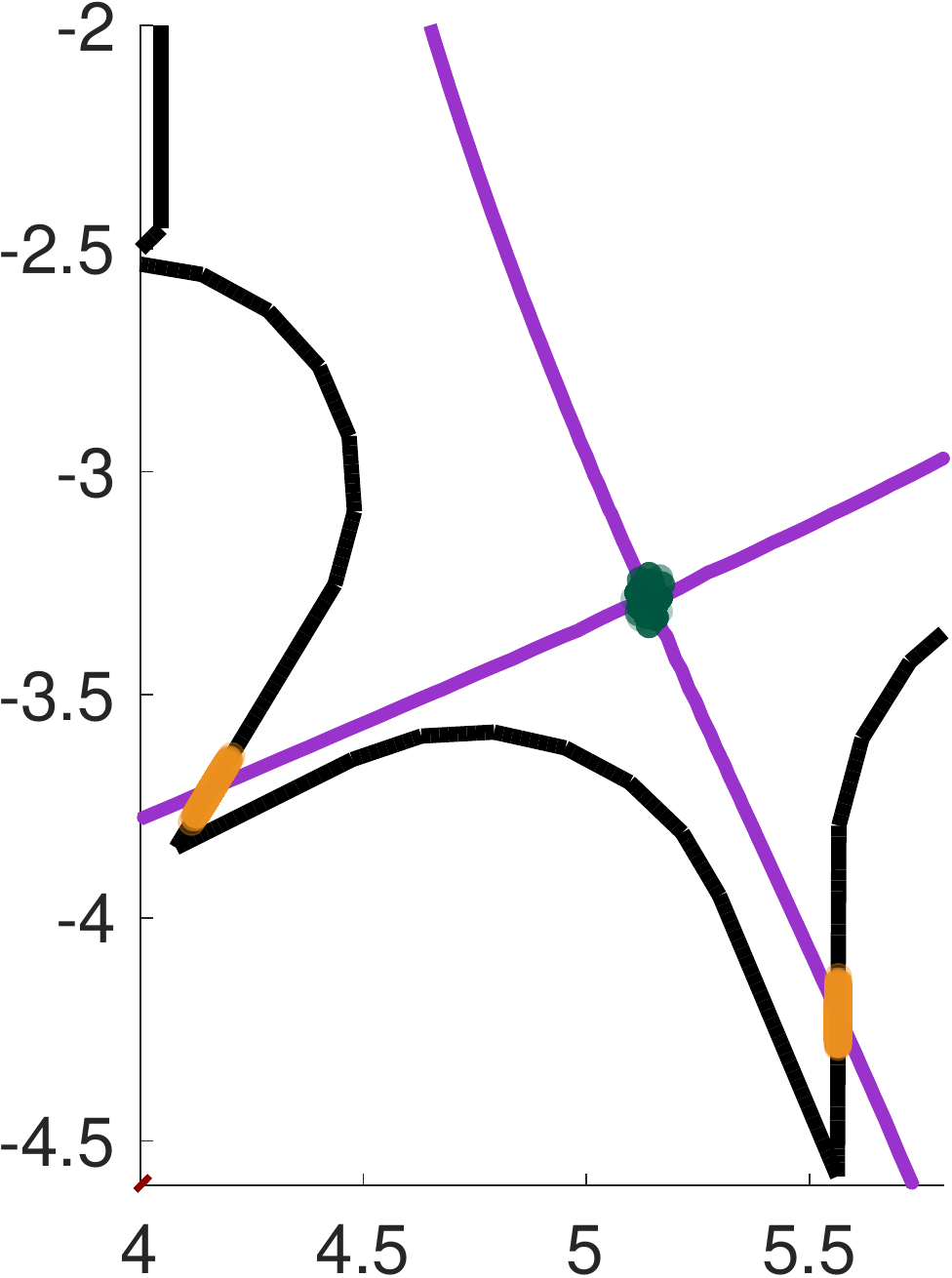} \quad&
\includegraphics[width=0.2\linewidth]{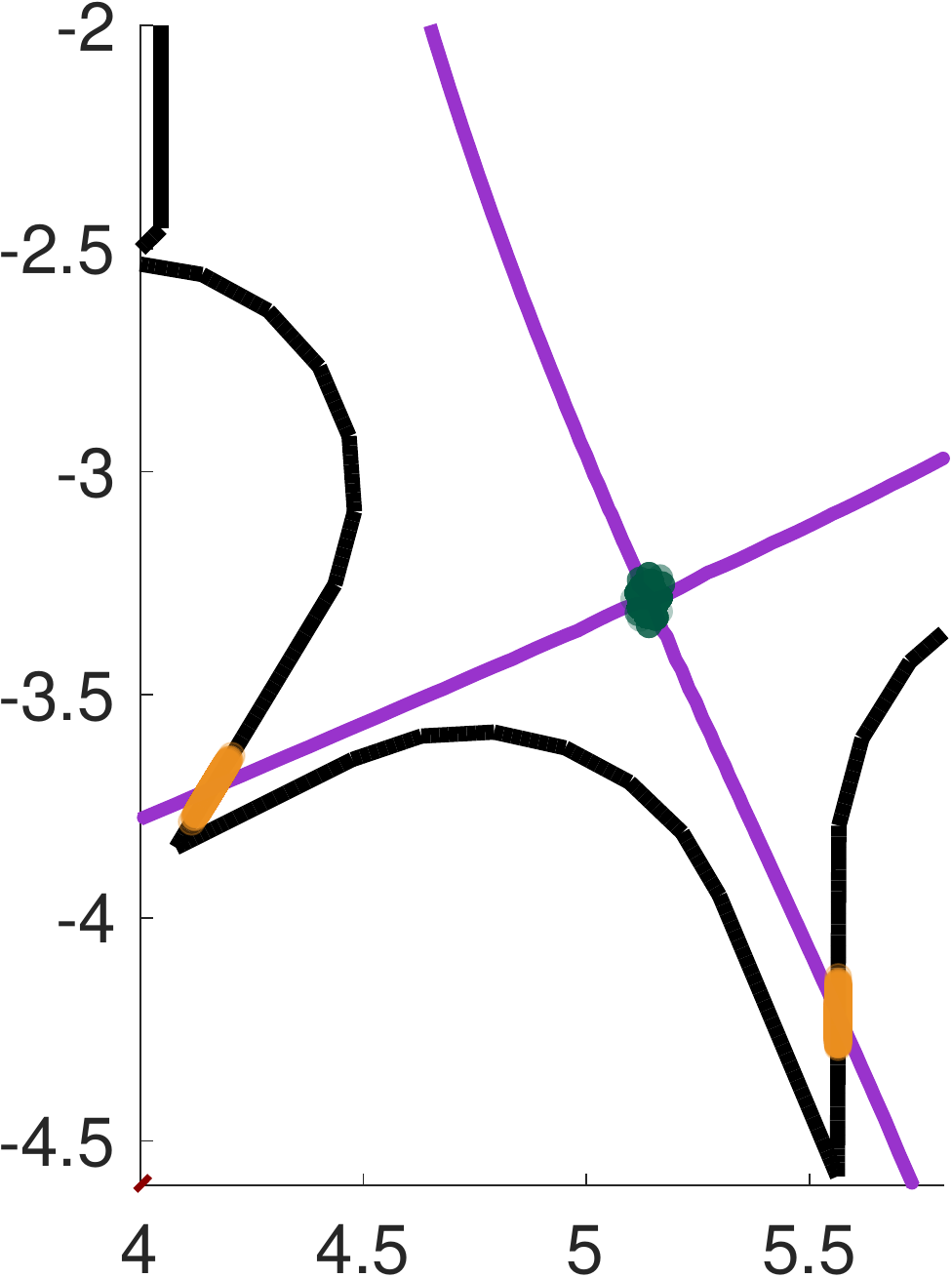} \quad&
\includegraphics[width=0.2\linewidth]{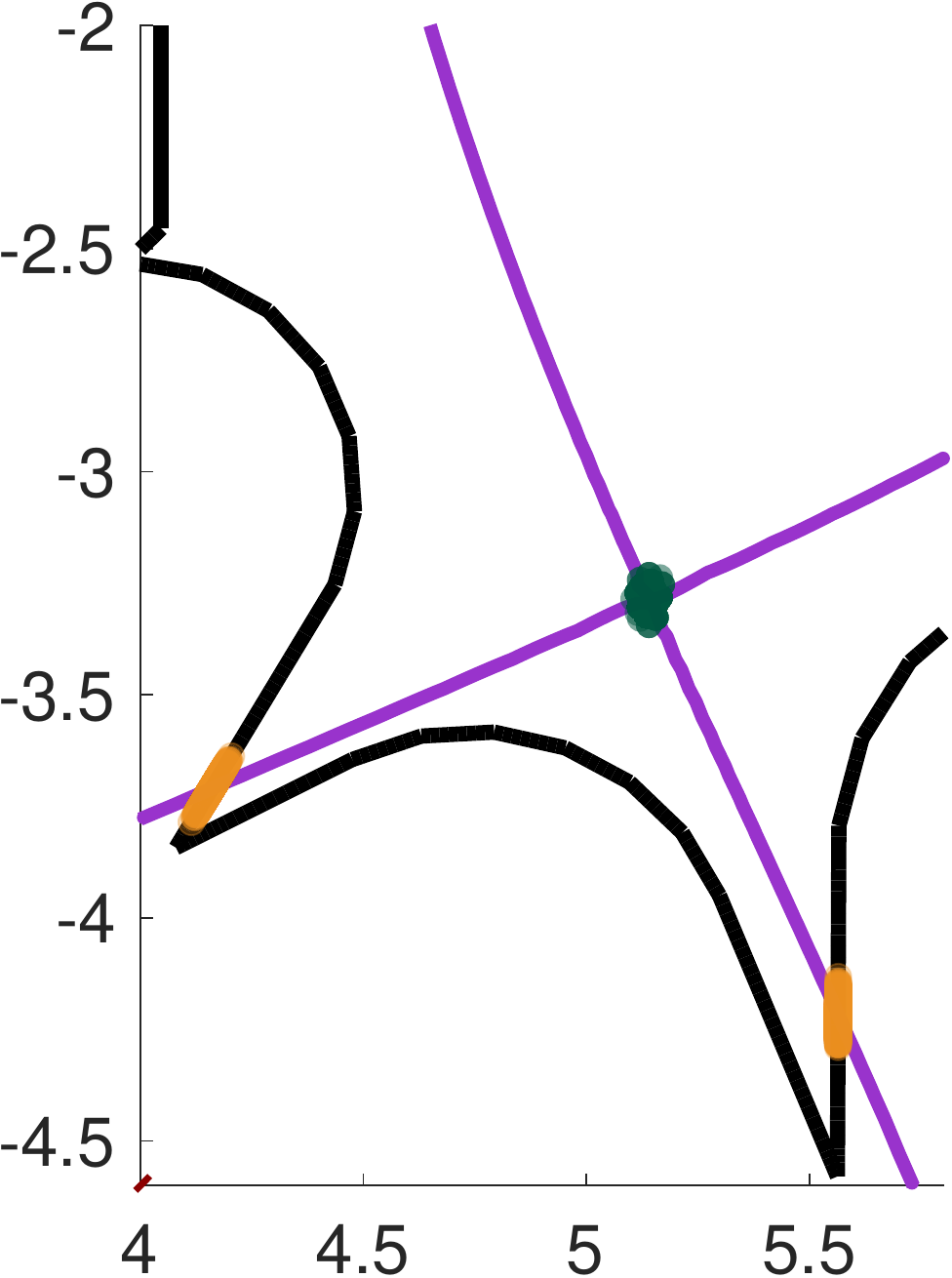} \quad&
\includegraphics[width=0.2\linewidth]{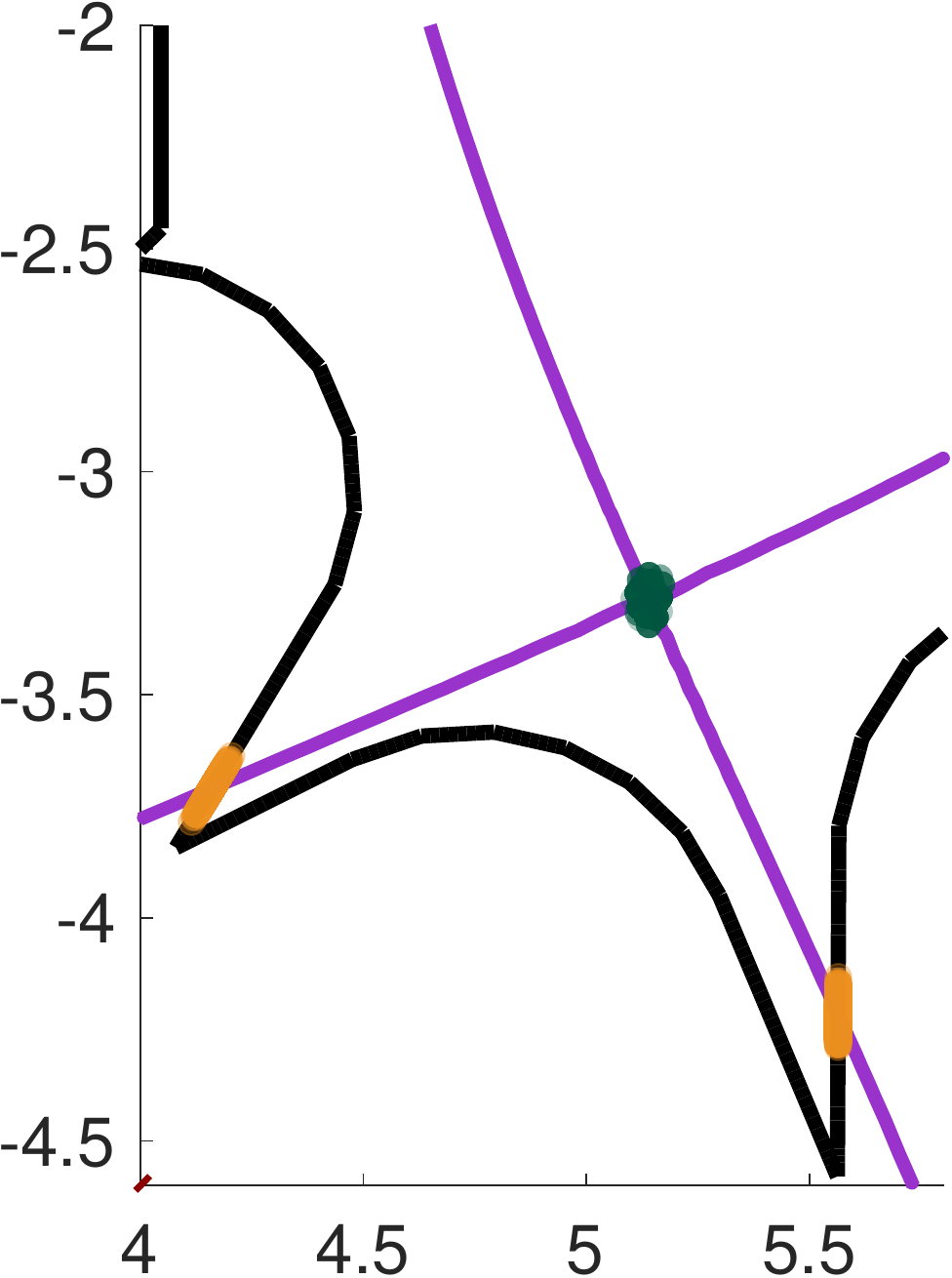}\\
\includegraphics[width=0.2\linewidth]{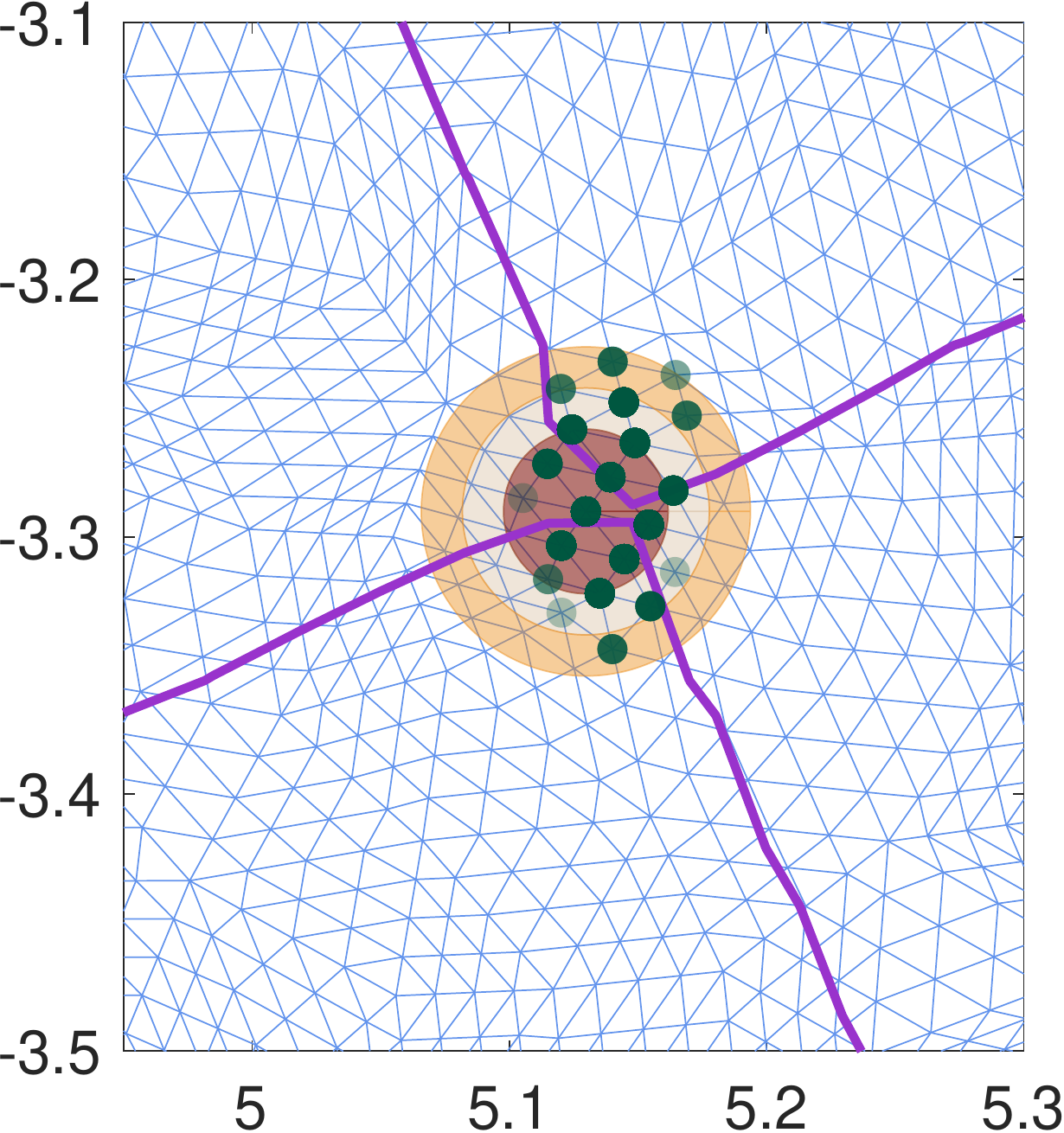} \quad&
\includegraphics[width=0.2\linewidth]{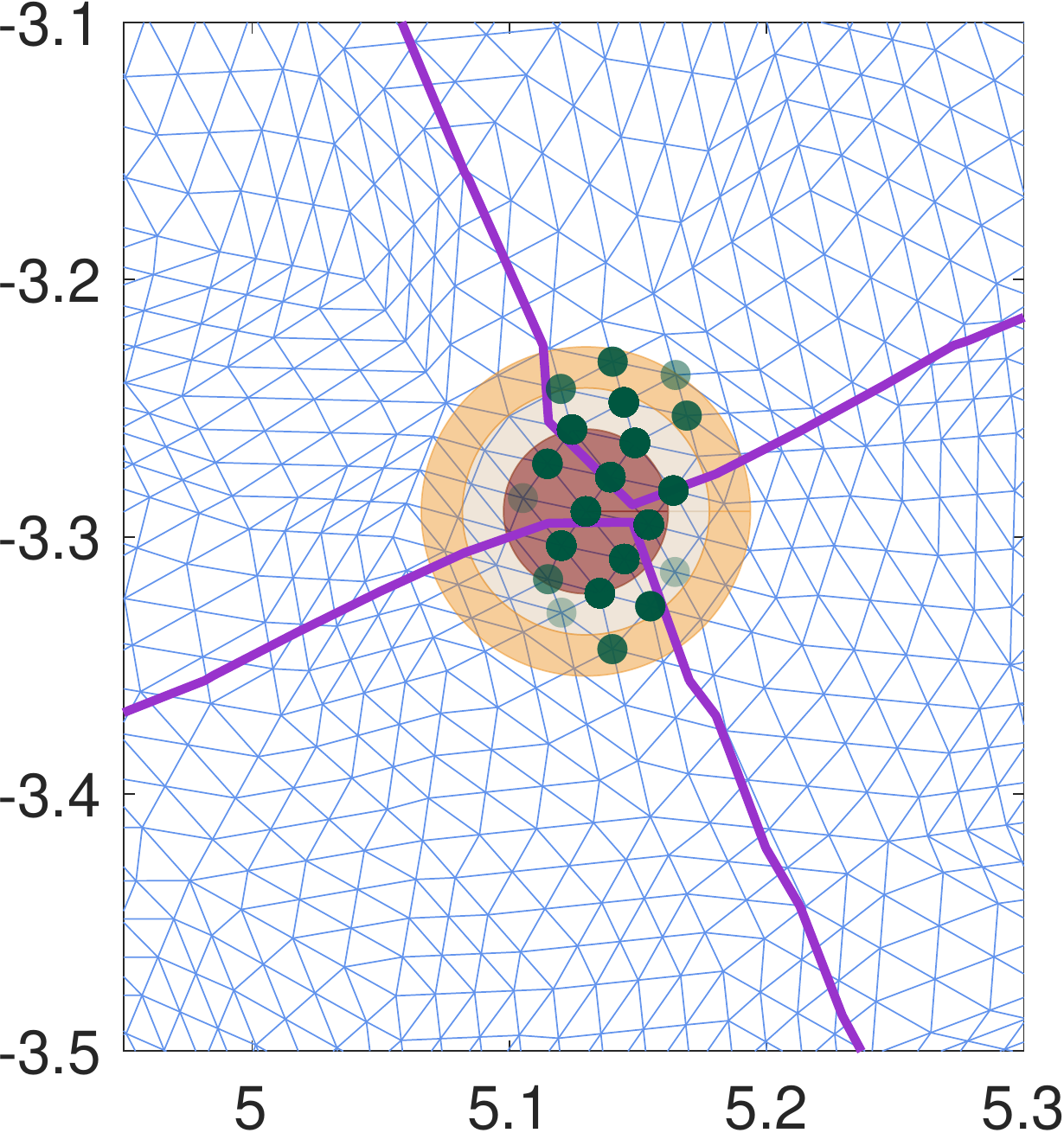} \quad&
\includegraphics[width=0.2\linewidth]{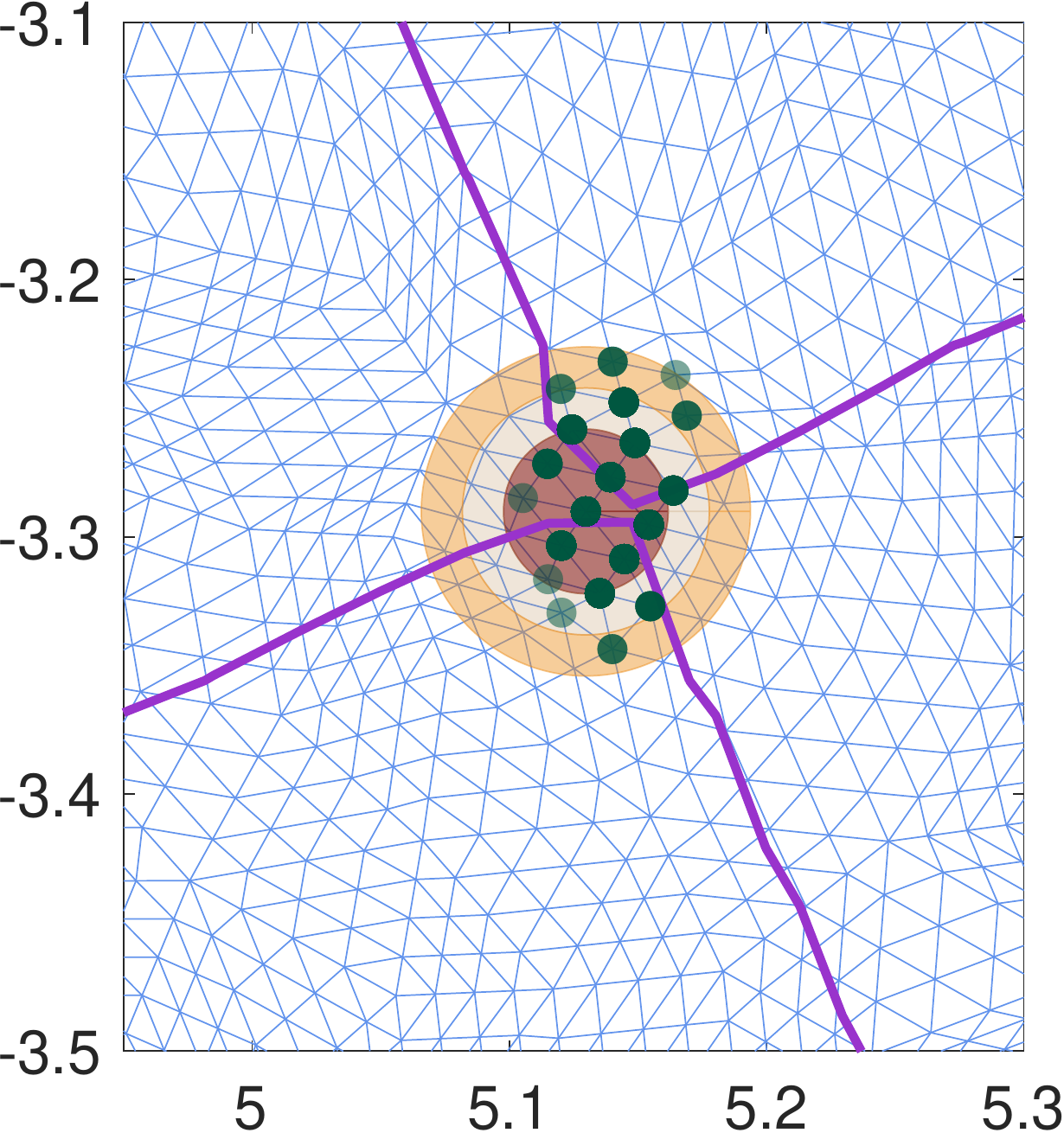} \quad&
\includegraphics[width=0.2\linewidth]{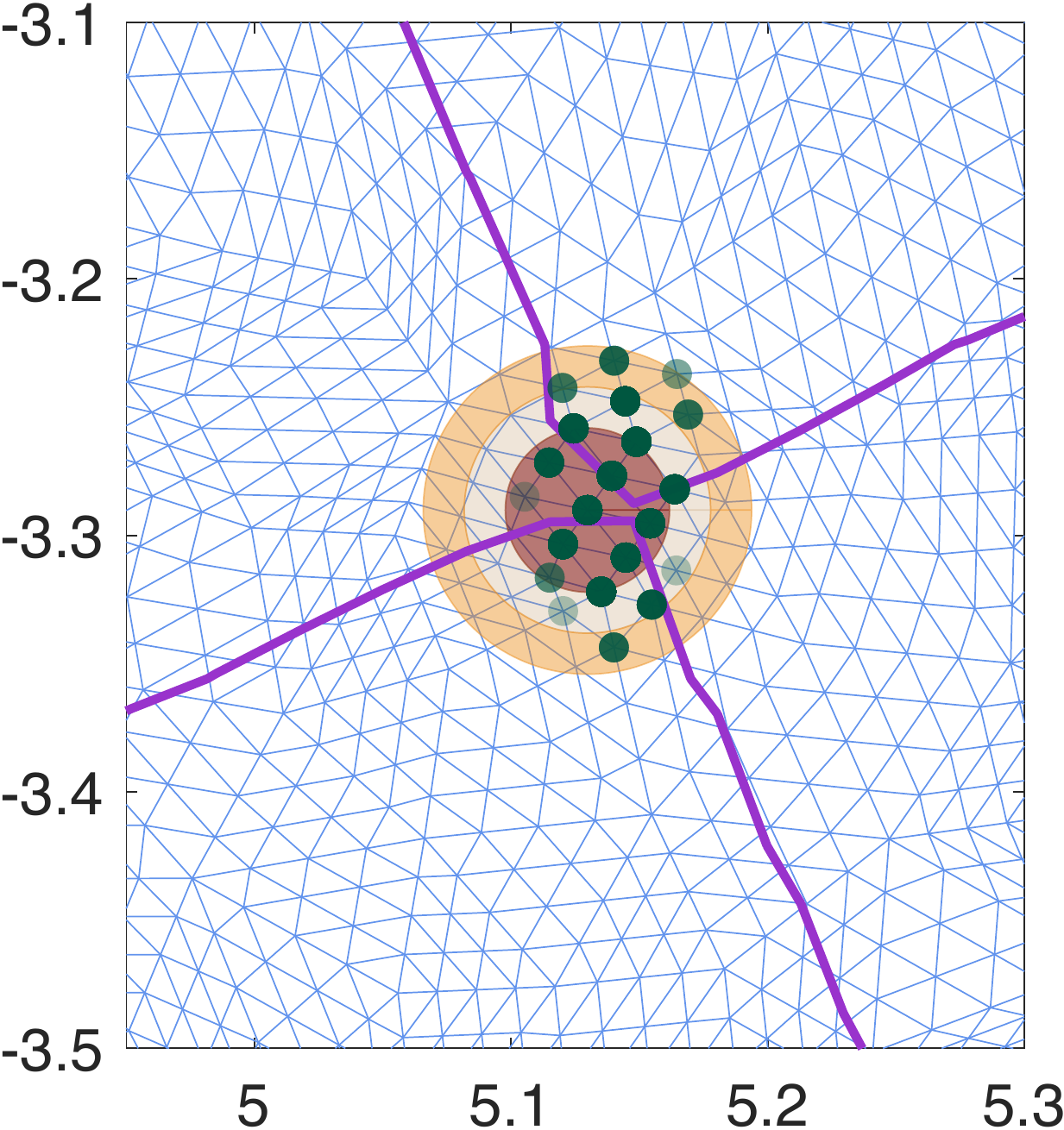}
\end{tabular}
\caption{Distributions of x-points (green) and strike points (orange) for 2000 samples of perturbed current values with 1\% noise. In all figures, scales are in meters, the reference plasma boundary is displayed in purple, the inner wall of the  reactor is displayed in black, and the outer walls of the reactor are displayed in dark red. From left to right, images correspond to surrogates obtained from sparse grid levels 2, 3 and 4, and direct solves. The top row shows (purple) the separatrix of the unperturbed problem. The second row shows the x-points and strike points in more detail. The third row shows further magnified neighborhoods of the x-points including the regions and perturbed x-points referenced in Table \ref{Tab:xpt-radii1}; frequency of x-points corresponds to intensity of their displays.}
\label{fig:pt-distribs1}
\end{figure}

\begin{figure}[h!]\centering
\begin{tabular}{cccc}
\includegraphics[width=0.2\linewidth]{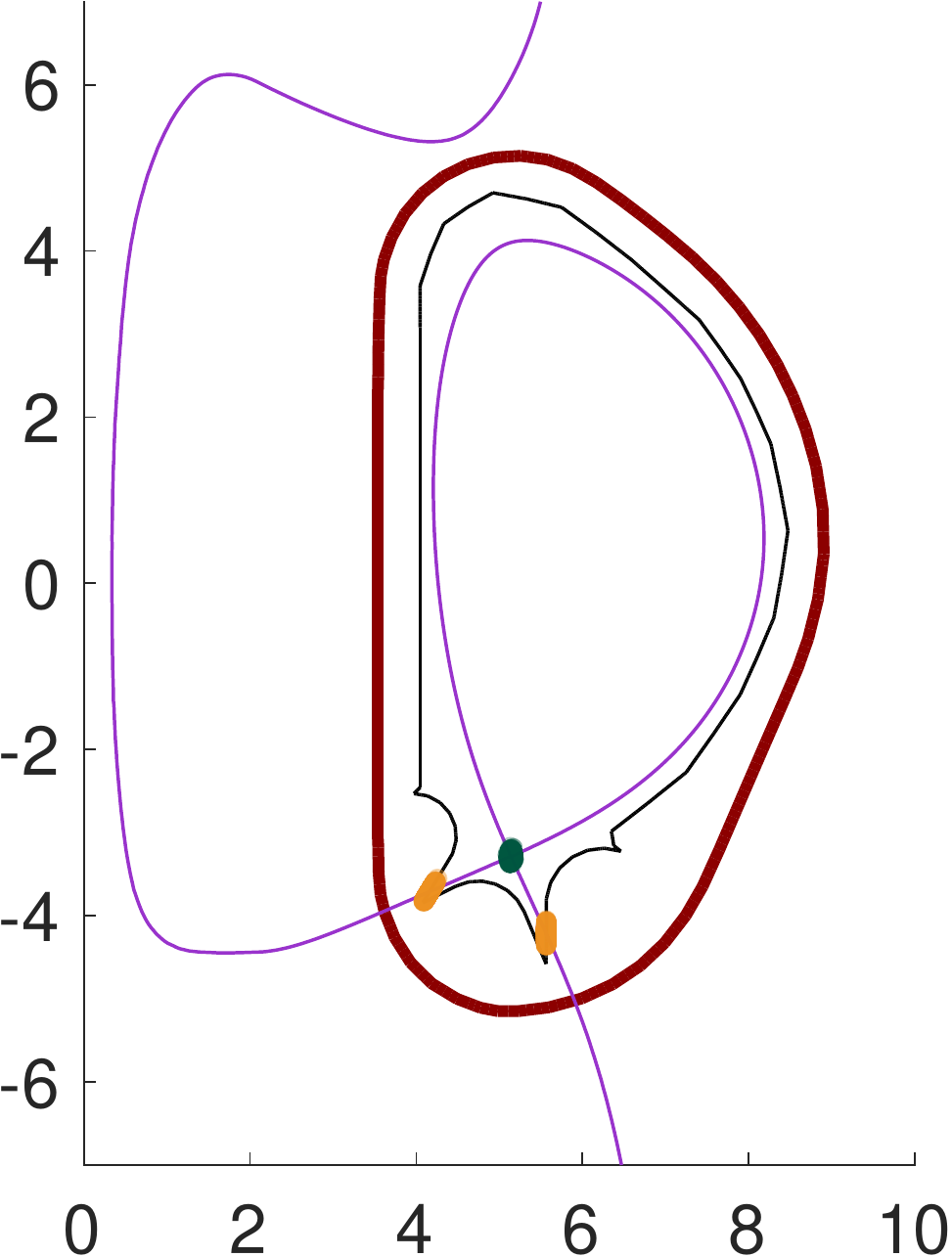} \quad&
\includegraphics[width=0.2\linewidth]{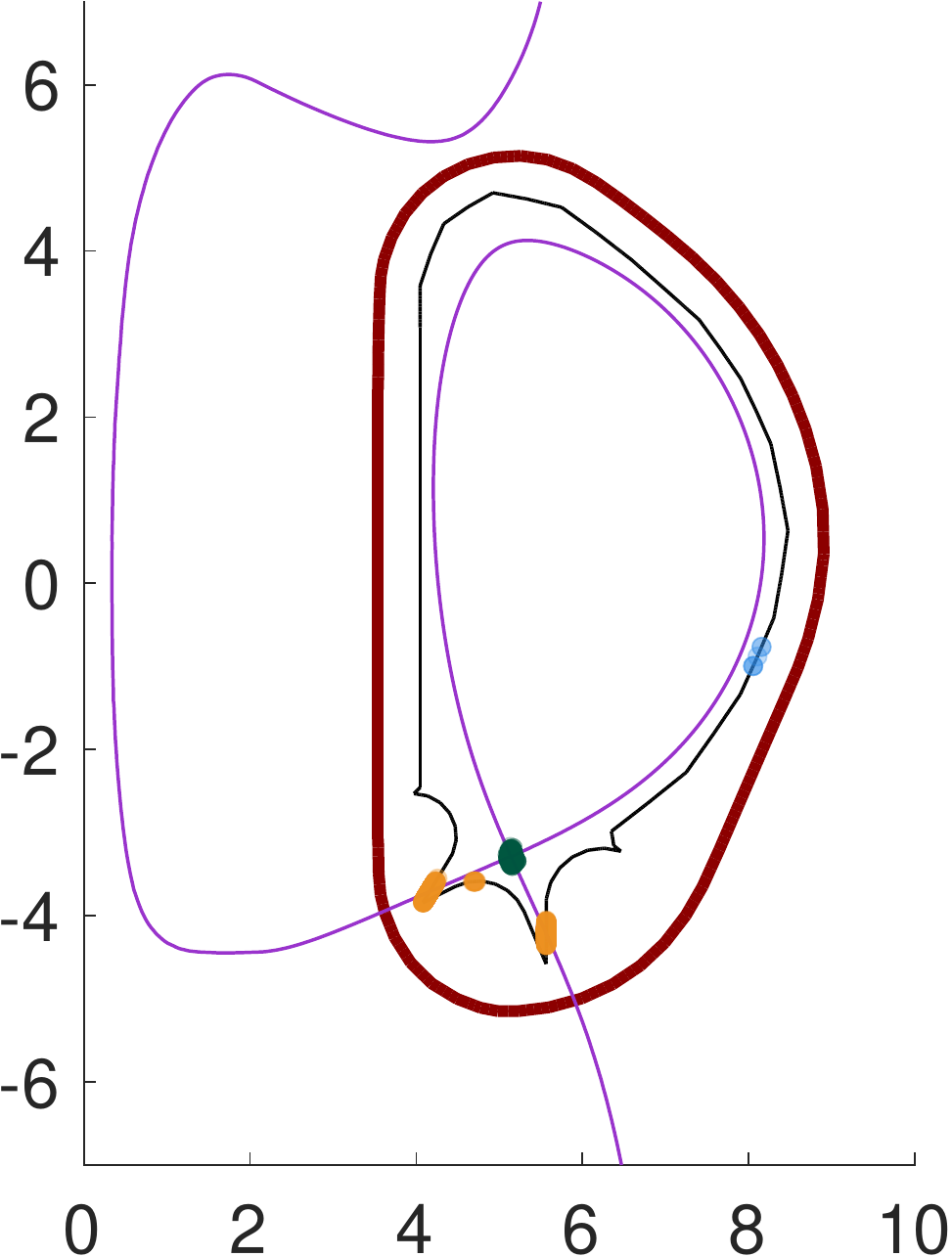} \quad&
\includegraphics[width=0.2\linewidth]{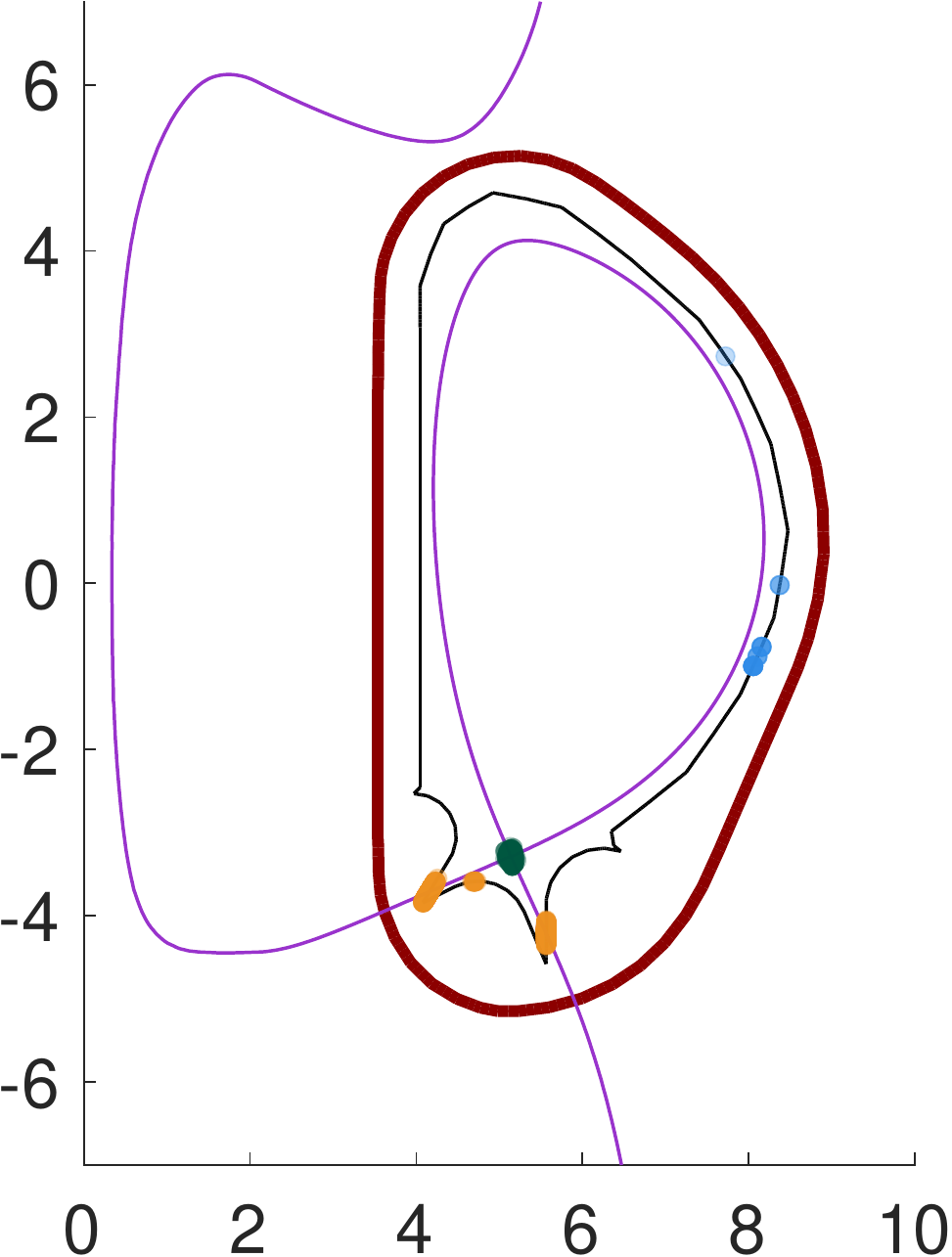} \quad&
\includegraphics[width=0.2\linewidth]{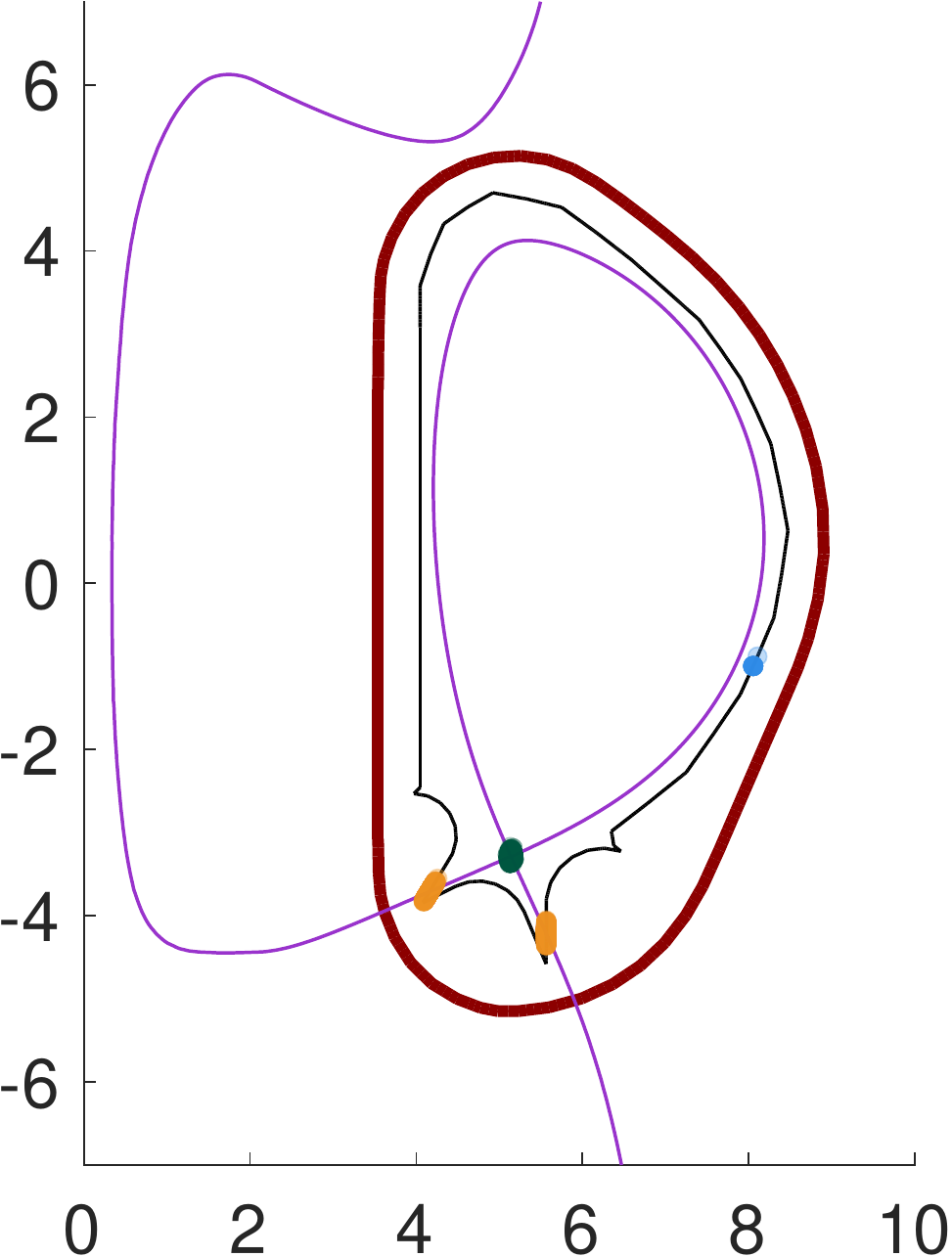}\\ 
\includegraphics[width=0.2\linewidth]{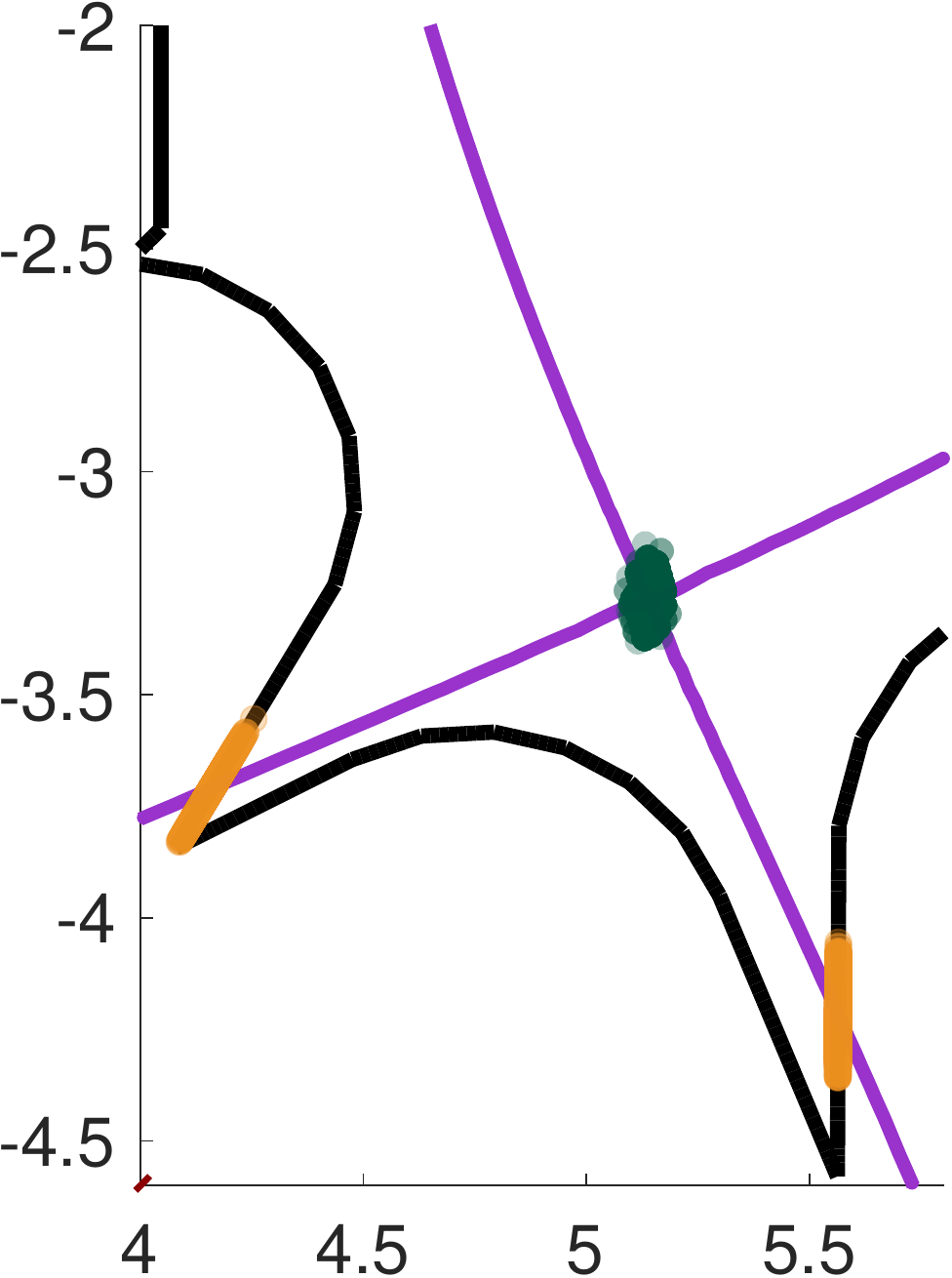} \quad&
\includegraphics[width=0.2\linewidth]{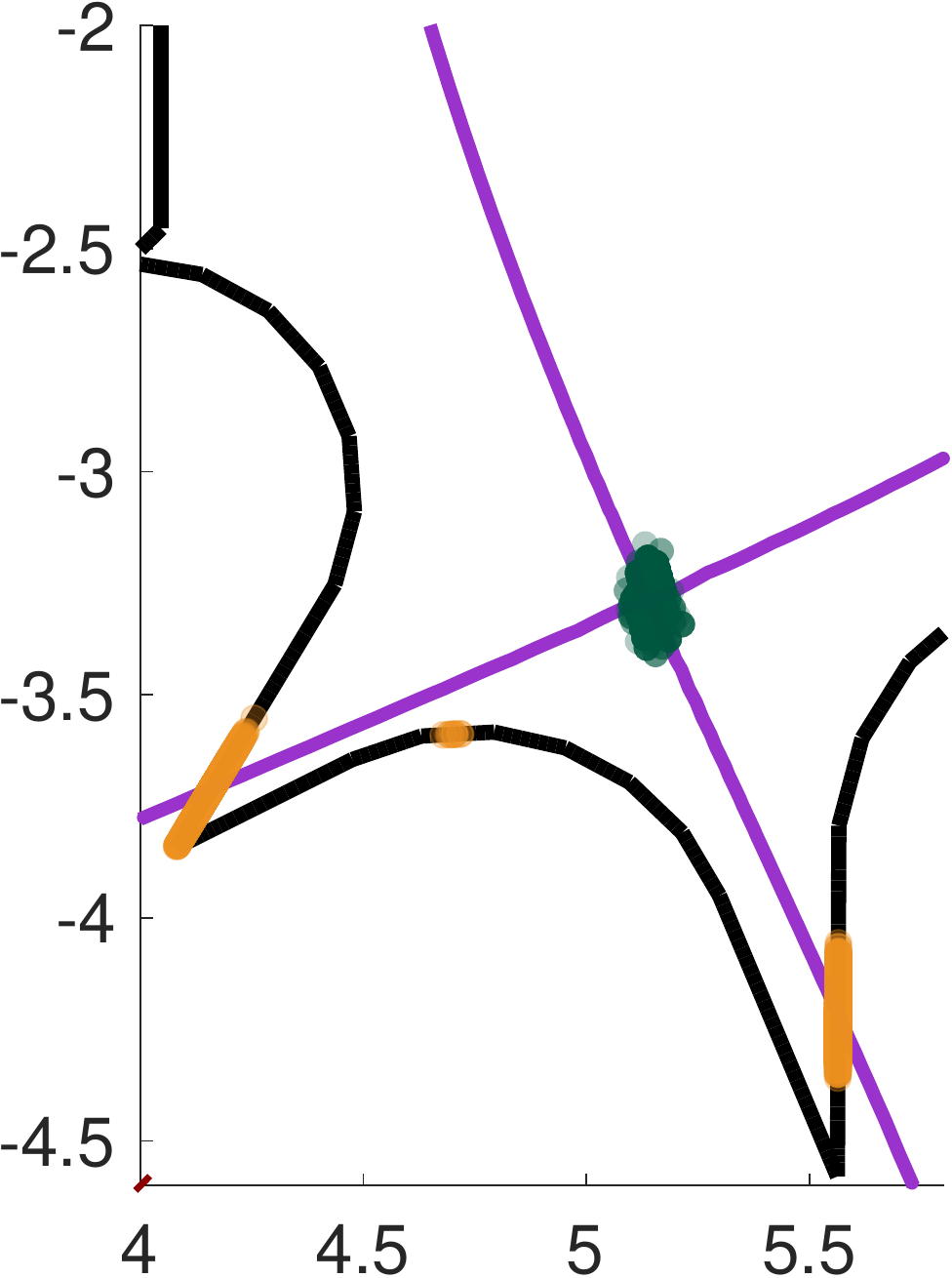} \quad&
\includegraphics[width=0.2\linewidth]{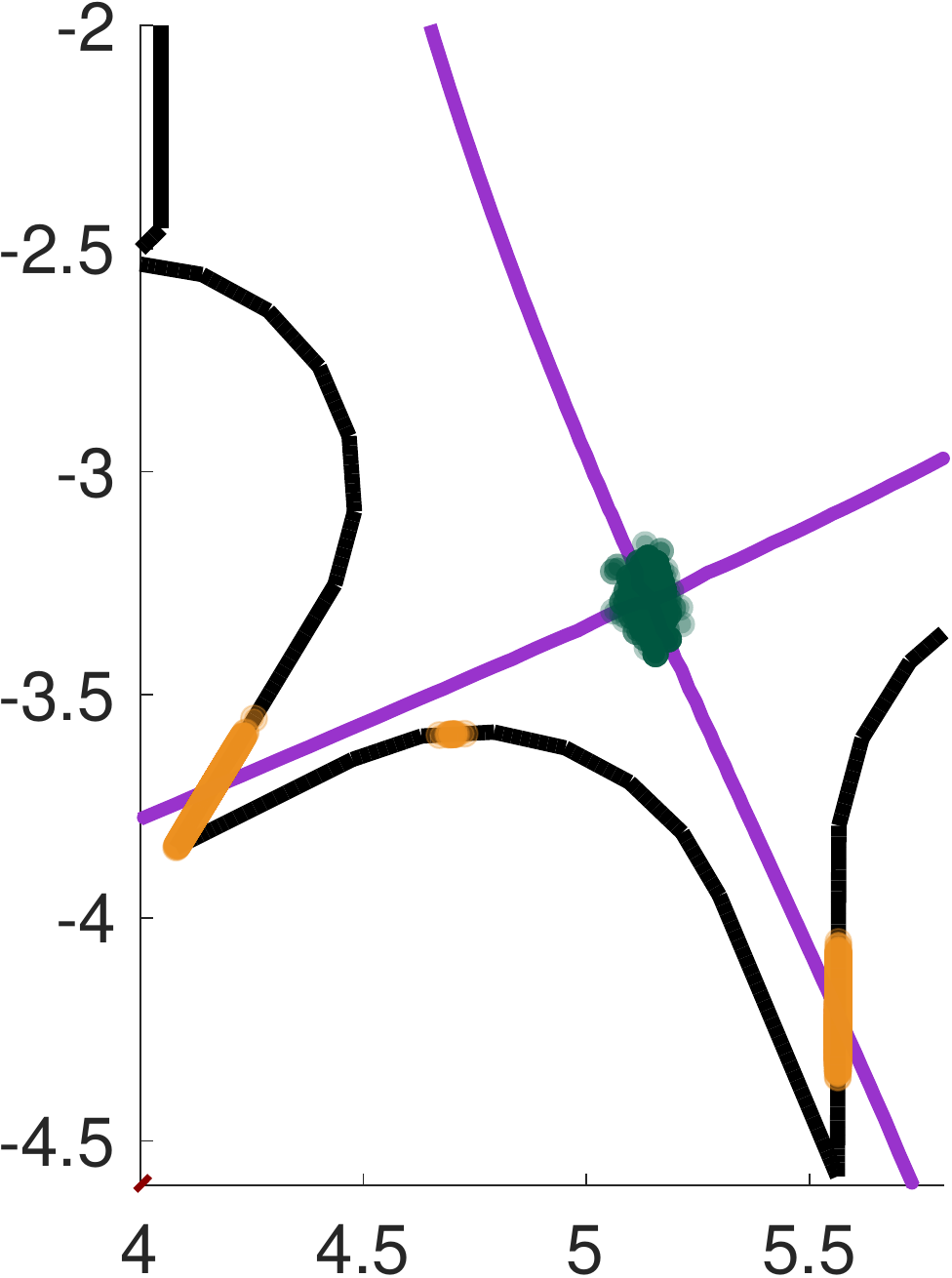} \quad&
\includegraphics[width=0.2\linewidth]{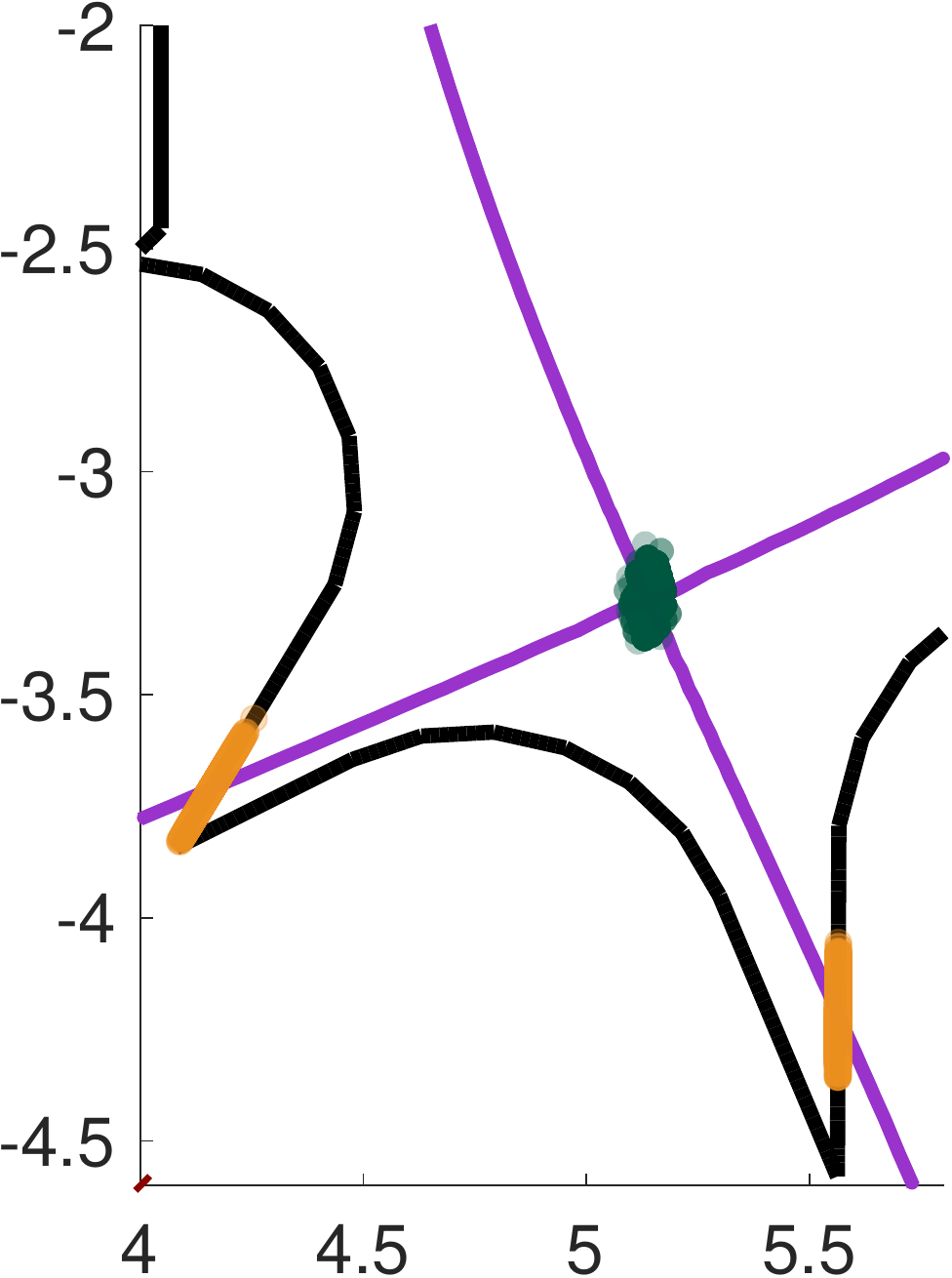}\\ 
\includegraphics[width=0.2\linewidth]{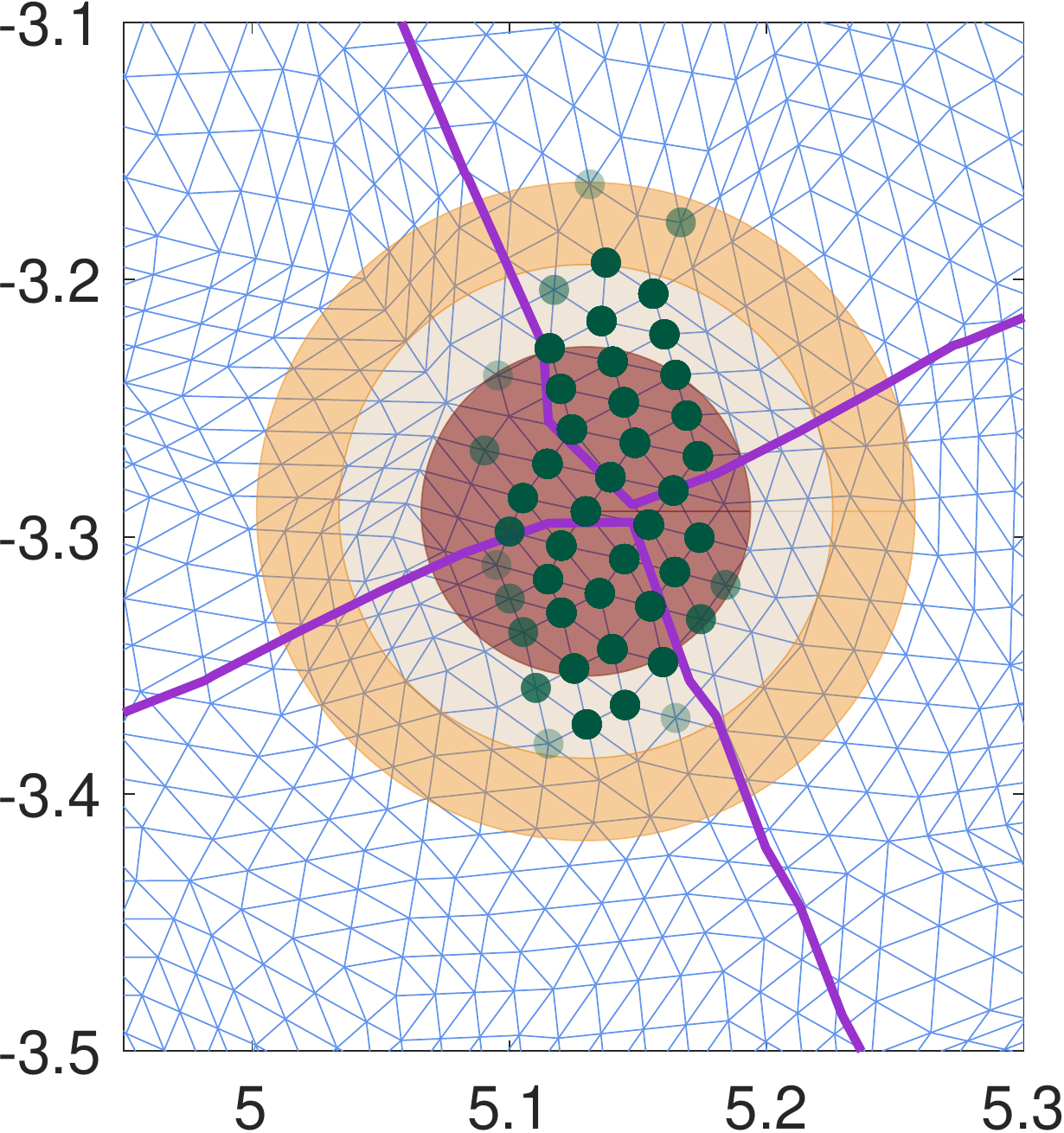} \quad&
\includegraphics[width=0.2\linewidth]{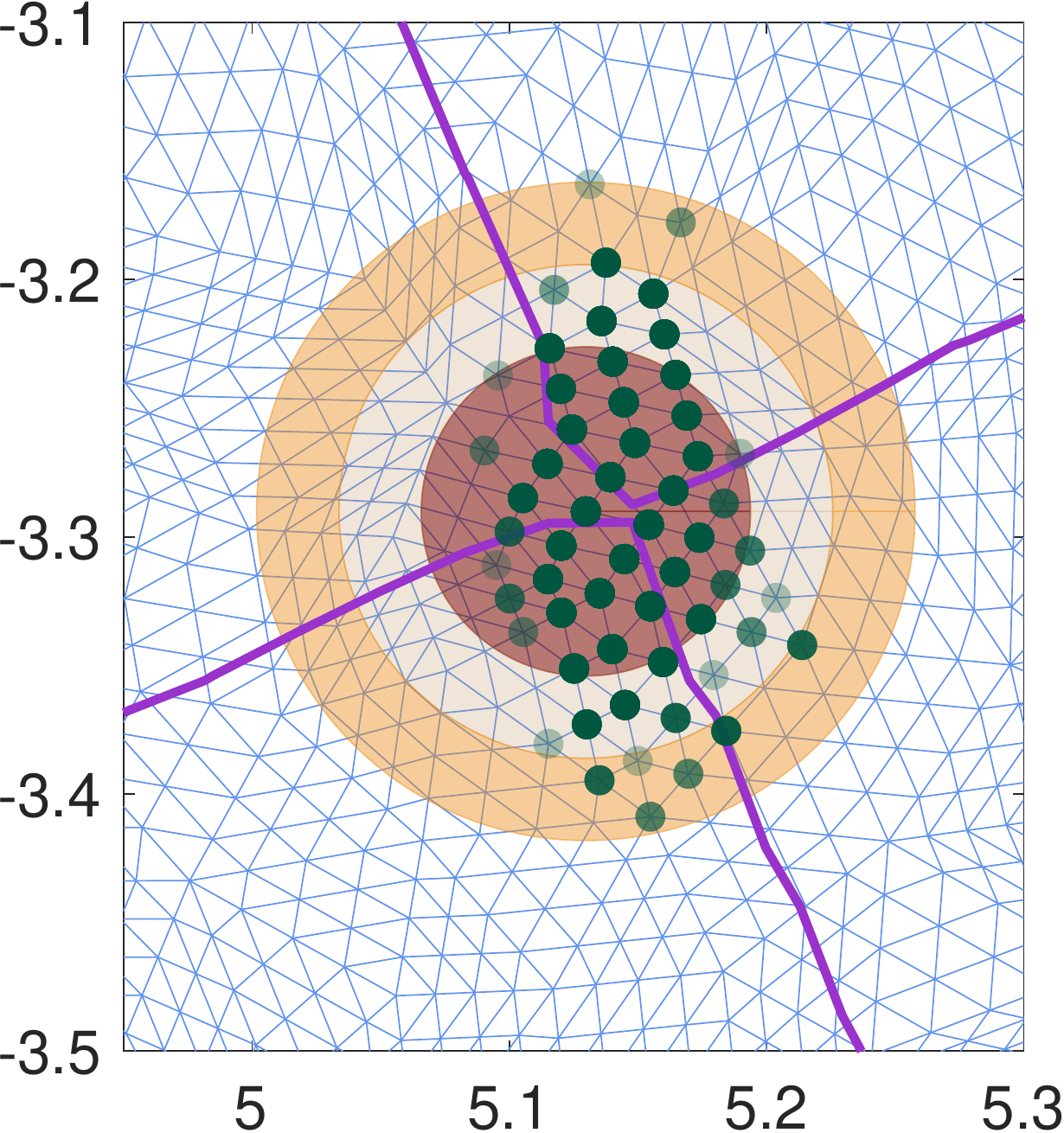} \quad&
\includegraphics[width=0.2\linewidth]{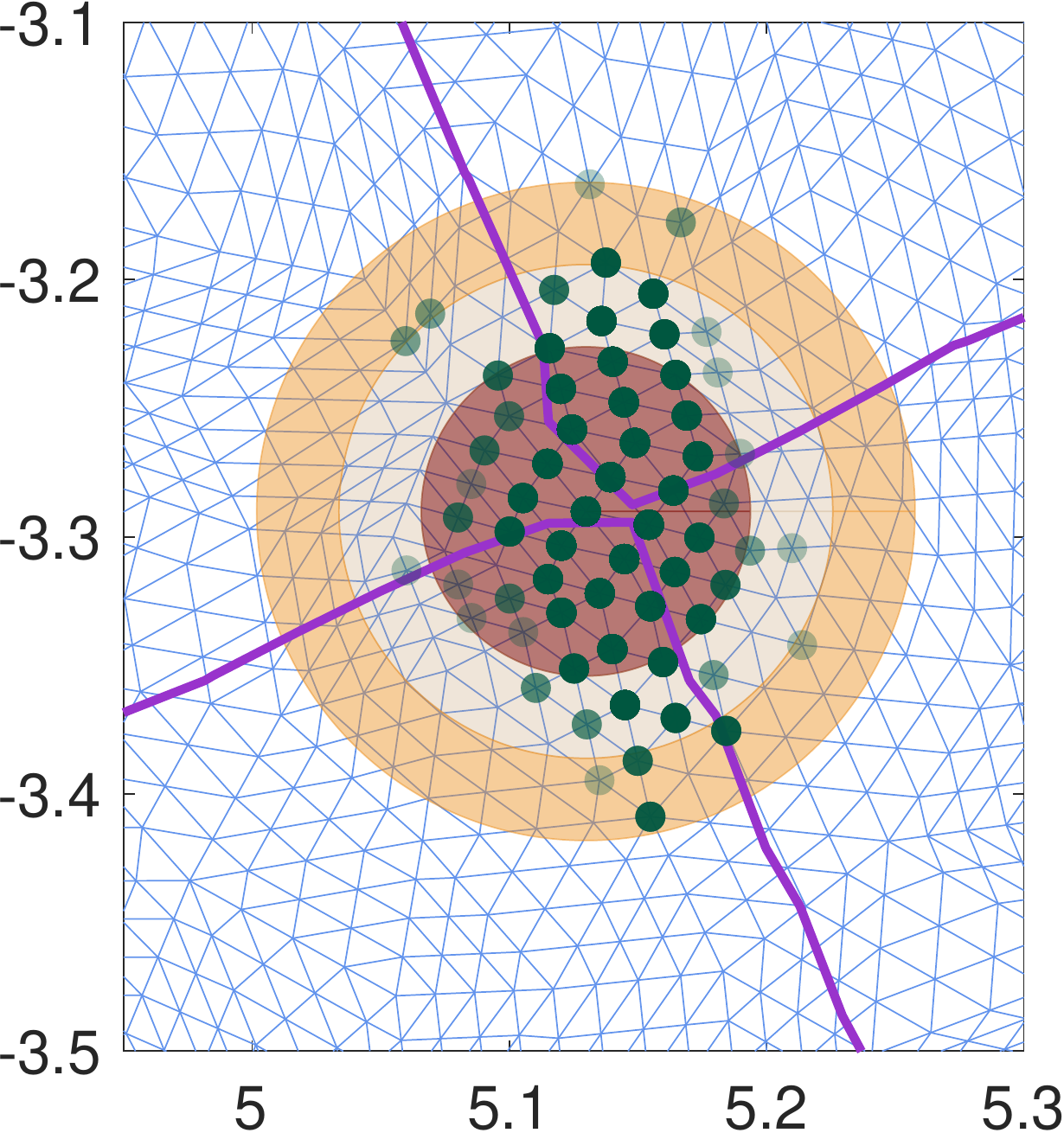} \quad&
\includegraphics[width=0.2\linewidth]{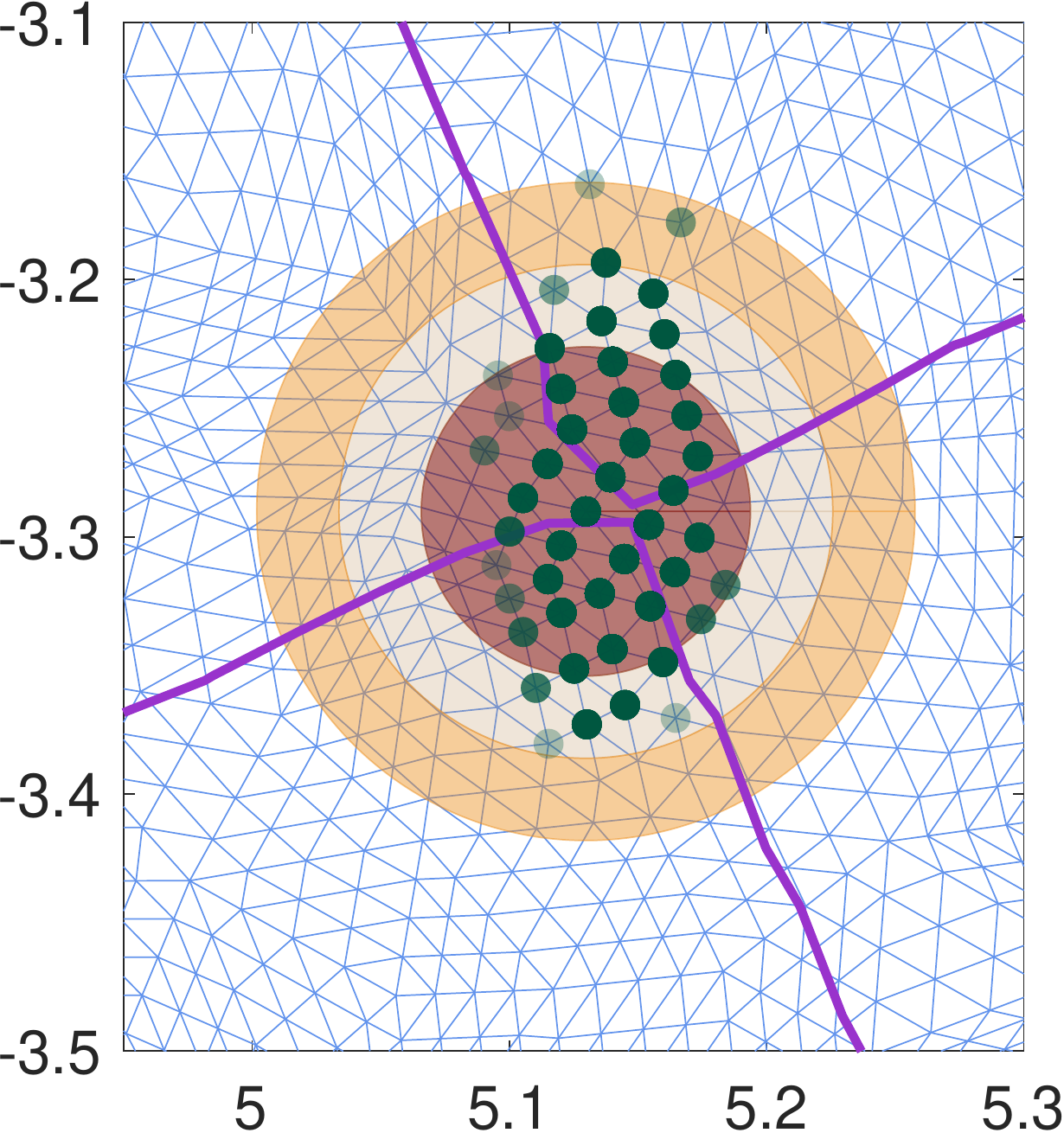}\\
\includegraphics[width=0.2\linewidth]{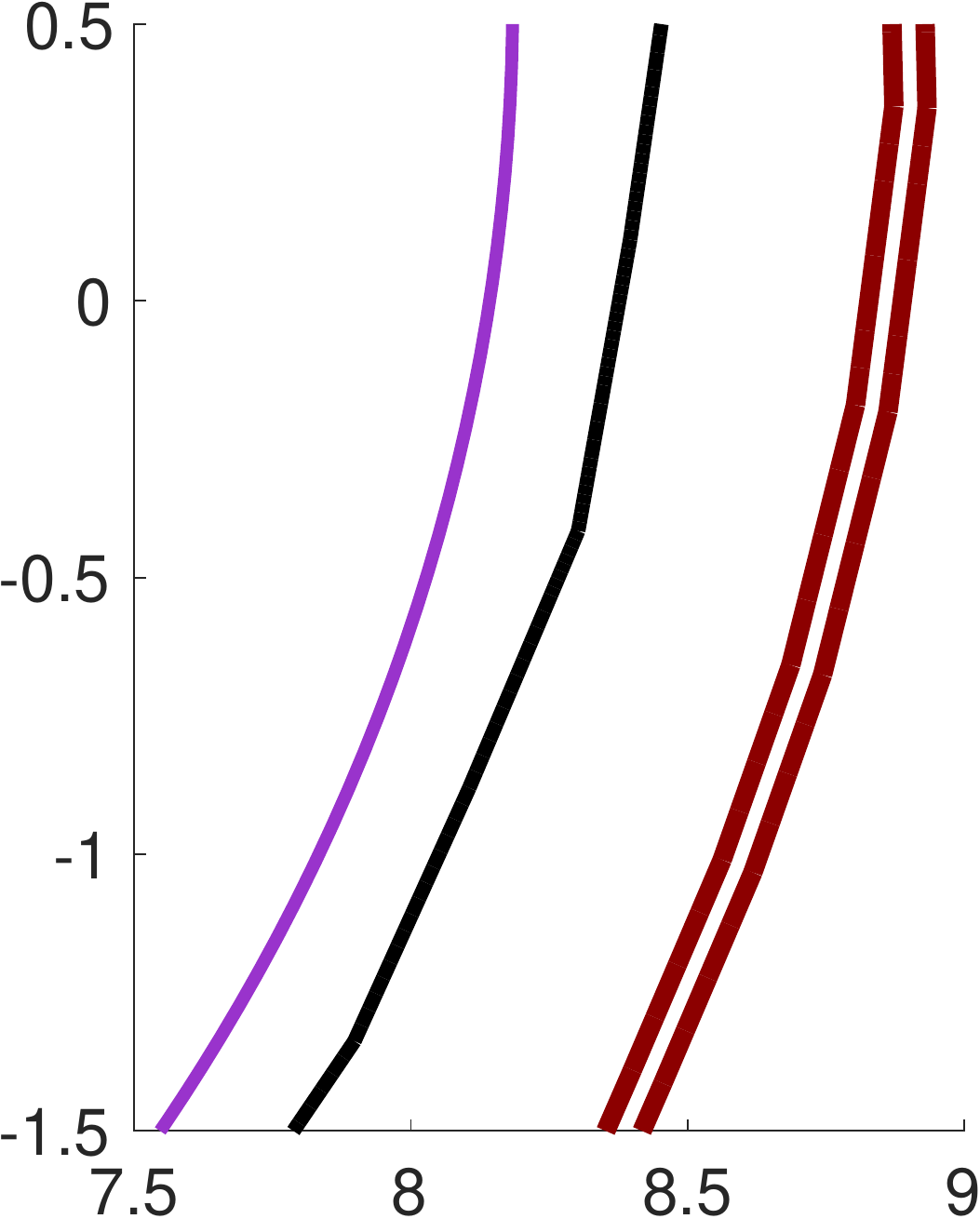} \quad&
\includegraphics[width=0.2\linewidth]{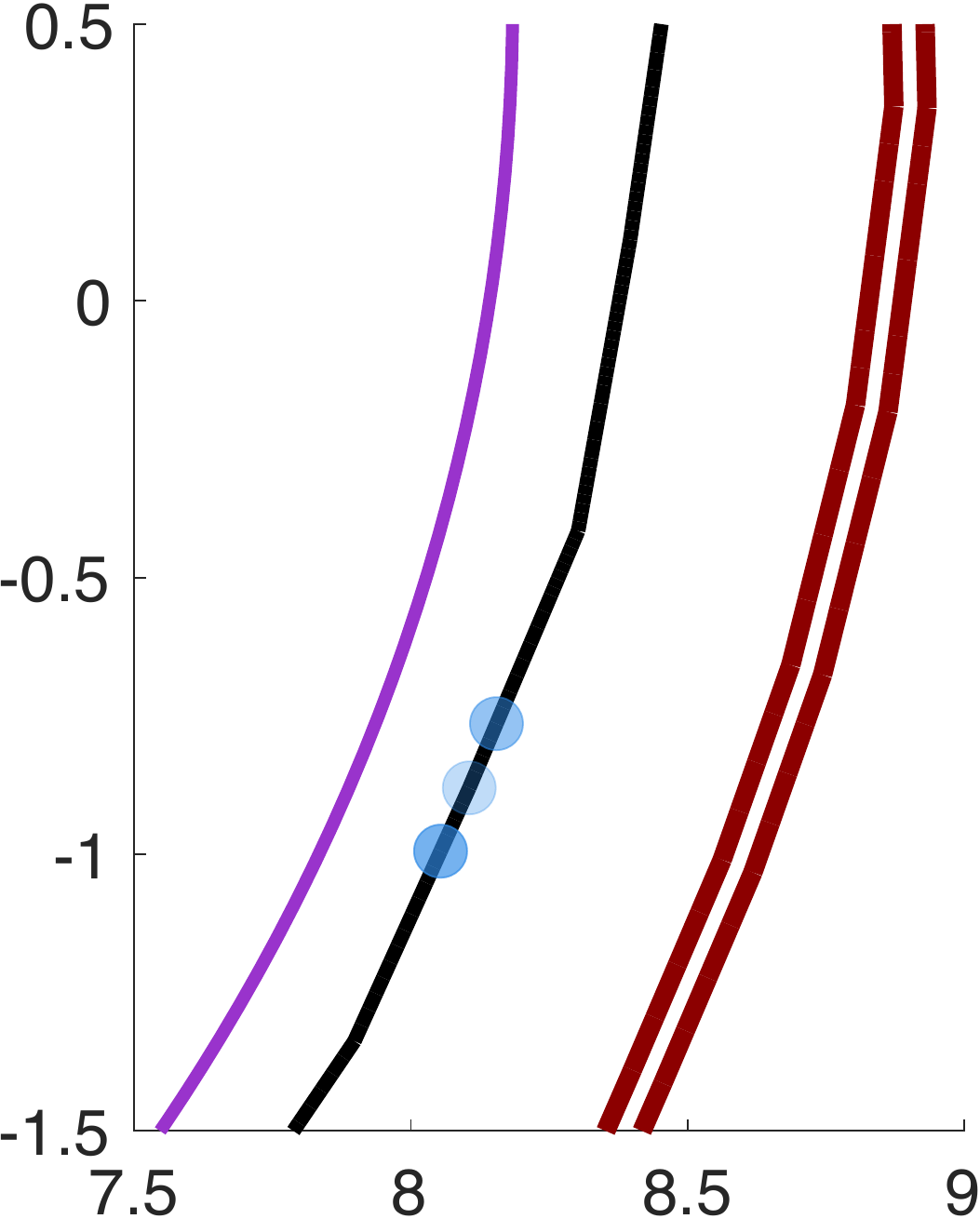} \quad&
\includegraphics[width=0.2\linewidth]{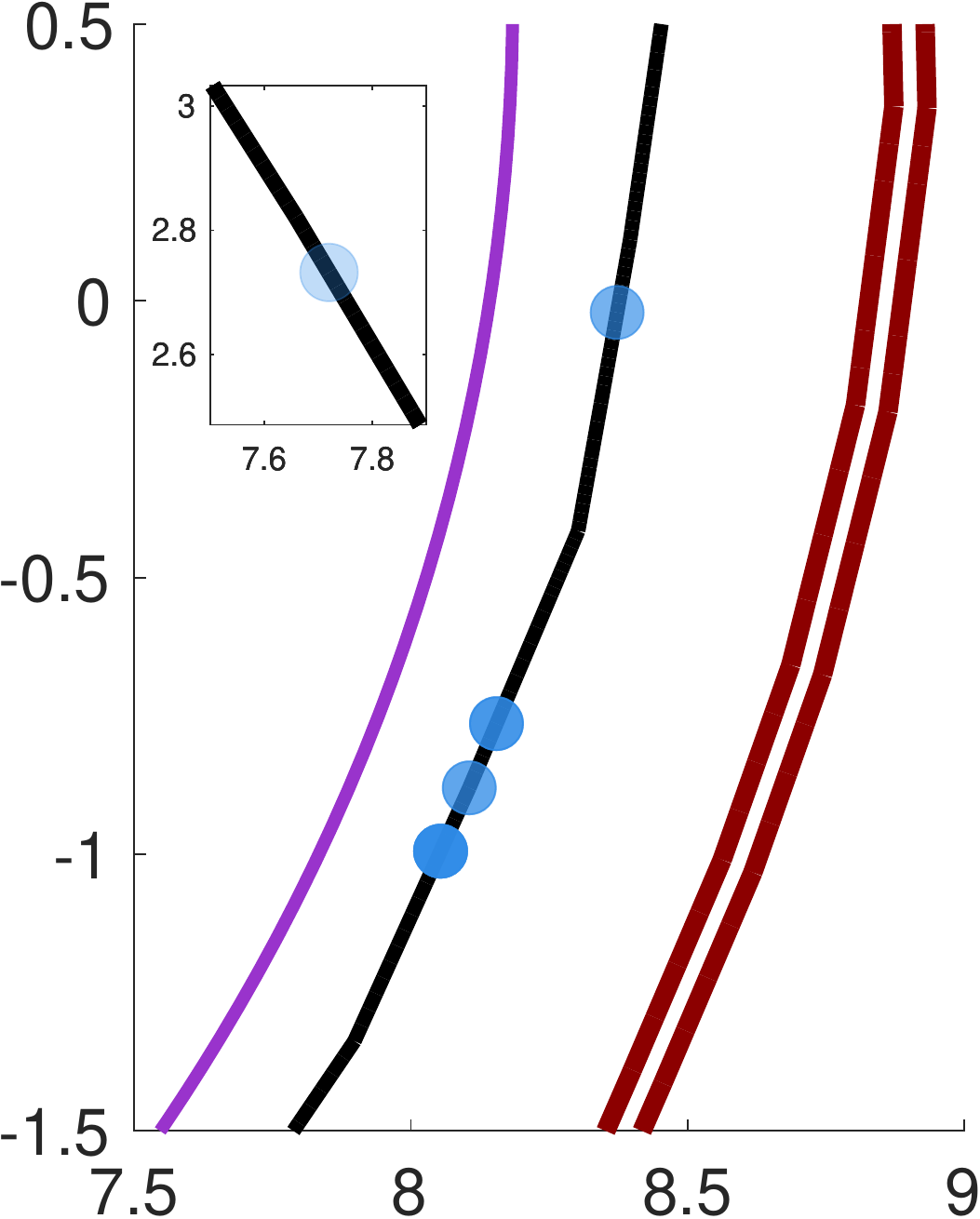} \quad&
\includegraphics[width=0.2\linewidth]{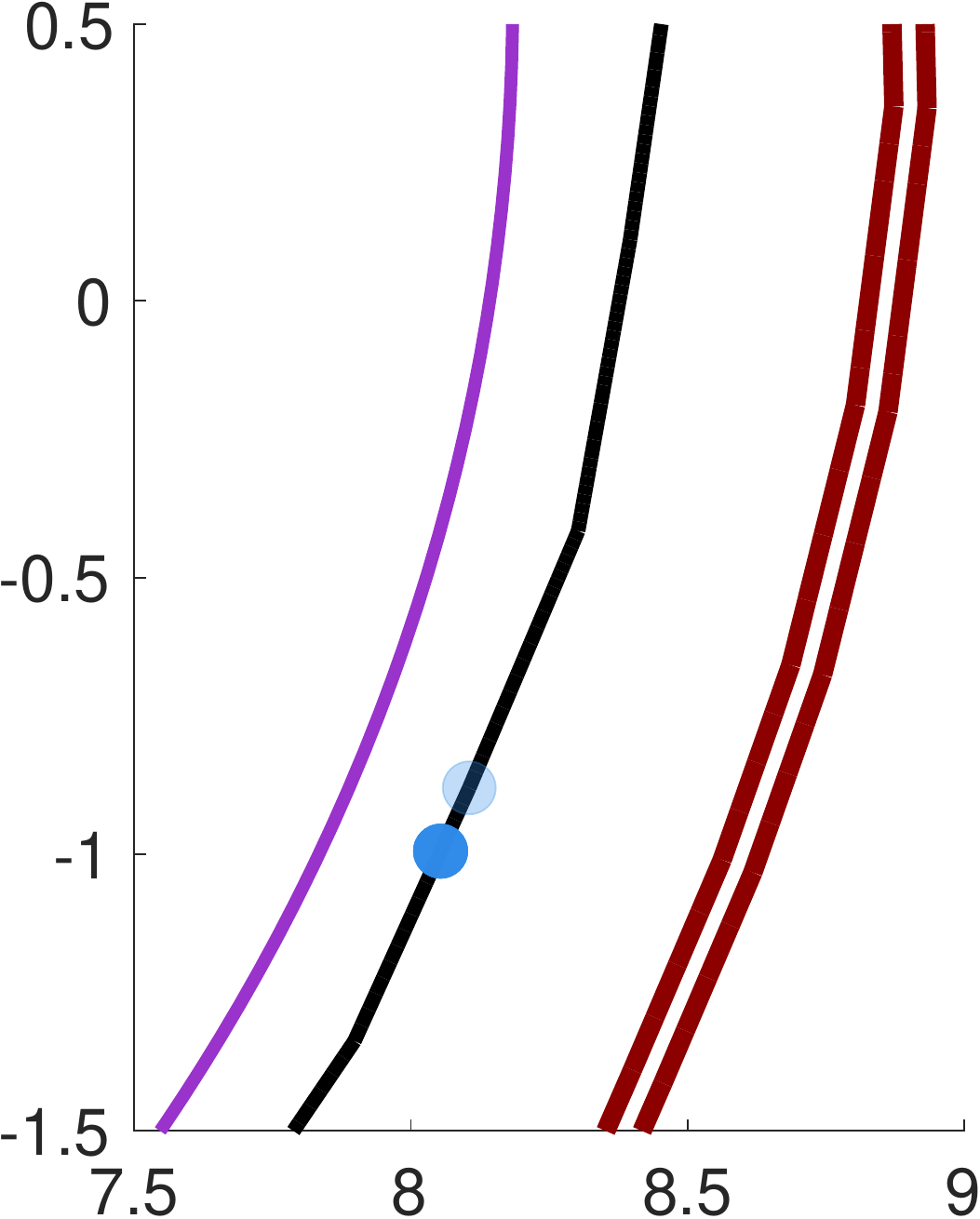}
\end{tabular}
\caption{Distributions of x-points (green), strike points (orange) and contact points (blue), for 2000 samples of perturbed current values with 2\% noise. In all figures, scales are in meters, the reference plasma boundary is displayed in purple, the inner wall of the  reactor is displayed in black, and the outer walls of the reactor are displayed in dark red. From left to right, images correspond to surrogates obtained from sparse grid levels 2, 3 and 4 and direct solves. The top row shows (in purple) the separatrix obtained from the unperturbed problem. The second row shows the x-points and strike points in more detail. The third row shows further magnified neighborhoods of the x-points including the regions and perturbed x-points referenced in Table \ref{Tab:xpt-radii2}, and the fourth row shows contact points. In the two bottom rows, frequency of x-points and contact points corresponds to intensity of their displays.} 
\label{fig:pt-distribs2}
\end{figure}

\begin{table}[ht]
	\centering
	\scalebox{0.8}{
		\begin{tabular}{c|c|c|c|c|c|c|c|c|c|c|c|c|c|c|}
			\cline{2-5}
			& \multicolumn{4}{|c|}{Number of x-points between $r_{i-1}$ and $r_i$ \;($r_0 = 0$)}\\
			\cline{2-5}
			&\multicolumn{1}{|c|}{Surrogate level 2} &\multicolumn{1}{c|}{Surrogate level 3} 
			      &\multicolumn{1}{c|}{Surrogate level 4} &  \multicolumn{1}{c|}{Direct solver} \\
			\hline
			\multicolumn{1}{|l|}{$r_1=0.032=.5r_3$} &1455 (72.75\%)&1457 (72.85\%)&1459( 72.95\%)&1457 (72.85\%)\\
			\multicolumn{1}{|l|}{$r_2=0.048=.75r_3$} &522 (26.10\%)&521 (26.05\%)&518 (25.90\%)&520 (26.00\%)\\
			\multicolumn{1}{|l|}{$r_3=0.064$}             & 23 (1.15\%)&22 (1.10\%)&23 (1.15\%)&23 (1.15\%)\\
			\hline
	\end{tabular}}
				\caption{Number (frequencies) of x-points found within various distances from the
				reference x-point, for 2000 perturbations with 1\% noise in 12 coils. 
				Counts correspond to number of points of distance at most $r_1$ or within annuli of width $r_2-r_1$ or $r_3-r_2$. Distances are measured in meters.}
		\label{Tab:xpt-radii1}
\end{table}

\begin{table}[ht]
	\centering
	\scalebox{0.8}{
		\begin{tabular}{|c|c|c|c|c|c|c|c|c|c|c|c|}
			\hline	
			Level &\# Nodes & $E_{2000}$ & Inverse Aspect & Elongation & Upper            & Lower & Strike & Strike & x-point\\
			           & & & Ratio & & Triangularity & Triangularity & Point 1 & Point 2 &\\
			\hline
			 0&1&6.39e-03&2.76e-03&2.61e-03&1.65e-02&1.19e-02&3.96e-03&3.36e-03&3.34e-03\\
			 1&25&1.58e-04&2.81e-05&9.70e-05&8.92e-04&4.24e-04&4.77e-05&8.94e-06&6.32e-05\\
			 2&313&1.94e-04&2.59e-05&4.86e-05&4.52e-04&3.37e-04&3.19e-05&4.48e-06&5.01e-05\\
			 3&2649&2.54e-04&2.77e-05&5.11e-05&6.75e-04&2.97e-04&3.72e-05&8.67e-06&4.22e-05\\
			 4&17265&3.96e-04&3.49e-05&6.15e-05&1.31e-03&3.76e-04&4.69e-05&9.97e-06&5.47e-05\\
			\hline	
	\end{tabular}}
		\caption{Mean relative errors of various geometric features. The relative approximation error $E_{n_s}$ is defined in \eqref{eq:MeanError}. The plasma shaping parameters are given in (\ref{eq:ShapeParameters}). The experiment used 2000 random samples with 1\% perturbation of all 12 currents, for increasing collocation levels. 
		} \label{Tab:Errors-1}
\end{table}

Next, we will explore how effective the surrogate model is for extracting information about the confinement properties of  the magnetic field. To assess the effect that the stochasticity may have on physically relevant quantities, such as the location of x-points and the incidence of interactions between the plasma and the walls of the reactor, we will make use of a ``reference equilibrium configuration" corresponding to the ideal case where no variability is present in the coil currents. In the experiments used for this discussion, the surrogate was built using two levels of local spatial adaptivity in the free boundary solver, and surrogate solutions were obtained on a common grid constructed from two levels of uniform refinement. This ensures that the size of the smallest elements in the mesh where the surrogate is built is comparable to that of the adaptive grids where the samples at collocation nodes are taken.

Figure \ref{fig:pt-distribs1} depicts the distributions of x-points and strike points found by surrogate approximations determined from 2000 perturbed current values, with 1\%, for three levels of sparse grids, together with the analogous results obtained from direct solution. Visually, the results obtained from all instances of the surrogate computations are indistinguishable from those obtained from direct simulation. These results are quantified in Table \ref{Tab:xpt-radii1}, which counts the numbers of x-points that are found at various distances from a reference x-point, defined to be the point obtained from unperturbed input. These data are also displayed in the bottom row of Figure \ref{fig:pt-distribs1}.  They indicate that the trends for the surrogates and direct solutions are in fact quantitatively similar and give additional credence for the accuracy of the surrogates. Table \ref{Tab:Errors-1} continues this exercise by showing the errors for a collection of geometric features of the plasma boundary 
(\ref{eq:ShapeParameters}). The results indicate that reasonable accuracy is obtained from surrogates constructed using sparse grid levels 2 or 3.

\begin{table}[ht]
	\centering
	\scalebox{0.8}{
		\begin{tabular}{c|c|c|c|c|c|c|c|c|c|c|c|c|c|c|}
			\cline{2-5}
			&\multicolumn{4}{|c|}{Number of x-points between $r_{i-1}$ and $r_{i} $\; $(r_0=0)$}\\
			\cline{2-5}	
			&\multicolumn{1}{|c|}{Surrogate level 2} &\multicolumn{1}{c|}{Surrogate level 3} 
			      &\multicolumn{1}{c|}{Surrogate level 4} &  \multicolumn{1}{c|}{Direct solver} \\			
			\hline
			\multicolumn{1}{|l|}{$r_1=0.064=.5r_3$}&1819 (90.95\%)&1773 (88.92\%)&1723 (87.24\%)&1814 (91.39\%)\\
		         \multicolumn{1}{|l|}{$r_2=0.096=.75r_3$}&168 (8.40\%)&179 (8.98\%)&209 (10.58\%)&157 (7.91\%)\\
			\multicolumn{1}{|l|}{$r_3=0.128$}&13 (0.65\%)&42 (2.11\%)&43 (2.18\%)&14 (0.71\%)\\
			\hline
			&\multicolumn{4}{|c|}{Number of plasma / reactor contacts}\\
			\cline{2-5}	
			 & 0 & 6 & 25 & 15 \\
			\cline{2-5}
	\end{tabular}}	
		\caption{ Incidence of contacts between the plasma and the reactor,
		 along with number (frequencies) of x-points found within various distances from the 
				reference x-point, for 2000 perturbations with 2\% noise in 12 coils. 
				Counts of x-points correspond to number of points of distance at most $r_1$ or 
				 within annuli of width $r_2-r_1$ or $r_3-r_2$.}  		
		\label{Tab:xpt-radii2}
\end{table}

\begin{table}[ht]
	\centering
	\scalebox{0.8}{
		\begin{tabular}{|c|c|c|c|c|c|c|}
			\hline	
			Level &\# Nodes & $E_{2000}$ & Inverse Aspect & Elongation & Upper & Lower \\
			& & & Ratio & & Triangularity & Triangularity \\
			\hline
			0&1&1.26e-02&5.46e-03&4.95e-03&2.93e-02&1.65e-02\\
			1&25&2.69e-04&8.49e-05&3.02e-04&3.05e-03&1.11e-03\\
			2&313&2.37e-04&4.19e-05&6.91e-05&3.56e-03&4.49e-04\\
		         3&2649&1.69e-03&7.98e-04&9.30e-04&6.44e-03&3.21e-03\\
			 4&17265&4.80e-03&2.41e-03&2.66e-03&9.10e-03&6.10e-03\\
			\hline	
	\end{tabular}}
		\caption{Mean relative errors of various geometric features. The relative approximation error $E_{n_s}$ is defined in \eqref{eq:MeanError}. The plasma shaping parameters are given in (\ref{eq:ShapeParameters}). The experiment used 2000 random samples with 2\% perturbation of all 12 currents, for increasing collocation levels.
		}		\label{Tab:Errors-2}
\end{table}

Figure \ref{fig:pt-distribs2} and Tables \ref{Tab:xpt-radii2} and \ref{Tab:Errors-2} show the results of similar experiments in which the size of the perturbations is increased to  2\%. We highlight several trends that are somewhat different than those observed for smaller perturbations. First, the  images of the resulting x-points found by the surrogate, while still very close to those obtained by direct solution, are no longer identical. This is confirmed by Table \ref{Tab:xpt-radii2}, which shows qualitative agreement in the locations of the perturbed x-points produced by surrogates and direct solutions but also reveals more x-points computed by surrogates that are somewhat further from the reference solution than those found by direct solution. Similar trends are evident for both strike points and contact points. In particular, Figure \ref{fig:pt-distribs2} shows that the surrogate computations yield a small number (9 for sparse grid level 3, 12 for level 4) of false positives for strike points. The third row of images in Figure \ref{fig:pt-distribs2} shows examples of contact points found by each of the simulations. The direct solver and surrogates for sparse grid levels 3 and 4 all find the majority of contact points near $y=1$; the surrogates for these levels also find a small number (on the order of 10 out of 2000) in other places in the domain.

Note that all these results identify some effects that variability has on properties of the solution such as location of the x-point and strike points. Not surprisingly, increased noise in the currents lead to wider distributions of such quantities, and, for 2\% noise, increased likelihood of the occurrence of contact points. The surrogate computations largely reproduce the results of the direct solution albeit with some loss of accuracy occurring with increased perturbation size. None of the examples lead to conclusions that would be dramatically misleading; for example, a false-positive strike point would lead to a slightly over-cautious approach when shielding the reactor. 

\begin{figure}[htb]\centering
\begin{tabular}{cccc}
\multicolumn{2}{c}{\Large{$1 \%$ current variability}} & \multicolumn{2}{c}{\Large{$2 \%$ current variability}} \\ [2ex]
\includegraphics[width=0.2\linewidth]{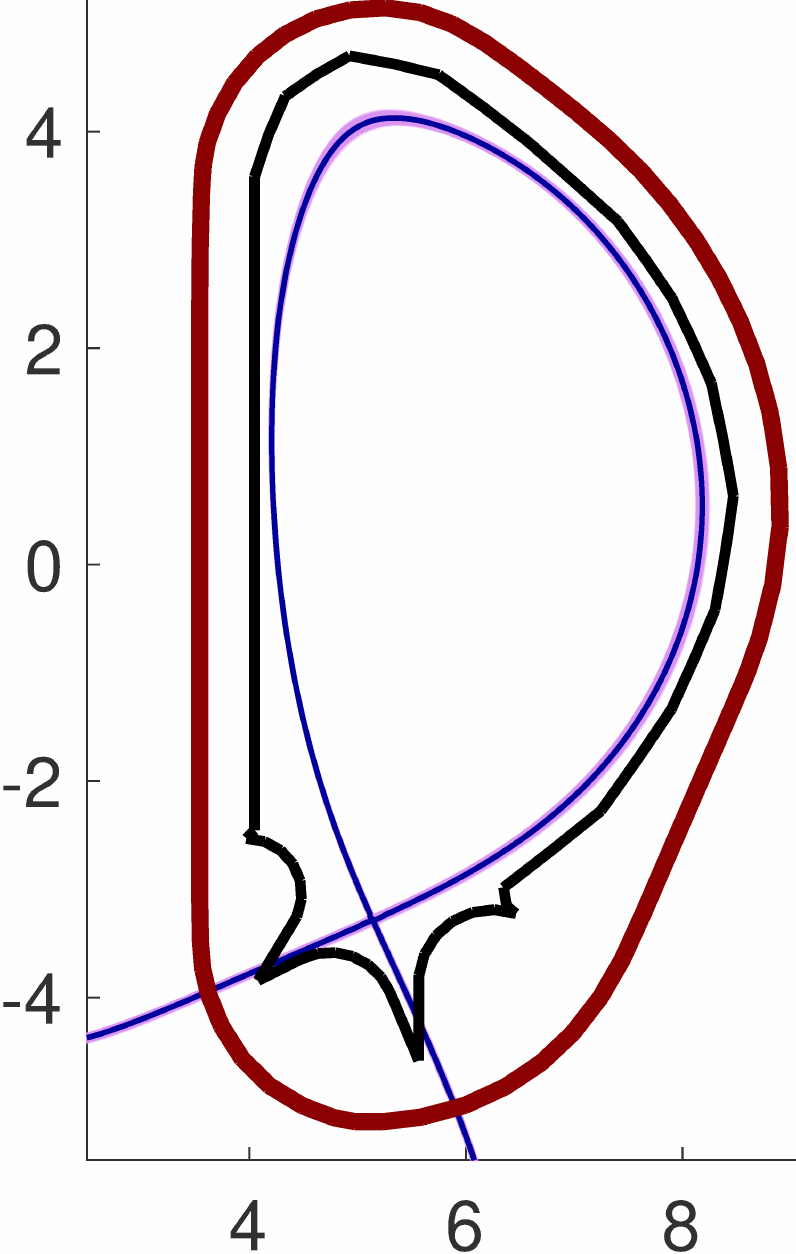} \quad &
\includegraphics[width=0.2\linewidth]{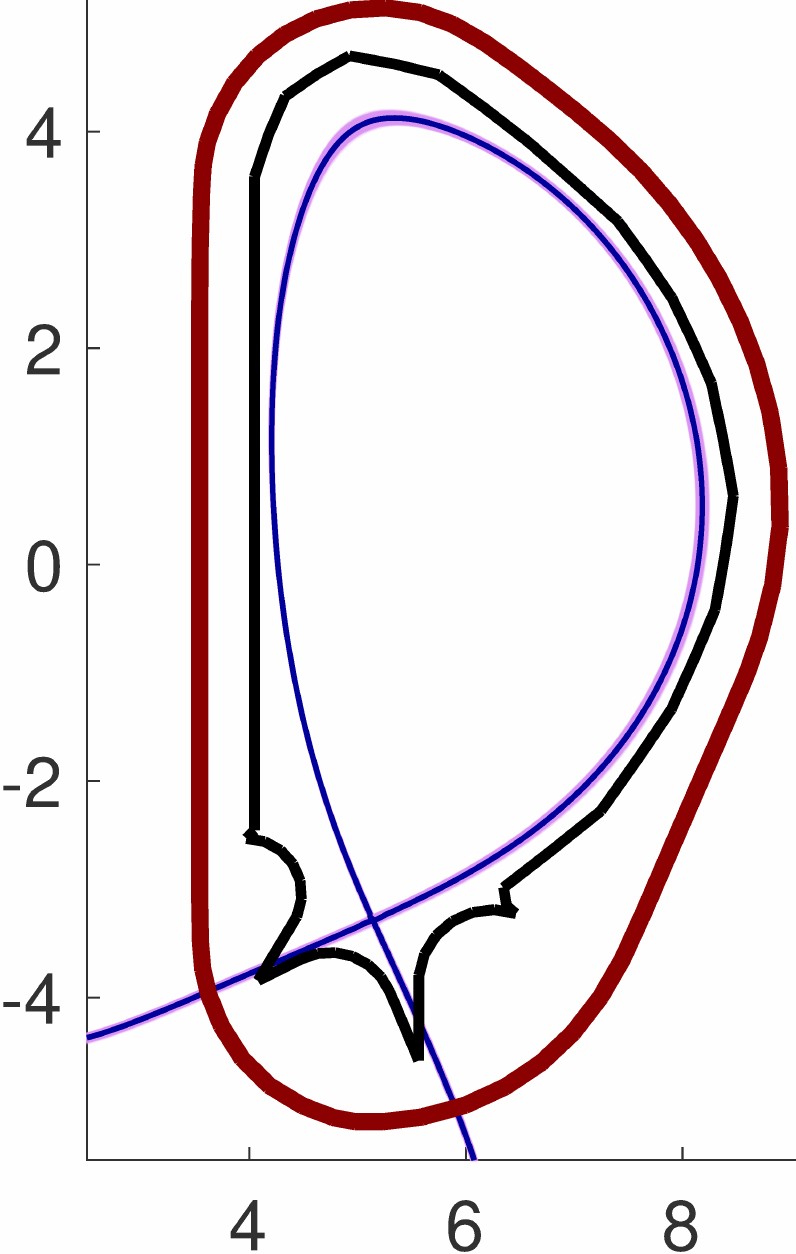} \quad &
\includegraphics[width=0.2\linewidth]{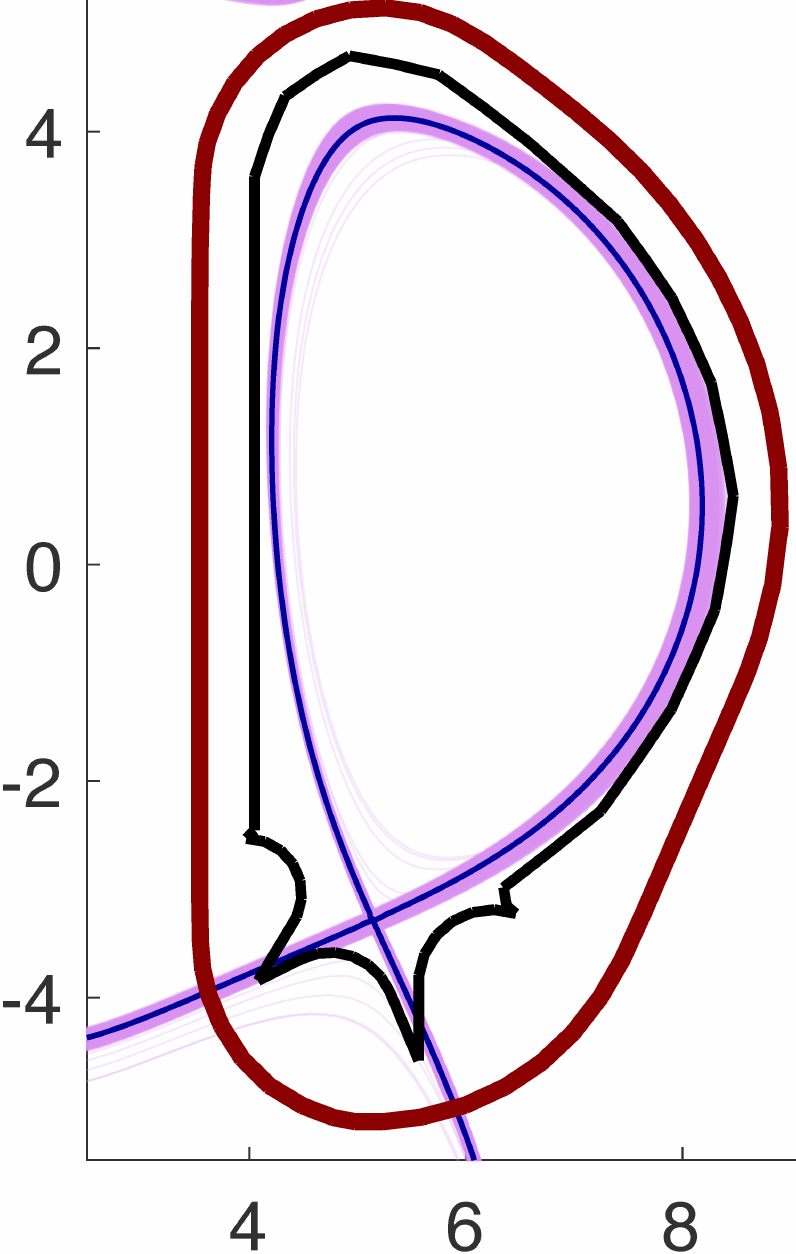} \quad &
\includegraphics[width=0.2\linewidth]{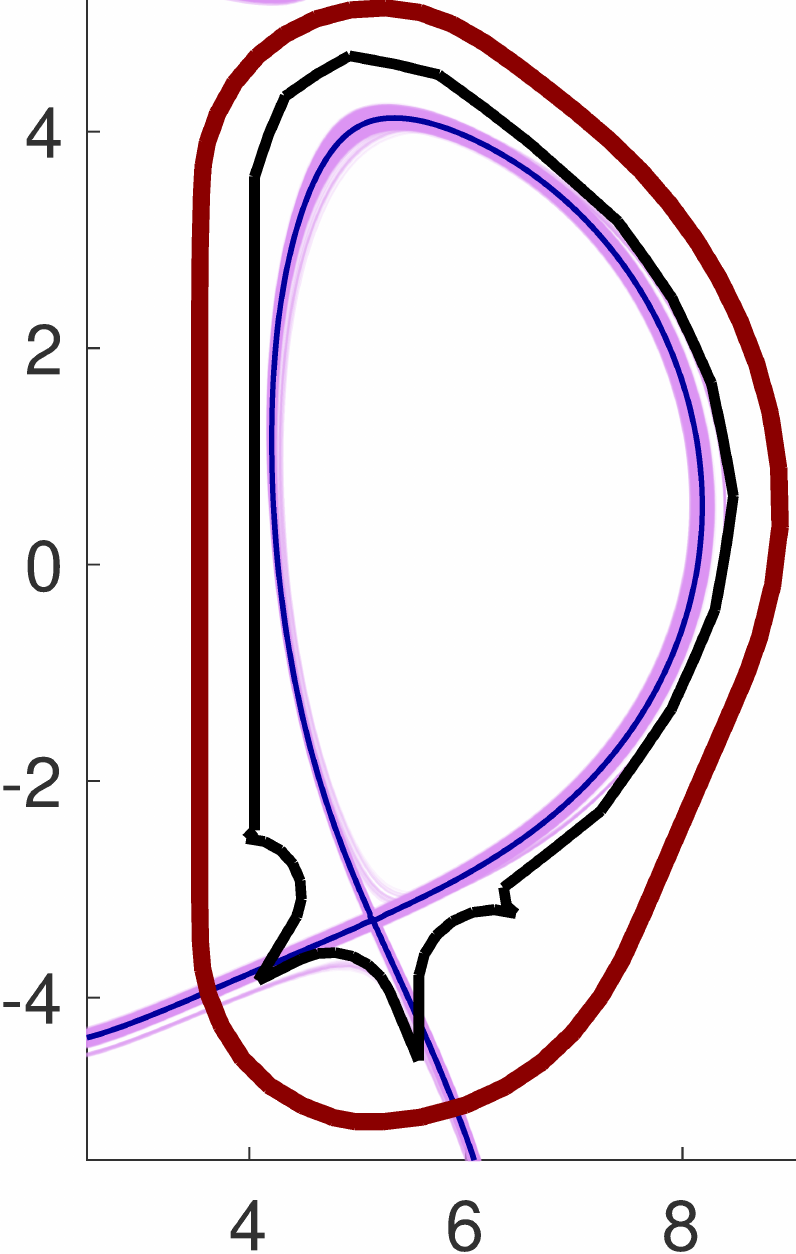} \\ [2ex]
\quad Surrogate level 3 & \quad Direct solver & \quad Surrogate level 3 & \quad Direct solver
\end{tabular}
\caption{Each of the panels displays the overlayed plasma boundaries resulting from 2000 random realizations of coil currents. The solid violet line is the plasma boundary resulting from the reference current values. The boundaries from the random realizations are displayed in light pink. The inner wall of the reactor is shown in solid black, while the outer wall is displayed in dark red. Left: Evaluations of the level 3 surrogate function with 1\% variability on the currents. Center left: Solutions obtained from the direct solver with 1\% variability on the currents. Center right: Evaluations of the level 3 surrogate function with 2\% variability on the currents. Right: Direct solver, 2\% current variability. In the figures, all scales are in meters. }
\label{fig:boundary-distributions}
\end{figure}

The variability of the currents also has an effect on the spread of locations for the plasma boundary. To show how the surrogate function captures this behavior accurately, 2000 equilibria were calculated using both direct solver and a surrogate built with 3 levels of the sparse grid. These computations were done using random perturbations of sizes 1\% and 2\% about the reference currents. The results of the experiment are shown on Figure \ref{fig:boundary-distributions}, where the plasma boundary from the reference current is plotted in dark violet, and the ones corresponding to the random perturbations are plotted in light pink and overlayed. In the figure it is evident that the qualitative behavior predicted by the surrogate function coincides with the one from the direct solution. The spread of the boundary focuses on the outer radius of the plasma (on the right side of the plot), increasing the likelihood for contacts on the outer side of the inner wall. This spread is more accentuated for surrogate evaluations, but located in the same regions as the direct solution, resulting in accurate---if slightly too cautions---predictions. This becomes even more evident for for 2\% noise, where it is also possible to observe isolated boundary curves resulting from plasma/wall contacts for both computational methods.

As the previous experiments have shown, the qualitative and quantitative agreement between direct computations and surrogate evaluations is significant. The key advantage of the use of these surrogate computations lies in their significantly smaller computational cost. These time savings come mostly from the fact that the surrogate bypasses the solution of the nonlinear system, but an additional source of savings comes from the fact that the surrogate can be designed to focus on smaller regions of interest. For the problem at hand, the direct solver is bound to provide an approximate solution over the entire computational domain, whereas for most applications, only the interior of the reactor is of interest. We have taken advantage of this fact by constructing a surrogate that is defined only inside of the reactor. This results in smaller vectors and shorter computation times.

We will demonstrate the time efficiency of the surrogate by comparing the performance of simulations done using the surrogate computations, for sparse grid levels 2 through 4, and the direct solver. In these tests, we used a dynamic approach for sampling based on the following ideas. The margin of error for the estimation of the mean value of a distribution based on a sample of size $n$ is given by $\epsilon=z^*\frac{\sigma}{\sqrt{n}}$, where $\sigma$ is the standard deviation of the distribution and $z^*$, known as the \textit{confidence coefficient}, is a constant depending on the desired confidence level for the estimate \cite{RoVa:1975}. For a confidence level of $95\%$, a commonly accepted standard, the value of the confidence coefficient is $z^*=1.96$. It follows that, with 95\% confidence, the margin of error of a sample of size $n$ can be estimated by
\begin{equation} \label{eq:RMSE-est}
\epsilon_n=1.96\frac{S_n}{\sqrt{n}}, 
\end{equation}
where the sample standard deviation $S_n = (n-1)^{-1}\sum_{i=1}^n {(x_i-\overline{x})^2}$ was used as an unbiased estimator for the unknown population standard deviation $\sigma$. The sample standard deviation can be computed dynamically using Welford's algorithm \cite{LF:1974,Welford:1962} as new samples are added to the data set. The process is as follows: the sample mean $m$ and the auxiliary variable $s$ are initialized as $m_0=0,s_0=0$. Then, as the $k$-th sample is drawn, these values are updated according to the formulas
\[ 
    m_k=m_{k-1}+\frac{x_k-m_{k-1}}{k},  \qquad
    s_k=s_{k-1} + (x_k-m_{k-1})(x_k - m_k).
\]
Finally, for a sample of size $k$, the sample standard deviation is given by $S_k=\sqrt{s_k/(k-1)}$. Combining this strategy with the approximation \eqref{eq:RMSE-est}, a stopping criterion for sampling is obtained. As new samples are drawn, the standard deviation is updated, $\epsilon_k$ is estimated, and the process is halted when $\epsilon_k$ falls below a prescribed threshold. In all our experiments, we will set $\epsilon=0.01$, thus requiring an estimate with a margin of error that, 95\% of the time, will be below 1\%.

We applied this strategy using the maximum entry-wise moving variance of both the solution $\psi_h$ obtained from the direct solver and the surrogate approximation $\widehat{\psi}_h$. The results are shown in Table \ref{Tab:timing}, which gives CPU times for current perturbations with 1\% and 2\% noise. We give several variants of these results.
In particular, the number of adaptively chosen samples varied somewhat with the choice of function used, $\psi_h$ or ${\widehat \psi}_h$.  For 1\% noise, the direct solver required more samples than the level 2 and 3 surrogates (see Table \ref{Tab: Table3}); for 2\% noise, it required fewer (see Table \ref{Tab: Table5}). It can be seen from these results that in all cases for level 2 and level 3 sparse grids, the cost of using the surrogate is significantly smaller than that of using the direct solver. In simulations with 1\% noise, the direct solver was approximately 33 times more expensive than the level 2 surrogate, and approximately 9 times more than the level 3 surrogate. With 2\% noise, the analogous factors were approximately 33 and 7. Here, we are considering the total time of the simulation with dynamic sampling which tends to benefit the direct solver---as it usually requires fewer samples to achieve a target RMS. As noted above, from the point of view of accuracy, there is essentially no difference between the surrogate and the direct solutions for 1\% noise; for 2\% noise, the significant cost reduction comes with some compromise in accuracy but little difference in conclusions concerning features of the plasma.
\begin{table}[ht]
	\centering
	\scalebox{0.7}{
		\begin{tabular}{|c|c|c|c|c|c|c|}
			\cline{3-6}
			\multicolumn{2}{c|}{}&{\ \ Surrogate\ \ } & \ \ Surrogate \ \ &\ \ Surrogate \ \ &{\ \ Direct\ \ }\\	
			\multicolumn{2}{c|}{}&\multicolumn{1}{c|}{Level 2} &\multicolumn{1}{c|}{Level 3} &\multicolumn{1}{c|}{Level 4}
					&\multicolumn{1}{c|}{ Solver }\\
			\hline
			\multirow{3}{*}{1\% Noise}
			& \# Samples &5100&5000&5000&5200\\
            \cline{2-6}
			&Time per Sample (s.) & 0.16 & 0.59 & 3.03 & 5.22\\
			\cline{2-6}
			&Total Time (s.) &831&2944&15141&27166\\
			\hline
			\hline
			
			\multirow{3}{*}{2\% Noise}
			& \# Samples &20000&25300&31100&19100\\
            \cline{2-6}
			&Time per Sample (s.) & 0.16 & 0.59 & 3.04 & 5.69\\
			\cline{2-6}
			&Total Time (s.) &3264& 14911&94602&108620\\
			\hline
	\end{tabular}}
		\caption{Average CPU times, in seconds, for evaluation of the surrogate function and direct solution of the free boundary problem.}
		\label{Tab:timing}
\end{table}
As mentioned above, the number of samples required to fall within a predetermined Monte Carlo tolerance varies depending on whether the direct solver or the surrogate are used to obtain an approximation of the stream function $\psi$. This points to a difference in the spread in the distributions generated by each of the methods. The means of the distributions, however, coincide to the prescribed tolerance by construction, regardless of whether the samples used are identical or not. For the experiment showcased in Table \ref{Tab: Table5}, two different sets of coil currents were generated considering a variability of 2\%. A surrogate of level 3 required 25300 randomly chosen current values for its Monte Carlo error to fall below 1\%, as given by \eqref{eq:RMSE-est}. These current values constitute sample A. When a direct solver was used instead, the dynamically estimated error needed only 19100 realizations to fall below the 1\% threshold. These currents constitute Sample B. As can be seen in Table \ref{Tab: Table5}, the variances associated with the quantities extracted from the surrogate are uniformly larger than those stemming from the direct solver. The mean values are nonetheless equal to within the prescribed tolerance. Despite the difference in required sample size, the level 3 surrogate was still faster than the direct solver by a factor of more than 7, as reported in Table \ref{Tab:timing}.

\begin{table}[ht]
	\centering
			\scalebox{0.62}{
		\begin{tabular}{c|c|c|c|c|c|c|c|c|}
			\cline{2-5}
				&\multicolumn{2}{c|}{Evaluation of surrogate}&\multicolumn{2}{c|}{Direct solver}\\
			&\multicolumn{2}{c|}{Sample A (25300 realizations)}&\multicolumn{2}{c|}{Sample B (19100 realizations)}\\
			\hline
			\multicolumn{1}{|c|}{Quantity of interest}&Mean&Variance&Mean&Variance\\
			\hline
			\multicolumn{1}{|c|}{x point}&(5.14,-3.29)&(2.89e-04,1.48e-03)&(5.14,-3.29)&(2.56e-04,1.31e-03)\\
			\hline
			\multicolumn{1}{|c|}{magnetic axis}&(6.41,0.61)&(1.53e-03,1.14e-03)&(6.41,0.61)&(1.19e-03,6.42e-04)\\
			\hline
			\multicolumn{1}{|c|}{strike} &(4.16,-3.71)&(2.42e-03,2.38e-03)&(4.16,-3.71)&(8.59e-04,2.09e-03)\\
			\multicolumn{1}{|c|}{points}&(5.56,-4.22)&(5.18e-08,3.20e-03)&(5.56,-4.21)&(5.06e-08,3.13e-03)\\
			\hline
			\multicolumn{1}{|c|}{inverse aspect ratio} &0.32&5.53e-06&0.32&4.85e-06\\
			\hline
			\multicolumn{1}{|c|}{elongation} &1.86&3.00e-04&1.86&3.12e-04\\
			\hline
			\multicolumn{1}{|c|}{upper triangularity} &0.43&2.87e-04&0.43&2.65e-04\\
			\hline
			\multicolumn{1}{|c|}{lower triangularity} &0.53&2.37e-04&0.53&2.15e-04\\
			\hline
	\end{tabular}}
	\caption{Comparison of sample means and variances for relevant physical quantities extracted from two samples by the surrogate function and the direct solver. The number of realizations for each of the samples was determined dynamically to achieve a Monte Carlo estimation accurate to at least 1\% with 95\% confidence. The noise level in the current values was of 2\% and the surrogate was built with a sparse grid of level 3.}
	\label{Tab: Table5}
\end{table}
In a similar vein, Table \ref{Tab: Table3} displays analogous results for an experiment where the coil currents were perturbed by 1\%. Samples for $\widehat\psi$  and $\psi_h$ were taken using either surrogates of levels 2 and 3 or the direct solver. The number of realizations was determined dynamically by setting the Monte Carlo tolerance to 1\% and using the estimator \eqref{eq:RMSE-est}. The mean values of all quantities coincide to within the prescribed tolerance, but the time required to generate the data using the surrogates is substantially shorter (see Table \ref{Tab:timing}).
\begin{table}[ht]
	\centering
			\scalebox{0.65}{
		\begin{tabular}{c|c|c|c|c|c|c|c|c|c|c|c|c|}
			\cline{2-7}
				&\multicolumn{2}{c|}{Evaluation of surrogate (Level 2)}&\multicolumn{2}{c|}{Evaluation of surrogate (Level 3)}&\multicolumn{2}{c|}{Direct solver}\\
			&\multicolumn{2}{c|}{5100 realizations}&\multicolumn{2}{c|}{5000 realizations}&\multicolumn{2}{c|}{5200 realizations}\\
			\hline
			\multicolumn{1}{|c|}{Quantity of interest }&Mean&Variance&Mean&Variance&Mean&Variance\\
			\hline
			\multicolumn{1}{|c|}{x point}&(5.14,-3.29)&(1.07e-04,3.58e-04)&(5.14,-3.29)&(1.06e-04,3.55e-04)&(5.14,-3.29)&(1.09e-04,3.57e-04)\\
			\hline
			\multicolumn{1}{|c|}{magnetic axis}&(6.41,0.61)&(2.92e-04,2.84e-04)&(6.41,0.61)&(2.84e-04,4.85e-04)&(6.41,0.61)&(2.98e-04,2.79e-04)\\
			\hline
			\multicolumn{1}{|c|}{strike} &(4.16,-3.71)&(1.98e-04,5.43e-04)&(4.16,-3.71)&(1.97e-04,5.40e-04)&(4.16,-3.71)&(2.02e-04,5.54e-04)\\
			\multicolumn{1}{|c|}{points}&(5.56,-4.22)&(1.28e-08,7.93e-04)&(5.56,-4.22)&(1.32e-08,8.14e-04)&(5.56,-4.22)&(1.27e-08,7.89e-04)\\
			\hline
			\multicolumn{1}{|c|}{inverse aspect ratio} &0.32&1.19e-06&0.32&1.17e-06&0.32&1.19e-06\\
			\hline
			\multicolumn{1}{|c|}{elongation} &1.86&3.39e-05&1.86&3.44e-05&1.86&3.47e-05\\
			\hline
			\multicolumn{1}{|c|}{upper triangularity} &0.43&6.41e-05&0.43&6.60e-05&0.43&6.58e-05\\
			\hline
			\multicolumn{1}{|c|}{lower triangularity} &0.53&3.71e-05&0.53&3.73e-05&0.53&3.73e-05\\
			\hline
	\end{tabular}}
	\caption{Sample means and variances for relevant physical parameters extracted using levels 2 and 3 of the surrogate and the direct solver. The quantities of interest are the coordinates of the x-point, magnetic axis and strike points, as well the values for plasma shaping parameters. The experiment was carried over with 1\% noise in the coil currents. The sample size was increased dynamically using batches of 100 samples until the estimated Monte Carlo estimation error fell below $0.01$. } 
	\label{Tab: Table3}
\end{table}

Note that a pre-processing step is required to build the surrogates; the required times are shown in Table \ref{Tab:timing-construction}. This ``offline'' step, while not negligible, is done only once and need not be repeated if multiple simulations are done of if the number of samples needs to be increased. Even if these one-time costs are factored in, the speedup for the level 3 surrogate is still a factor of 3 for 2\% perturbations.
 
\begin{table}[ht]
	\centering
	\scalebox{0.8}{
		\begin{tabular}{c|c|c|c|c|c|c|}
			\cline{2-4}	
			&\multicolumn{1}{c|}{Level 2} &\multicolumn{1}{c|}{Level 3} &\multicolumn{1}{c|}{Level 4} \\
			\hline
			\multicolumn{1}{|c|}{1\% Total Time (s.)} &2289&19375&126276\\
			\hline
			\multicolumn{1}{|c|}{2\% Total Time (s.)} &2378&20129&131194\\
			\hline
	\end{tabular}}
		\caption{Processing time, in seconds, required for the construction of surrogate functions of different levels.}
		\label{Tab:timing-construction}
\end{table}

\section{Concluding remarks}\label{sec:Conclusions}
We conclude with some comments highlighting the new features and properties of the methods introduced in this study. The main contribution is the development of an efficient computational algorithm for approximating the solution of the nonlinear system of equations arising from Grad-Shafranov free boundary problem. This model depends on parameters, the intensity of electric currents in coils that control the equilibrium configuration of the plasma in the model. Simulation of the influence of these parameters on the configuration requires multiple solutions for multiple realizations of parameters, and the use of surrogate approximations reduces the computational costs of simulation. For the surrogate computations, we used the sparse grid collocation method, where the approximation interpolates the solution at a distinguished set of points in the parameter domain, a so-called sparse grid. We implemented the surrogate computations using two existing codes,  {\tt FEEQS.M}, a {\tt MATLAB} implementation of CEDRES++ \cite{Heumann:feeqsm,CEDRES} for the Grad-Shafranov equation, and {\tt spinterp} \cite{spinterp-v5.1,KlBa:2005} for sparse grid collocation. Both codes were modified slightly, to facilitate adaptive grid refinement for the equations and to enhance vectorized performance in {\tt MATLAB} for collocation.

Our results demonstrate the enhanced efficiency of the surrogate computations, with significantly smaller costs (by factors ranging from 7 to 30) for constructing approximate solutions of good accuracy. We did observe, however, that with the number of parameters in the model (12), there was some reduction in accuracy for the surrogates, although the general trends revealed using them were the same as those seen with true solutions. As for all surrogate approximations, there is overhead associated with the offline computations needed to {\it construct} the surrogate; if this overhead is amortized over multiple simulations, its impact will be small.

\section{Acknowledgements}\label{sec:Acknowledgements}
All the free boundary computations were carried over using the code {\tt FEEQS.M} \cite{Heumann:feeqsm,CEDRES}. The authors are deeply grateful to Holger Heumann, the INRIA CASTOR team, and all the development team of the CEDRES++ free boundary solver for sharing access to the code and for helping us get up to speed with its usage. The authors would also want to thank Antoine Cerfon (NYU-Courant) for his valuable insights. 

Howard Elman has been partially supported by the U. S. Department of Energy under grant DE-SC0018149 and by the National Science Foundation under grant DMS1819115. Jiaxing Liang has been partially funded by the U. S. Department of Energy under grant DE‐SC0018149. Tonatiuh S\'anchez-Vizuet has been partially funded by the U. S. Department of Energy under Grant DE-FG02-86ER53233.

\bibliographystyle{abbrv}
\bibliography{references}

\end{document}